\documentclass[preprint,12pt]{elsarticle}

\usepackage{amssymb,latexsym,amsmath} %Standard packages
\usepackage{bm}
\usepackage{graphicx} 
\usepackage{caption}
\usepackage{subcaption}
\usepackage{siunitx}
\usepackage{textcomp}
\usepackage{color}
\usepackage{todonotes}
\usepackage{blindtext}
\usepackage{comment}
\DeclareMathOperator*{\argmax}{arg\,max}
\DeclareMathOperator*{\argmin}{arg\,min}

\journal{Journal of Computational Physics}

\begin{document}

% \tableofcontents

\begin{frontmatter}



%\title{Data assimilation by conditioning the Gaussian random field for stocahstic PDEs}

%\title{Machine learning for data assimilation, uncertainty quantification and reduction}
\title{Conditional Karhunen-Lo\`eve expansion for uncertainty quantification and active learning in partial differential equation models}


%\author{}

%\address{}

%\author{Ramakrishna Tipireddy\fnref{myfootnote}}
%\address{Pacific Northwest National Laboratory, P.O. Box 999, MSIN K7-90, Richland, WA 99352}

\author[mymainaddress]{Ramakrishna Tipireddy}
\ead{Ramakrishna.Tipireddy@pnnl.gov}

\author[mymainaddress]{David A Barajas-Solano}
\ead{David.Barajas-Solano@pnnl.gov}

\author[mymainaddress]{Alexandre Tartakovsky\corref{mycorrespondingauthor}}
\cortext[mycorrespondingauthor]{Corresponding author}
\ead{Alexandre.Tartakovsky@pnnl.gov}

\address[mymainaddress]{Pacific Northwest National Laboratory, P.O. Box 999, MSIN K7-90, Richland, WA 99352}

%\begin{itemize}
%	\item Abstract and key words
%	\item Introduction
%	\item Measurement of the permeability
%	\item Estimate parameters of covariance model of GRF
%	\item KL expansion and then condition the GRF 
%	\begin{itemize}
%		\item MC simulations
%		\item SG simulations with reduced dimension
%	\end{itemize}
%	\item condition the GRF and then KL expansion
%	\begin{itemize}
%		\item MC simulations
%		\item SG simulations with reduced dimension
%	\end{itemize}
%	\item Numerical examples
%
%	\item Discussion of numerical results
%	\item Conclusions
%\end{itemize}

\begin{abstract}
  We use a conditional Karhunen-Lo\`eve (KL) model to quantify and reduce uncertainty in a stochastic partial differential equation (SPDE) problem with partially-known space-dependent coefficient, $Y(x)$.
  We assume that a small number of $Y(x)$ measurements are available and model $Y(x)$ with a KL expansion. 
  We achieve reduction in uncertainty by conditioning the KL expansion coefficients on measurements. 
  We consider two approaches for conditioning the KL expansion: In Approach 1, we condition the KL model first and then truncate it.
  In Approach 2, we first truncate the KL expansion and then condition it.

  We employ the conditional KL expansion together with Monte Carlo and sparse grid collocation methods to compute the moments of the solution of the SPDE problem.
  Uncertainty of the problem is further reduced by adaptively selecting additional observation locations using two active learning methods.
  Method 1 minimizes the variance of the PDE coefficient, while Method 2 minimizes the variance of the solution of the PDE. 
 
  We demonstrate that conditioning leads to dimension reduction of the KL representation of $Y(x)$.
  For a linear diffusion SPDE with uncertain log-normal coefficient, we show that Approach 1 provides a more accurate approximation of the conditional log-normal coefficient and solution of the SPDE than Approach 2 for the same number of random dimensions in a conditional KL expansion. 
  Furthermore, Approach 2 provides a good estimate for the number of terms of the truncated KL expansion of the conditional field of Approach 1.
  Finally, we demonstrate that active learning based on Method 2 is more efficient for uncertainty reduction in the SPDE's states (i.e., it leads to a larger reduction of the variance) than active learning using Method 2.  
  
\end{abstract}


\end{frontmatter}


\section{Introduction}
\label{sec:intro}
% \todo[inline]{remove symbols and equations, add more literature, describe goals, contributions, outcomes }

Uncertainty quantification in partial differential equations (PDE) problems with partially known parameters (e.g., coefficients and source terms) is often performed by modeling these partially known quantities as random variables with appropriate probability distribution functions.
Spectral methods such as Polynomial chaos (PC)-based stochastic Galerkin~\cite{RK:Ghanem1991,Babuka2002} and stochastic collocation~\cite{RK:Xiu2002,Babuka2010} are commonly used for solving PDEs with random parameters. In spectral methods, random fields are represented in terms of their Karhunen-Lo\`eve (KL) expansions. While infinite KL expansions are necessary to exactly represent the two-point statistics of a random field, numerical treatment of PDEs requires truncating KL expansions.
These KL expansions are truncated based on the decay of their eigenvalues, so that the random fields can be reconstructed with sufficient accuracy using the retained KL terms.

The computational cost of spectral methods exponentially increases with the dimensionality of the stochastic problem (i.e., the number of terms in the truncated KL expansion)~\cite{Venturi2013JCP,Lin2010JSC,Lin2010JCP,Lin2009AWR,barajassolano-2016-stochastic}.
Various approaches have been recently proposed to address this issue: by finding the solution in reduced dimensional spaces using basis adaption~\cite{RK:Tipireddy2014}, domain decomposition methods~\cite{TIPIREDDY2017203, tipireddy2017stochastic}, sliced inverse regression~\cite{li2016inverse, yang2017sliced}, and sparsity enhancing together with the active subspace method~\cite{yang2016enhancing}, among others.
In this work, we propose reducing the computational cost of spectral methods by conditioning the KL expansion on available data.
We also demonstrate that conditioning on data reduces uncertainty of predictions, i.e., reduces the variance of quantities of interest of the governing stochastic PDEs.

Most of existing UQ methods, including all work referenced above, is based on stochastic models of unknown fields with constant variance and stationary covariance functions, that is, covariance functions of the form $C(x,y) = C(|x-y|)$.
The unconditional statistics of fields is estimated from data using the so-called ergodicity assumption, where fields are treated as realizations of a random process, and the spatial statistics of the fields is assumed to be the same as the ensemble statistics of the generating random process. However, there is no reason to assume that the variance at locations where measurements are available should be the same as the variance in locations with no measurements.

Gaussian process (GP) regression, also known as kriging, has been used in geostatistics and hydrogeological modeling since its introduction in the sixties~\cite{matheron1963principles} to represent partially observed properties of materials (e.g., permeability of geological porous media) as random field conditioned on observations~\cite{zhang-2002-stochastic,cressie-2015-geostatistics}.
A characteristic feature of GP models of random fields is that their \emph{conditional} statistics, i.e., statistics conditioned on data, are not stationary.
For example, the conditional variance of a GP model is a function of space (e.g., in the absence of measurement errors, it is zero at the observation locations), and the conditional covariance function $C^c(x,y)$ depends explicitly on both $x$ and $y$.
There are very few studies of PDEs with non-stationary random fields conditioned on data, including the conditional Moment Method~\cite{neuman1993prediction,morales2006non} limited to parameters with small variances, and a few papers on \emph{conditional} PC methods~\cite{LU2004859,ossiander2014,li2014conditional}. 

Two main approaches have been proposed for conditional spectral methods: (i) first truncating the KL expansion of the random parameters and then conditioning the resulting truncated expansion on data~\cite{ossiander2014}, and (ii) first conditioning the infinite KL expansion and then truncating it~\cite{li2014conditional}.
Here, we demonstrate that the truncating and conditioning operations do not commute as the two approaches produce different results.
A detailed analysis of spectral methods for PDEs with non-stationary random inputs is clearly lacking.
In this work, we study the application of conditional KL models to quantify and reduce uncertainty in physical systems modeled with stochastic PDEs (SPDEs). 
We compare the solution of the SPDE in terms of its conditional mean and variance obtained using both the conditioning first and truncating first approaches, and discuss the merits of both constructions. 
Our results show that the approach of conditioning first and then truncating (with $N$ terms) approach is more accurate than the approach of truncating first (with $N$ terms) and then conditioning.
We also show that for the truncating-first approach the final dimension of the conditional model is reduced from $N$ to $N-N_s$, where $N_s$ is the number of measurements. 

Furthermore, we adopt active learning~\cite{zhao2017active, raissi2017inferring} to identify additional observation locations in order to efficiently reduce predictive uncertainty in terms of the variance of the solution of the SPDE.
We consider two criteria for identifying observation locations.
In the first criterion, we choose the location that minimizes a norm of the variance of the conditional KL expansion of the model parameters.
In the second criterion, we propose a novel, GP regression-based approximation to the conditional variance of the solution of the SPDE, and choose locations that minimize this approximation.
We demonstrate that the second strategy leads to higher reduction of predictive uncertainty for the same number of additional measurements.

This manuscript is organized as follows: In Section \ref{sec:spde-gp} we formulate a steady-state stochastic diffusion equation with random coefficient and the GP model for the random coefficient.
In Section~\ref{sec:cond_rand}, we describe two approaches for constructing finite-dimensional conditional KL models of the random coefficient.
In Section~\ref{sec:active_learning}, we present active learning criteria for further reducing uncertainty in conditional KL models.
Numerical examples are given in Section~\ref{sec:numerical} and conclusions are presented in Section~\ref{sec:conclusions}.

%\section{Measurements of a permeability field}

\section{Governing Equations}
\label{sec:spde-gp}

We study a two-dimensional steady-state diffusion equation with a random diffusion coefficient $k(\mathbf{x}, \omega) : D \times \Omega \to \mathbb{R}$, $D \subset \mathbb{R}^2$ that is bounded and strictly positive:
\begin{equation}\label{RK:eq:rf_bound}
  0 < k_l \leq k(\mathbf{x},\omega) \leq k_u < \infty \quad \text{a.e.} \quad \text{in} \quad D \times \Omega.
\end{equation}
We seek the stochastic solution $u(\mathbf{x},\omega) : D \times \Omega \to \mathbb{R}$ to the problem
\begin{equation}
  \begin{aligned}
    \label{RK:eq:spde}
    \nabla \cdot \left [ k(\mathbf{x},\omega) \nabla u(\mathbf{x},\omega) \right ] &= 0 && \text{in}~D\times \Omega,\\
    u(\mathbf{x},\omega) &= u_{\Gamma} && \text{on}~\Gamma \subset \partial D,\\
    \nabla u(\mathbf{x},\omega) \cdot n &= u_{\Gamma'} && \text{on}~\Gamma' = \partial D \backslash \Gamma,
  \end{aligned}
\end{equation}
where the boundary conditions $u_{\Gamma}$ and $u_{\Gamma'}$ are deterministic and known, and $n$ denotes the outward-pointing unit vector normal to $\partial D$.
In many practical applications, the full probabilistic characterization of the coefficient $k(\mathbf{x},\omega)$ is not known, but measurements of $k$ are available at a few spatial locations.
In this work, we assume that the distribution of $k$ is known and is log-normal~\cite{Ghanem1999,zhang-2002-stochastic}, i.e., $k(\mathbf{x},\omega) \equiv \exp[g(\mathbf{x},\omega)]$, where $g(\mathbf{x},\omega) \equiv \log k(\mathbf{x}, \omega)$ is a Gaussian random field.
We construct the GP prior or \emph{unconditional} (i.e., not conditioned on measurements) model for $g$ employing the following two-step approach common in geostatistics~\cite{cressie-2015-geostatistics}: we first select a parameterized GP prior covariance kernel $C_g(\mathbf{x}, \mathbf{x}' \mid \bm{\theta}) : D \times D \to \mathbb{R}$, with hyperparameters $\bm{\theta}$; and next, we compute an estimate $\hat{\bm{\theta}}$ for the hyperparameters from available measurements of $\log k$ by type-II maximum likelihood estimation~\cite{cressie-2015-geostatistics,rasmussen2006gaussian}.
In the first step, we assume that $g$ is wide-sense stationary with zero mean and squared exponential covariance function
\begin{gather}
  \label{gp_prior}
  g \sim \mathcal{GP}(0, C_g(\mathbf{x}, \mathbf{x}' \mid \bm{\theta})),\\
  \label{cov_g}
  C_g(\mathbf{x}, \mathbf{x}' \mid \bm{\theta}) = \sigma^2 \exp \left [ -\frac{(x_1-x'_1)^2}{l_1^2} -\frac{(x_2-x'_2)^2}{l_2^2} \right ],
\end{gather}
where $\bm{\theta} = \{ \sigma, l_1, l_2 \}$ is the set of hyperparameters of the covariance kernel:
$\sigma$ is the standard deviation, and $l_1$ and $l_2$ are the correlation lengths along the $x_1$ and $x_2$ spatial coordinates, respectively.
To estimate $\bm{\theta}$, we assume that measurements are of the form $y_i = g(\mathbf{x}^{*}_i) + \epsilon_i$, $i \in [1, N_s]$, where $N_s$ is the number of observations, $\mathbf{x}^{*}_i$ is the measurement location for the $i$th measurement, and the $\epsilon_i \sim \mathcal{N}(0, \sigma_{\epsilon})$ are iid measurement errors with standard deviation $\sigma_{\epsilon}$ independent of $g(\mathbf{x}, \omega)$.
The observations and observation locations are arranged into the vector $\mathbf{y} = (y_1, \dots, y_{N_s})^{\top}$ and the matrix $X = (\mathbf{x}^{*}_1, \dots, \mathbf{x}^{*}_{N_s})$, respectively.
Similarly, the values of $g$ at the observation locations are arranged into the vector $\mathbf{g} = (g(\mathbf{x}^{*}_1), \dots, (\mathbf{x}^{*}_{N_s}))^{\top}$, so that $\mathbf{y}$ and $\mathbf{g}$ are related by
\begin{equation}
  \label{eq:g-obs}
  \mathbf{y} = \mathbf{g} + \bm{\epsilon}, \quad \bm{\epsilon} \sim \mathcal{N}(0, \sigma^2_{\epsilon} \mathbf{I}), \quad \mathbb{E}[ \mathbf{g} \bm{\epsilon}^{\top}] = 0.
\end{equation}

\begin{comment}
%
By the Bayes rule, the posterior of $\mathbf{g}$ given the measurements $(\mathbf{x}, \mathbf{y})$ and the hyperparameters $\bm{\theta}$ is given by
%
\begin{equation*}
  p(\mathbf{g} \mid X, \mathbf{y}, \bm{\theta}) = \frac{p( \mathbf{y} \mid \mathbf{g}, X, \bm{\theta}) p(\mathbf{g} \mid X, \bm{\theta})}{p(\mathbf{y} \mid X, \bm{\theta})},
\end{equation*}
%
where $p(\mathbf{y} | \mathbf{g}, X, \bm{\theta}) \equiv \mathcal{N}(\mathbf{y} \mid \mathbf{g}, \sigma^2_{\epsilon} I)$ is the likelihood of the observations given $\mathbf{g}$ and $\bm{\theta}$, and $p(\mathbf{g} \mid X, \bm{\theta})$ is the GP prior for $\mathbf{g}$.
By construction, we have $p(\mathbf{g} \mid \bm{\theta}) \equiv \mathcal{N}(\mathbf{g} \mid 0, C_s(\bm{\theta}))$, where $C_s(\bm{\theta})$ is the $N_s \times N_s$ matrix with entries $C_{s, ij} = C_g(x^{*}_i, x^{*}_j \mid \bm{\theta})$.
The marginal likelihood $p(\mathbf{y} \mid X, \bm{\theta})$ is the normalizing constant in the Bayes rule, given by
%
\begin{equation}
  \label{eq:ml}
  p(\mathbf{y} \mid X, \bm{\theta}) = \int p(\mathbf{y} \mid \mathbf{g}, X, \bm{\theta}) p(\mathbf{g} \mid X, \bm{\theta}) \, \mathrm{d} \mathbf{g}.
\end{equation}
\end{comment}

As stated above, in this work we employ type II maximum likelihood estimation \cite{cressie-2015-geostatistics,rasmussen2006gaussian} to estimate the hyperparameters of the GP prior and of the likelihood of the observations:
\begin{equation}
  \label{eq:typeII-ML}
  (\bm{\theta}, \sigma_{\epsilon}) \equiv \argmax_{\bm{\theta}, \sigma_{\epsilon}} L(\bm{\theta}, \sigma_{\epsilon}, X, \mathbf{y}),
\end{equation}
where $L(\bm{\theta}, \sigma_{\epsilon}, X, \mathbf{y})$ is the log-marginal likelihood function
\begin{equation}
  \label{eq:LML}
  L(\bm{\theta}, \sigma_{\epsilon}, X, \mathbf{y}) = - \frac{1}{2} y^{\top} | C_s(\bm{\theta}) + \sigma^2_{\epsilon} I |^{-1} y - \frac{1}{2} \ln |C_s(\bm{\theta}) + \sigma^2_{\epsilon} I| - \frac{N_s}{2} \ln 2 \pi,
\end{equation}
and $C_s(\bm{\theta})$ is the observation covariance matrix with $ij$th entry given by $C_g(\mathbf{x}^{*}_i, \mathbf{x}^{*}_j \mid \bm{\theta})$.

To compute the solution of the SPDE~\eqref{RK:eq:spde} conditioned on the observations $\mathbf{y}$, we adopt the stochastic collocation approach, which requires a finite-dimension stochastic representation of the random field $g(\mathbf{x}, \omega)$.
For constructing this representation, we will consider two strategies that rely on GP regression and the KL expansion of random fields:
\begin{enumerate}
\item In the first approach~\cite{li2014conditional}, the conditional random field $g^c(\mathbf{x},\omega)$ is obtained using GP regression and then discretized by calculating its KL expansion in terms of standard Gaussian random variables.
  Then, the KL expansion is truncated to an appropriate number of random dimensions $d^c$.
\item In the second approach~\cite{ossiander2014}, the unconditioned random field $g(\mathbf{x},\omega)$ is first discretized in terms of its KL expansion and unconditioned standard Gaussian random variables $\bm{\xi}$.
  Then, conditional Gaussian random variables $\tilde{\bm{\xi}}$, conditioned on the observations $\mathbf{y}$, are obtained by projection.
\end{enumerate}
We describe these two approaches in details in Section~\ref{sec:cond_rand}.

\section{Conditional KL models}
\label{sec:cond_rand}

\subsection{Approach 1: truncated KL expansion of the conditioned GP field}
\label{sec:trunc-cond-kl}

In the first approach, the conditional random field $g^c(\mathbf{x},\omega)$ is approximated with a KL expansion written in terms of standard Gaussian random variables.
We assume that the hyperparameters of the prior covariance function $C_g(\mathbf{x}, \mathbf{x}')$ have been estimated and are thus dropped from the notation.
The mean and covariance of conditioned Gaussian random field $g^c(\mathbf{x},\omega) \sim \mathcal{GP} \left (\mu_g^c, C_g^c(\mathbf{x}, \mathbf{x}') \right )$ are computed using GP regression~\cite{cressie-2015-geostatistics,rasmussen2006gaussian}:
\begin{align} 
  \label{eq:cond_mean}
  \mu_g^c(\mathbf{x}) &= C_g (\mathbf{x}, X) \left [ C_s + \sigma^2_{\epsilon} \right ]^{-1} \mathbf{y},\\
  \label{eq:cond_cov}
  C_g^c(\mathbf{x}, \mathbf{x}') &= C_g(\mathbf{x}, \mathbf{x}') - C_g(\mathbf{x}, X) \left [ C_s + \sigma^2_{\epsilon} \right ]^{-1} C_g(\mathbf{x}, \mathbf{x}').
\end{align}

The conditional field is then expanded using a truncated KL expansion~\cite{RK:Ghanem1991} as
\begin{equation} 
  \label{eq:cond_then_kl}
  g^c(\mathbf{x},\omega) = \mu^c_g(\mathbf{x}) + \sum_{i=1}^{d^c} \sqrt{\lambda_i^c} \phi^c_i(\mathbf{x}) \xi_i,  
\end{equation}
where $\{ \xi_i\sim \mathcal{N}(0,1) \}^{d^c}_{i = 1}$ are standard Gaussian random variables, and $\{ \lambda_i^c,\phi_i^c \}^{d^c}_{i = 1}$ are the first $d^c$ pairs of eigenvalues and eigenfunctions stemming from the eigenvalue problem
\begin{equation*}
  \int_D C_g^c(\mathbf{x}, \mathbf{x}') \phi^c(\mathbf{x}') \, \mathrm{d} \mathbf{x}' = \lambda^c \phi^c(\mathbf{x}),
\end{equation*}
where $C^c_g$ is given by~\eqref{eq:cond_cov}.
The conditional solution of the SPDE, $u^c$, can then be computed using MC or sparse grid collocation methods by sampling $g^c(\mathbf{x},\omega)$ using~\eqref{eq:cond_then_kl}, or by the PC method by constructing a spectral approximation of $u^c$ in terms of the $\xi_i$, $i \in [1, d^c]$.

\subsection{Approach 2. Conditioning truncated KL expansion of the unconditioned field}
\label{sec:cond-trunc-kl}

In this approach, introduced in~\cite{ossiander2014}, the KL expansion of $g(\mathbf{x}, \omega)$ is first truncated, and the resulting set of random variables are conditioned on the observations $\mathbf{y}$.
Here, we demonstrate that this approach reduces the number of random dimensions of the representation of $g^c$ by the number of observations and propose a method for rewriting the representation of $g^c$ in terms of the reduced number of random variables.
Effectively, this reduces dimensionality of the SPDE solution.
   
To present our approach we first summarize the conditional KL construction presented in~\cite{ossiander2014} and next describe our proposed reduced-dimension conditional KL representation.
We start with the KL expansion of the unconditional Gaussian random field $g(\mathbf{x},\omega)$:
\begin{equation}
  \label{eq:g-KL}
  g(\mathbf{x},\omega) = \sum_{i=1}^{\infty} \sqrt{\lambda_i} \phi_i(\mathbf{x}) \xi_i(\omega), 
\end{equation}
where the $\{\xi_i \sim \mathcal{N}(0, 1)\}^{\infty}_{i = 1}$ are iid standard Gaussian random variables.
In matrix notation,
\begin{equation}
  \label{eq:g-KL-matrix}
  g(\mathbf{x},\omega) = \Phi^{\top}(\mathbf{x}) \Lambda^{1 / 2} \bm{\xi}(\omega), 
\end{equation}
where $\Phi(\mathbf{x}) = (\phi_1(\mathbf{x}), \phi_2(\mathbf{x}), \dots )^{\top}$, $\Lambda = \operatorname{diag}(\lambda_1, \lambda_2, \dots)$, and $\bm{\xi}(\omega)$ is the infinite-dimensional random vector $\bm{\xi}(\omega) = (\xi_1(\omega), \xi_2(\omega), \dots)^{\top}$.
By Mercer's theorem, the covariance function of $g$ can be represented as the infinite sum
\begin{equation}
  \label{eq:cov-g-mercer}
  C_g(\mathbf{x}, \mathbf{x}') = \sum_{i=1}^{\infty} \lambda_i \phi_i(\mathbf{x}) \phi_i(\mathbf{x}'),
\end{equation}
where the eigenpairs $\{ \lambda_i, \phi_i \}^{\infty}_{i = 1}$ are the solutions to the eigenvalue problem
\begin{equation*}
  \int_D C_g(\mathbf{x}, \mathbf{x}') \phi(\mathbf{x}') \, \mathrm{d} \mathbf{x}' = \lambda \phi(\mathbf{x}).
\end{equation*}
In matrix notation, the Mercer expansion~\eqref{eq:cov-g-mercer} can be written as
\begin{equation}
  \label{eq:cov-g-mercer-matrix}
  C_g(\mathbf{x}, \mathbf{x}') = \Phi^{\top}(\mathbf{x}) \Lambda \Phi(\mathbf{x}).
\end{equation}

In \cite{ossiander2014}, the KL expansion~\eqref{eq:g-KL} is conditioned on the observations $\mathbf{y}$ by conditioning $\bm{\xi}$ on these observations as described in the following Eqs \eqref{eq:y-KL}--\eqref{M_matrix}.
Evaluating~\eqref{eq:g-KL-matrix} at $X$ and substituting into~\eqref{eq:g-obs} yields
\begin{equation}
  \label{eq:y-KL}
  \mathbf{y} = \Phi^{\top}(X) \Lambda^{1 / 2} \bm{\xi} + \bm{\epsilon}.
\end{equation}
Therefore, the joint distribution of $\mathbf{y}$ and $\bm{\xi}$ is given by
\begin{equation}
  \label{eq:y-xi-joint}
  \begin{bmatrix}
    \mathbf{y} \\ \bm{\xi}
  \end{bmatrix} \sim \mathcal{N} \left (
    \begin{bmatrix}
      0 \\ 0
    \end{bmatrix},
    \begin{bmatrix}
      C_s + \sigma^2_{\epsilon} I & \Phi^{\top}(X) \Lambda^{1 / 2} \\
      \Lambda^{1 / 2} \Phi^{\top}(X) & I
    \end{bmatrix}
  \right ),
\end{equation}
where we have employed the relation $\mathbb{E}[ \mathbf{g} \bm{\epsilon}^{\top}] = 0$ and introduced the notation $C_s \equiv \mathbb{E}[ \mathbf{g} \mathbf{g}^{\top}]$.

It follows that the distribution of $\bm{\xi}$ conditioned on the measurements is $\bm{\xi} \mid X, \mathbf{y} \sim \mathcal{N}(\bm{\mu}, M)$, where
\begin{align*}
  \bm{\mu} &= \Lambda^{1 / 2} \Phi(X) \mathbf{y},\\
  M &= I - \Lambda^{1 / 2} \Phi(X) (C_s + \sigma^2_{\epsilon} I)^{-1} \Phi^{\top}(X) \Lambda^{1 / 2}.
\end{align*}

For simplicity, we denote the GP conditional on $(X, \mathbf{y})$ by $\tilde{g}$, and the conditional random vector by $\tilde{\bm{\xi}} = (\tilde{\xi}_1, \tilde{\xi}_2, \dots)^{\top}$.
The conditional GP then reads
\begin{equation}
  \label{sec:g-KL-cond}
  \tilde{g}(\mathbf{x},\omega) = \Phi^{\top}(\mathbf{x}) \Lambda^{1 / 2} \tilde{\bm{\xi}}(\omega).
\end{equation}

We now apply the process of conditioning random variables to the truncated KL expansion of the unconditional random field.
The KL expansion of the unconditional random field $g(\mathbf{x},\omega)$, truncated to $d$ terms, reads
\begin{align}
  \label{eq:trunc-KL-g}
  g^d(\mathbf{x}, \omega) &= \sum_{i = 1}^{d} \sqrt{\lambda_i} \phi_i(\mathbf{x}) \xi_i(\omega), \\
                 &= (\Phi^d)^{\top}(\mathbf{x}) (\Lambda^d)^{1 / 2} \bm{\xi}^d(\omega),
\end{align}
where $\Phi^d(\mathbf{x}) = (\phi_1(\mathbf{x}), \dots, \phi_d(\mathbf{x}))^{\top}$, $\Lambda^d = \operatorname{diag}(\lambda_1, \dots, \lambda_d)$, and $\bm{\xi}^d(\omega) = (\xi_1(\omega), \dots, \xi_d(\omega))^{\top}$.
This expansion corresponds to the truncated covariance
\begin{equation}
  \label{eq:cov-g-mercer-trunc}
  C^d_g(\mathbf{x}, \mathbf{x}') = (\Phi^d)^{\top}(\mathbf{x}) \Lambda^d \Phi^d(\mathbf{x}'),
\end{equation}
obtained by substituting $\Phi^d(\mathbf{x})$ and $\Lambda^d$ for $\Phi$ and $\Lambda$ in Eq \eqref{eq:cov-g-mercer}.
The truncated representation~\eqref{eq:trunc-KL-g} can be conditioned on the data $(X, \mathbf{y})$ by following the procedure outlined above, resulting in the conditional model
\begin{equation}
  \label{eq:trunc-KL-g-cond}
  \tilde{g}^d(\mathbf{x}, \omega) = (\Phi^d)^{\top}(\mathbf{x}) (\Lambda^d)^{1 / 2} \tilde{\bm{\xi}}^d,
\end{equation}
where $\tilde{\bm{\xi}}^d \equiv \bm{\xi}^d \mid X, \mathbf{y} \sim \mathcal{N}(\bm{\mu}^d, M^d)$, with $\bm{\mu}^d$ and $M^d$ given by
\begin{align}
  \label{M_matrix}
  \bm{\mu}^d &= (\Lambda^d)^{1 / 2} \Phi^d(X) \mathbf{y},\\
  M^d &= I - (\Lambda^d)^{1 / 2} \Phi^d(X) (C^d_s + \sigma^2_{\epsilon} I)^{-1} (\Phi^d)^{\top}(X) (\Lambda^d)^{1 / 2},
\end{align}
and where $C^d_s \equiv (\Phi^d)^{\top}(X) \Lambda^d \Phi^d(X)$ is the truncated measurement covariance matrix.

Note that both \eqref{sec:g-KL-cond} and~\eqref{eq:trunc-KL-g-cond} employ the same set of eigenpairs derived from the unconditioned covariance function $C_g$.
Nevertheless, $\tilde{\bm{\xi}}^d$ and $(\tilde{\xi}_1, \dots, \tilde{\xi}_d)$ in~\eqref{sec:g-KL-cond} have different joint distribution, which implies that the conditional truncated model \eqref{eq:trunc-KL-g-cond} is different than the model~\eqref{sec:g-KL-cond} after truncating to $d$ terms.

Due to the conditioning on $N_s$ measurements, the rank of $M ^d$ is
\begin{equation}
  \label{eq:r-dim}
\operatorname{rk} (M ^d) = r \equiv d - N_s, 
\end{equation}
so that the model~\eqref{eq:trunc-KL-g-cond} is effectively of dimension $r$.
In other words, conditioning the truncated model results in dimension reduction of the GP model for $g$, and reduction of the stochastic dimensionality of the SPDE problem.

To leverage this dimension reduction, we propose rewriting the model~\eqref{eq:trunc-KL-g-cond} in terms of $r$ random variables.
We write the eigendecomposition of $M^d$ as $M^d = Q D Q^{-1}$, where $D$ is the diagonal matrix of the form
\begin{equation}
  \label{eq:D-block-form}
  D = \begin{bmatrix}
    D^r& 0\\
    0 & 0
  \end{bmatrix},
\end{equation}
since $\operatorname{rk}(M^d) = r$.
Substituting the eigendecomposition of $M^d$ into~\eqref{eq:trunc-KL-g-cond}, we obtain
\begin{equation}
  \label{eq:trunc-KL-g-cond-r}
  \tilde{g}^d(\mathbf{x}, \omega) = \hat{g}^d(\mathbf{x}) + (\Phi^d)^{\top}(\mathbf{x}) (\Lambda^d)^{1 / 2} Q D Q^{-1} \bm{\zeta}^d(\omega),
\end{equation}
where $\hat{g}^d(\mathbf{x}) \equiv (\Phi^d)^{\top}(\mathbf{x}) (\Lambda^d)^{1 / 2} \bm{\mu}^d$ and $\bm{\zeta}^d \sim \mathcal{N}(0, I_d)$.

Let $\bm{\eta} = Q^{-1}\bm{\zeta}^d$; then, substituting $\bm{\eta}$ into \eqref{eq:trunc-KL-g-cond-r}, we obtain
\begin{multline}
  \tilde{g}^d(\mathbf{x},\omega)- \hat{g}^d(\mathbf{x}) = \\
  \begin{bmatrix}
    (\Phi^d_r)^{\top}(\mathbf{x}) & (\Phi^d)^{\top}_{r'}(\mathbf{x})
  \end{bmatrix}
  \begin{bmatrix}
    \Lambda^d_r & 0\\
    0 & \Lambda^d_{r'}
  \end{bmatrix}^{1 / 2}
  \begin{bmatrix}
    Q_r & Q_{rr'}\\
    Q_{r'r} & Q_{r'}
  \end{bmatrix}
  \begin{bmatrix}
    D_r\bm{\eta}_r \\
    0
  \end{bmatrix}.
\end{multline}
Similarly, let $\tilde{\bm{\eta}}_r = D_r \bm{\eta}_r$; then, $\tilde{g}^d$ can be expanded in terms of the reduced set of random variables $\tilde{\bm{\eta}}_r$ as
\begin{align*}	
  \tilde{g}^d(\mathbf{x},\omega) - \hat{g}(\mathbf{x}) &= \begin{bmatrix} (\Phi^d_r)^{\top}(\mathbf{x}) & (\Phi^d)^{\top}_{r'}(\mathbf{x})
  \end{bmatrix} \begin{bmatrix}
    \Lambda^d_r & 0\\
    0 & \Lambda^d_{r'}
  \end{bmatrix}^{1 / 2} \begin{bmatrix}
    Q_r & Q_{rr'}\\
    Q_{r'r} & Q_{r'}
  \end{bmatrix}  \begin{bmatrix}
    \tilde{\bm{\eta}}_r \\
    0  
  \end{bmatrix}\\
&= [(\Phi^d_r)^{\top}(\mathbf{x}) (\Lambda^d_r)^{1 / 2} Q_r + (\Phi^d_{r'})^{\top}(\mathbf{x}) (\Lambda^d_{r'})^{1 / 2} Q_{r'r}] \tilde{\bm{\eta}}_r. 
\end{align*}
Here, the components of the random vector $\tilde{\bm{\eta}}_r$ are correlated with one another; nevertheless, $\tilde{\bm{\eta}}_r$ can be converted to a set of uncorrelated random variables using an orthogonalization method such as Gram-Schmidt process, so that $\bm{\zeta}^r = A \tilde{\bm{\eta}}_r$ and $\tilde{\bm{\eta}}_r = A^{-1} \bm{\zeta}^r,$ where $\bm{\zeta}^r \sim \mathcal{N}(0, I_r)$.
Now, $\tilde{g}^d$ can be expanded in terms of the new set of random variables $\bm{\zeta}^r$ as
\begin{equation}
  \label{eq:trunc-KL-g-cond-Psi-r}
  \tilde{g}^d(\mathbf{x},\omega) = \hat{g}^d(\mathbf{x}) + (\tilde{\Psi}^r)^{\top}(\mathbf{x}) \bm{\zeta}^r,
\end{equation}
where
\begin{equation}
  \label{eq:Psi-r}
  (\tilde{\Psi}^r)^{\top}(\mathbf{x}) = (\Phi^d_r)^{\top}(\mathbf{x}) (\Lambda^d_r)^{1 / 2} Q_r A^{-1} + (\Phi^d_{r'})^{\top}(\mathbf{x}) (\Lambda^d_{r'})^{1 / 2} Q_{r'r} A^{-1}.
\end{equation}
The $r$ components of $\bm{\zeta}^r$ are uncorrelated Gaussian random variables and hence are also independent. Therefore, the conditional KL expansion~\eqref{eq:trunc-KL-g-cond-Psi-r} can be used in combination with standard PC-based or stochastic collocation-based methods for solving stochastic PDEs.

\section{Active learning for uncertainty reduction}
\label{sec:active_learning}

The conditioning of $g(\mathbf{x},\omega)$ KL expansion presented in Section~\ref{sec:cond_rand} leads to reduced uncertainty in $u(\mathbf{x},\omega)$.
This uncertainty can be further reduced by collecting additional measurements of $g$.
Giving limited ability to collect additional $g$ measurements, it is important to identify locations for additional measurements that optimally minimize the uncertainty in $u$.
This process of optimal data acquisition is referred to as \emph{active learning}~\cite{cohn-1996-active,zhao2017active,raissi2017inferring}.

In this section, we discuss the standard active learning method (Method 1) based on minimizing the variance of $g^c$ ~\cite{zhao2017active,raissi2017inferring} and propose an alternative approach (Method 2) based on minimizing the variance of $u^c$.
In Method 1, the variance of $u^c$ is reduced by minimizing the variance of $g^c$, but, as we show below, it does not lead to the maximum reduction of the $u^c$ variance.

It is worth noting that active learning also reduces the cost of computing the conditional solution $u^c$, as solving SPDE problems with coefficients having smaller variance requires smaller number of MC simulations, lower order of polynomial chaos in the stochastic Galerkin method, and lower sparse grid level in the stochastic collocation method.

\subsection{Method 1: minimization of the conditional variance of $g^c$}
\label{sec:active-learning-g}

In the standard active learning method~\cite{zhao2017active,raissi2017inferring}, a new location $\mathbf{x}^{*}$ for sampling $g$ is selected to minimize the variance of the conditional model $g^c$ conditioned on the full set of observations, including the new observation.
This can be done in closed form as, per Eq.~\eqref{eq:cond_cov}, the conditional variance $C^c_g$ depends only on the observation locations and not the observed values.
Therefore, the new observation location $\mathbf{x}^{*}$ is selected following the acquisition policy
\begin{equation}
  \label{eq:active-learning-criteria-g}
  \mathbf{x}^* \equiv \argmin_{\mathbf{x}'} \int_D C^c_g \left ( \mathbf{x} , \mathbf{x}  \mid [X, \mathbf{x}'] \right ) \, \mathrm{d} \mathbf{x},
\end{equation}
where $C^c_g(\cdot, \cdot \mid [X, \mathbf{x}'])$ denotes the covariance of $g$ conditioned on the original set of observation locations, $X$, and a new location $\mathbf{x}'$, and is computed using~\eqref{eq:cond_cov}.
In practice, the minimization problem~\eqref{eq:active-learning-criteria-g} is approximately solved by choosing $\mathbf{x}^*$ as $\argmax_{\mathbf{x}} C^c_g \left (\mathbf{x}, \mathbf{x} \mid X \right)$~\cite{zhao2017active,raissi2017inferring}.

Reducing the variance of $g^c$ also reduces the variance of $u^c$.
Nevertheless, there is no reason to assume that the observation locations provided by the acquisition strategy~\eqref{eq:active-learning-criteria-g} will lead to optimal variance reduction for $u^c$.

\subsection{Method 2: By minimizing the conditional covariance of $u^c$}
\label{sec:active-learning-u}

In Method 2, we chose location for measurements that minimize the conditional variance of $u$, that is,
\begin{equation}
  \label{eq:active-learning-criteria-u}
  \mathbf{x}^{*} \equiv \argmin_{\mathbf{x}'} \int_D {C}^c_u(\mathbf{x}, \mathbf{x}  \mid [X, \mathbf{x}']) \, \mathrm{d} \mathbf{x},
\end{equation}
where $C^c_u(\cdot, \cdot \mid [X, \mathbf{x}'])$ denotes the covariance of $u$ conditioned on the original set of observation locations, $X$, and a new location $\mathbf{x}'$.
Note that solving this minimization problem is significantly more challenging than solving~\eqref{eq:active-learning-criteria-g}, as it requires (i) constructing the conditional model $g^c$ conditioned on the new observation, and (ii) forward uncertainty propagation of $g^c$ through the SPDE problem to estimate $C^c_u$.
Here, step (i) is the most challenging as it requires knowledge of $g(\mathbf{x}')$ in order to construct the conditioned model $g^c$ ~\cite{cohn-1996-active}. Note that  $g(\mathbf{x})$ is only partially known; otherwise, this problem would not be uncertain.

As an alternative, we propose modeling $g^c$ and $u^c$ as components of a bivariate Gaussian field $h^c$ in order to derive an approximation to $C^c_u(\mathbf{x}, \mathbf{x} \mid [X, \mathbf{x}'])$ that does not require knowledge of $g(\mathbf{x}')$.
Specifically, let $g^c$ and $u^c$ denote the fields conditioned on the original data set $(\mathbf{x}, \mathbf{y})$, and let $h^c(\mathbf{x}, \omega) = [g^c(\mathbf{x}, \omega), u^c(\mathbf{x}, \omega)]^{\top}$ be a bivariate Gaussian field.
It is possible to compute the conditional covariance of $h^c$ conditioned on an additional measurement location $x'$ by employing~\eqref{eq:cond_cov}.
After marginalizing the $u^c$ component, we obtain the approximation
\begin{multline}
  \label{eq:GP_cov_u}
  C^c_u(\mathbf{x}, \mathbf{x} \mid [X, \mathbf{x}']) \approx C^c_u(\mathbf{x}, \mathbf{x} \mid X )\\
  - C^c_{ug}(\mathbf{x}, \mathbf{x}' \mid X ) [ C^c_g(\mathbf{x}', \mathbf{x}' \mid X) ]^{-1} C^c_{gu}(\mathbf{x}', \mathbf{x} \mid X),
\end{multline}
where $C^c_{ug}(\mathbf{x}, \mathbf{x}' \mid X)$ denotes the $u$--$g$ cross-covariance conditioned on $X$.
Note that this approximation only requires the observation location $x'$ and not the observation value $g(\mathbf{x}')$.

In order to apply the approximation~\eqref{eq:GP_cov_u}, it is necessary to compute the conditional covariances $C^c_{u}(\mathbf{x}, \mathbf{x}' \mid X)$ and $C^c_{ug}(\mathbf{x}, \mathbf{x}' \mid X)$, which, unlike $C^c_{g}(\mathbf{x}, \mathbf{x}' \mid X)$ given by~\eqref{eq:cond_cov}, are not available in closed form; therefore, we compute sample approximations of these covariances.
For this purpose, we draw $M$ realizations of $g^c(\mathbf{x}, \omega)$ with conditioned mean and covariance given by \eqref{eq:cond_mean} and \eqref{eq:cond_cov}, resulting in the ensemble of synthetic fields $\{ g^c_{(i)}(\mathbf{x}) \equiv g^c(\mathbf{x}, \omega_{(i)}) \}^M_{i = 1}$.
For each member of the ensemble we solve the corresponding deterministic PDE problem, resulting in the ensemble of solutions $\{ u^c_{(i)} \}^M_{i = 1}$.
Employing both ensembles, we compute the sample covariances $\hat{C}^c_g(\cdot, \cdot)$, $\hat{C}^c_u(\cdot, \cdot)$ and $\hat{C}^c_{ug}(\cdot, \cdot)$.
Substituting these sample covariances into~\eqref{eq:GP_cov_u}, and the result into~\eqref{eq:active-learning-criteria-u}, leads to the acquisition policy
\begin{equation}
  \mathbf{x}^{*} \equiv \argmin_{\mathbf{x}'} \int_D \left \{ \hat{C}^c_u(\mathbf{x}, \mathbf{x}) - \hat{C}^c_{ug}(\mathbf{x}, \mathbf{x}') [ \hat{C}^c_g(\mathbf{x}', \mathbf{x}') ]^{-1} \hat{C}^c_{gu}(\mathbf{x}', \mathbf{x}) \right \} \, \mathrm{d} \mathbf{x}.
  \label{minimization}
\end{equation}

\section{Numerical Experiments}
\label{sec:numerical}

In this section, we apply the conditional KL modeling approaches presented in Section~\ref{sec:cond_rand} and the active learning methods presented in Section~\ref{sec:active_learning} to solve the stochastic diffusion equation problem \eqref{RK:eq:spde}.
We consider the following two-dimensional steady state diffusion equation with a random diffusion coefficient over the domain $D = [0,2] \times [0,1]$, subject to Dirichlet and Neumann boundary conditions:
\begin{equation}
  \label{RK:eq:spde_nr}
  \begin{aligned}
    -\nabla \cdot (k(\mathbf{x},\omega) \nabla u(\mathbf{x},\omega)) &= f(\mathbf{x},\omega) && u(\mathbf{x},\omega)\in D\times \Omega,\\
    u(\mathbf{x},\omega) &= 1 &&  x_1 = 0,\\
    u(\mathbf{x},\omega) &= 0 &&  x_1 = 2,\\
    n \cdot \nabla u(\mathbf{x},\omega) &= 0 && x_2 = \{ 0, 1 \}.\\
  \end{aligned}
\end{equation}

\begin{figure}[ht!]
    \centering
    \begin{subfigure}[t]{0.48\textwidth}
        \centering
        \includegraphics[scale=.32]{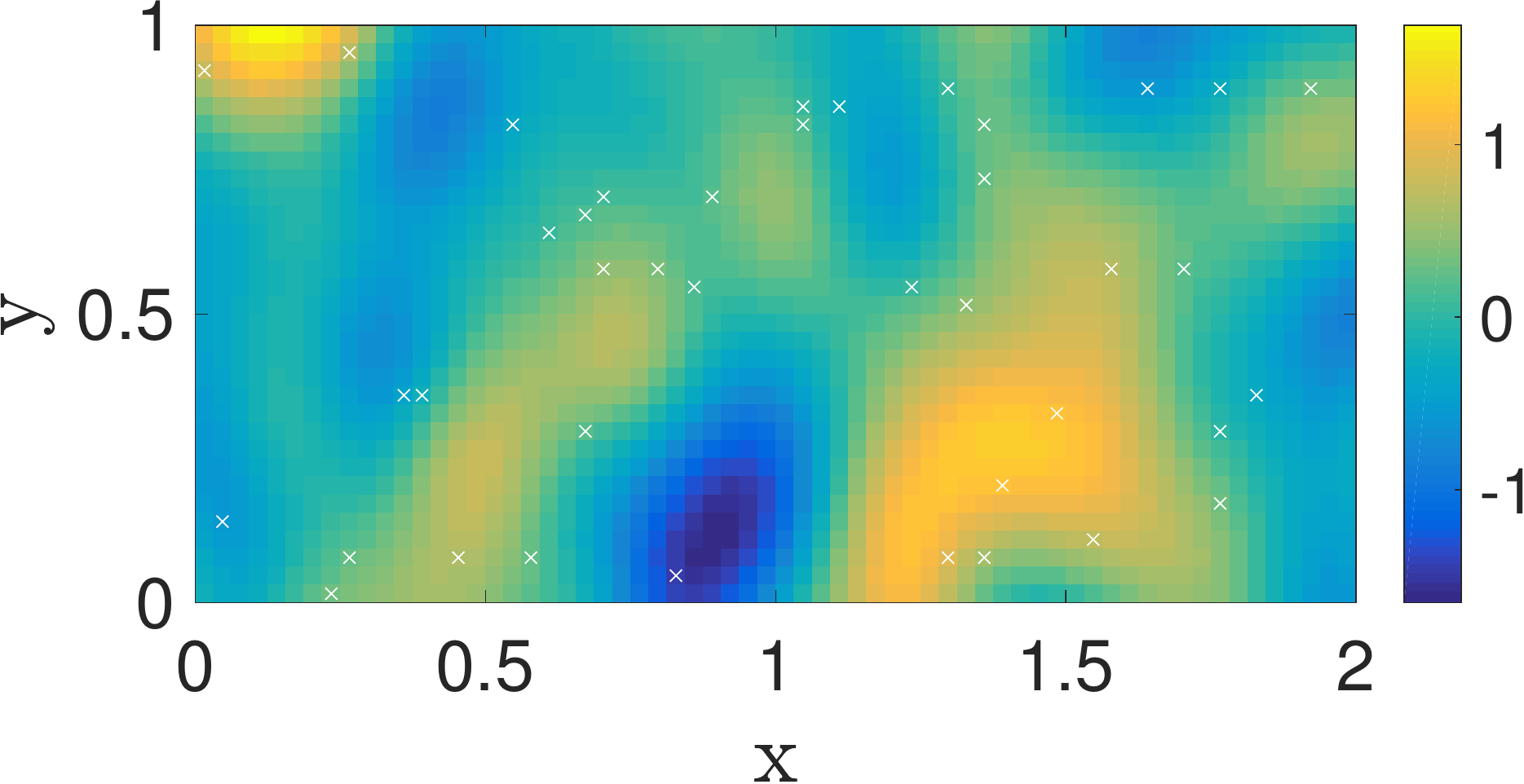}
        \caption{$g(\mathbf{x})$ field} \label{RK:fig:g_obs_sample_nobs40_mrst}
    \end{subfigure}        
    \begin{subfigure}[t]{0.48\textwidth}
        \centering
        \includegraphics[scale=.32]{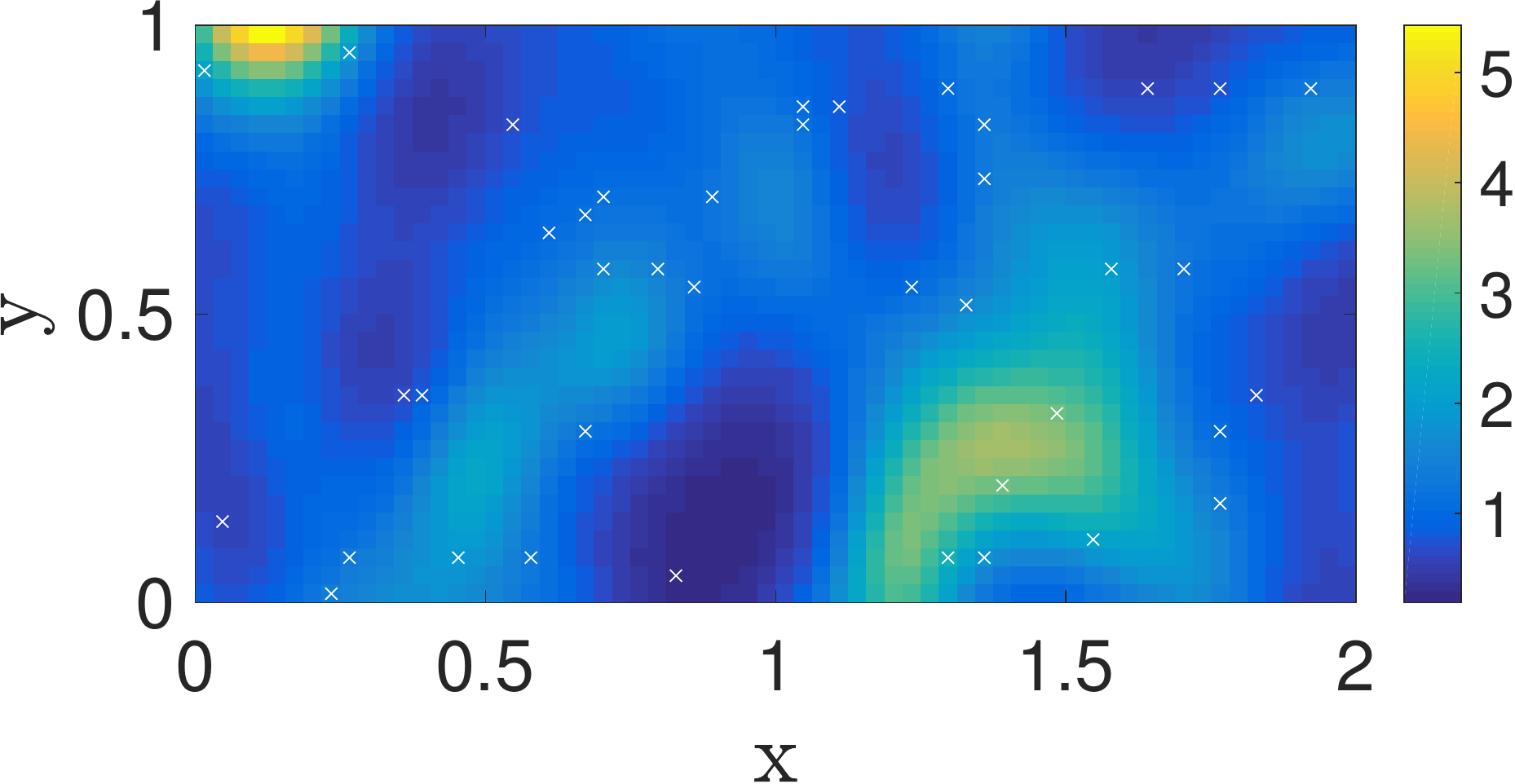}
        \caption{$k(\mathbf{x})$ field} \label{RK:fig:k_obs_sample_nobs40_mrst}
    \end{subfigure}    
       \caption{Synthetic field $k=\exp(g)$ and locations of the observations.}  \label{obs_sample_nobs40_mrst}
\end{figure}

The reference $g = \log k$ field is constructed synthetically by drawing a realization of the GP process \eqref{gp_prior} and \eqref{cov_g} with the parameters $\sigma_g = 0.65, l_x=0.15$ and $l_y=0.2$, shown in Figure~\ref{obs_sample_nobs40_mrst}.
From this reference field, we draw $40$ observations at random locations to be used for constructing the conditional GP model $g^c$ and computing conditional solution of Eq \eqref{RK:eq:spde_nr}.

\subsection{Conditional GP models}
\label{sec:numerical-conditional-gp}

We construct finite-dimensional conditional GP models for the field $g$ based on the synthetic dataset using the two approaches introduced in Section~\ref{sec:cond_rand} and employ these models to compute the mean and standard deviation of the conditional solution $u^c$ of the SPDE problem~\eqref{RK:eq:spde_nr}.
We employ MC simulations to compute reference unconditional and conditional mean and standard deviation of $g$ and $u$.
We use Approach 2 of Section~\ref{sec:cond-trunc-kl}, which provides a quantitative measure of dimension reduction due to conditioning, to estimate the dimensionality of the reduced conditional model $g^c$.
We construct finite-dimensional models for $g^c$ of the estimated dimension using both Approaches 1 and 2 and propagate uncertainty through the SPDE problem~\eqref{RK:eq:spde_nr} using the stochastic collocation method~\cite{RK:Xiu2002,Babuka2010}.

\subsubsection{Reference Monte Carlo unconditional and conditional solutions}
\label{sec:numerical-reference}

\begin{figure}[ht!]
  \centering
  \begin{subfigure}[t]{0.48\textwidth}
    \centering%
    \includegraphics[scale=.32]{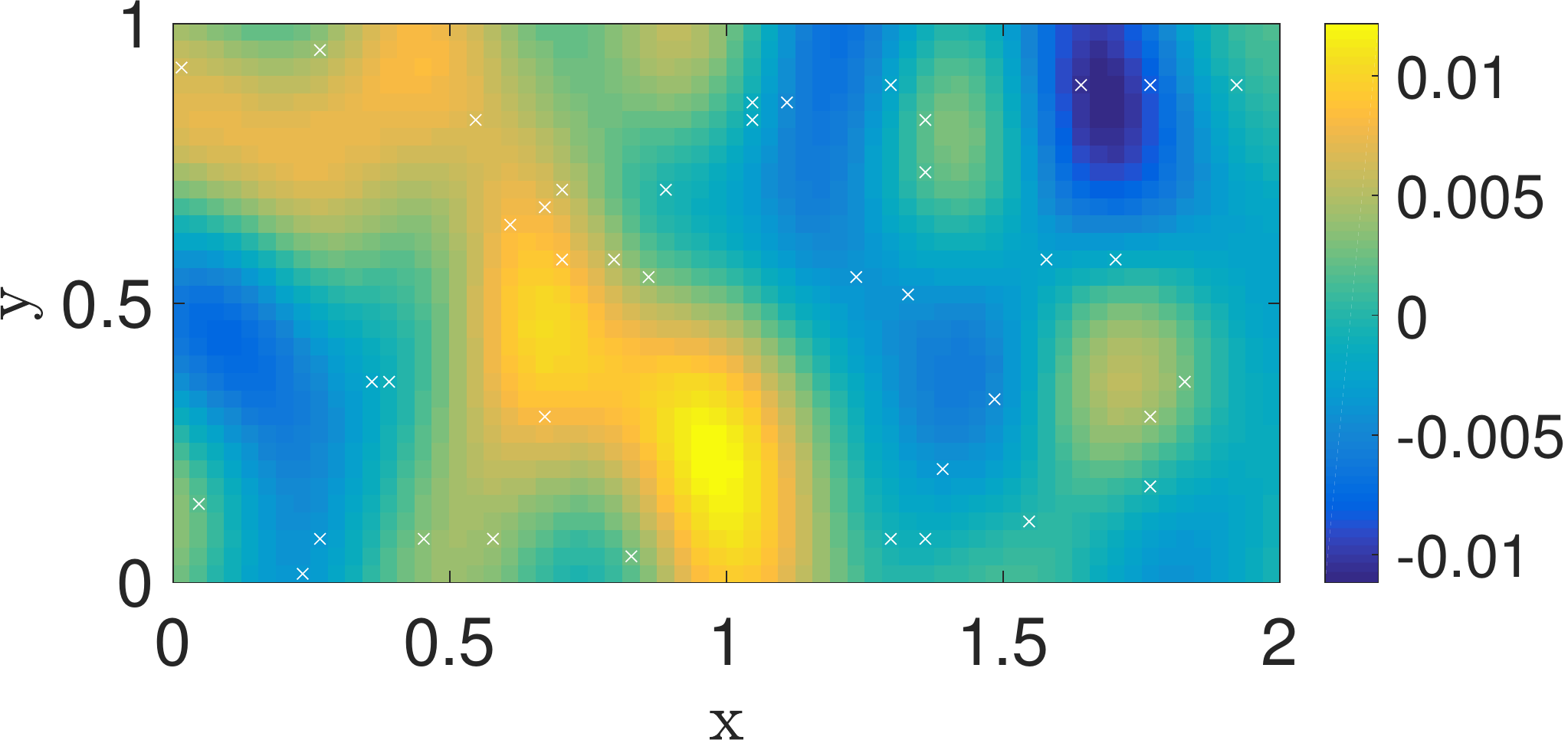}
    \caption{mean $g(\mathbf{x},\omega)$}
    \label{RK:fig:g_mean_mc_001_uncond_mrst}
  \end{subfigure}
  \begin{subfigure}[t]{0.48\textwidth}
    \centering%
    \includegraphics[scale=.32]{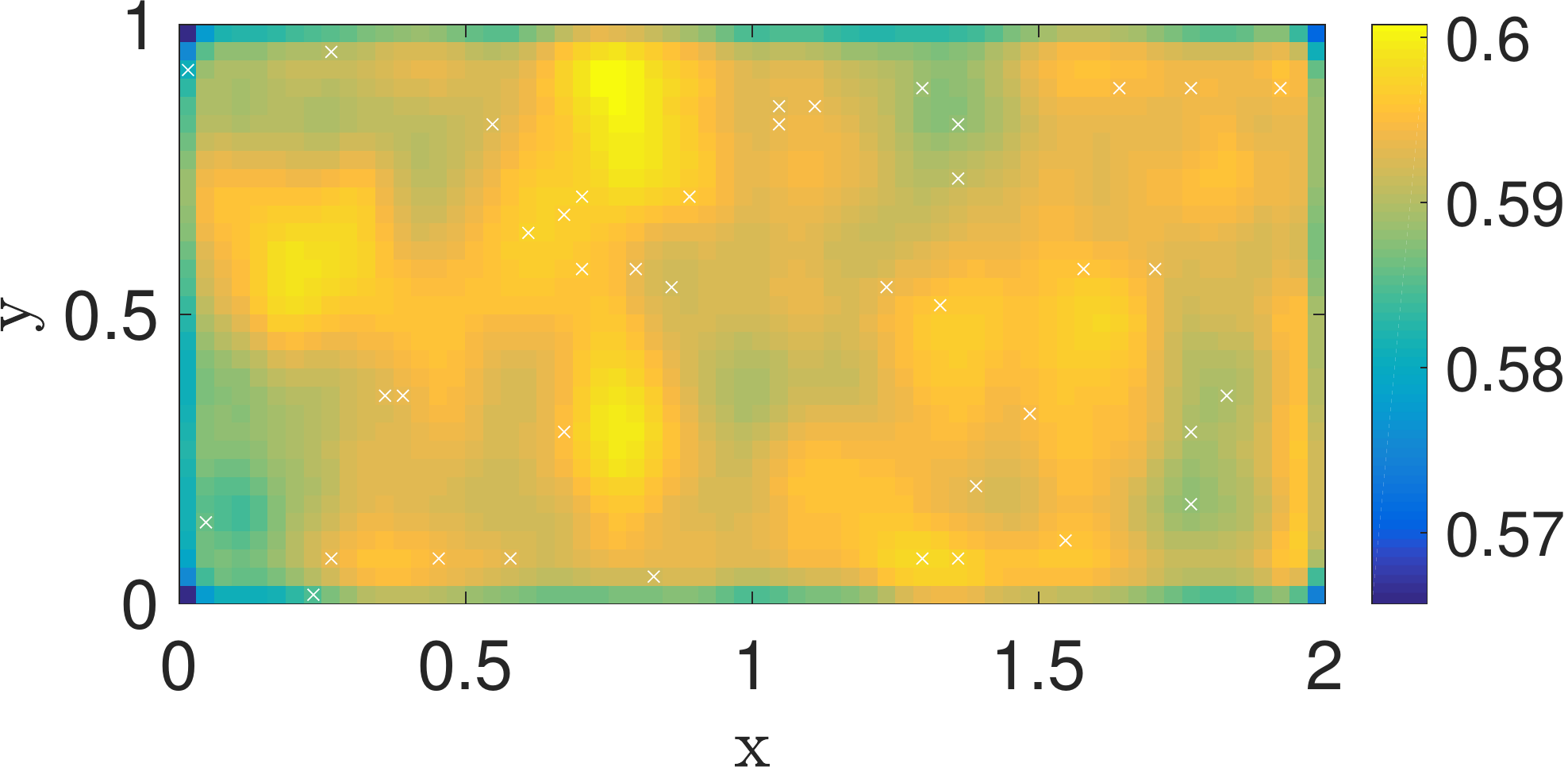}
    \caption{st. deviation of $g(\mathbf{x},\omega)$}
    \label{RK:fig:g_sdev_mc_001_uncond_mrst}
  \end{subfigure}
  \begin{subfigure}[t]{0.48\textwidth}
    \centering%
    \includegraphics[scale=.32]{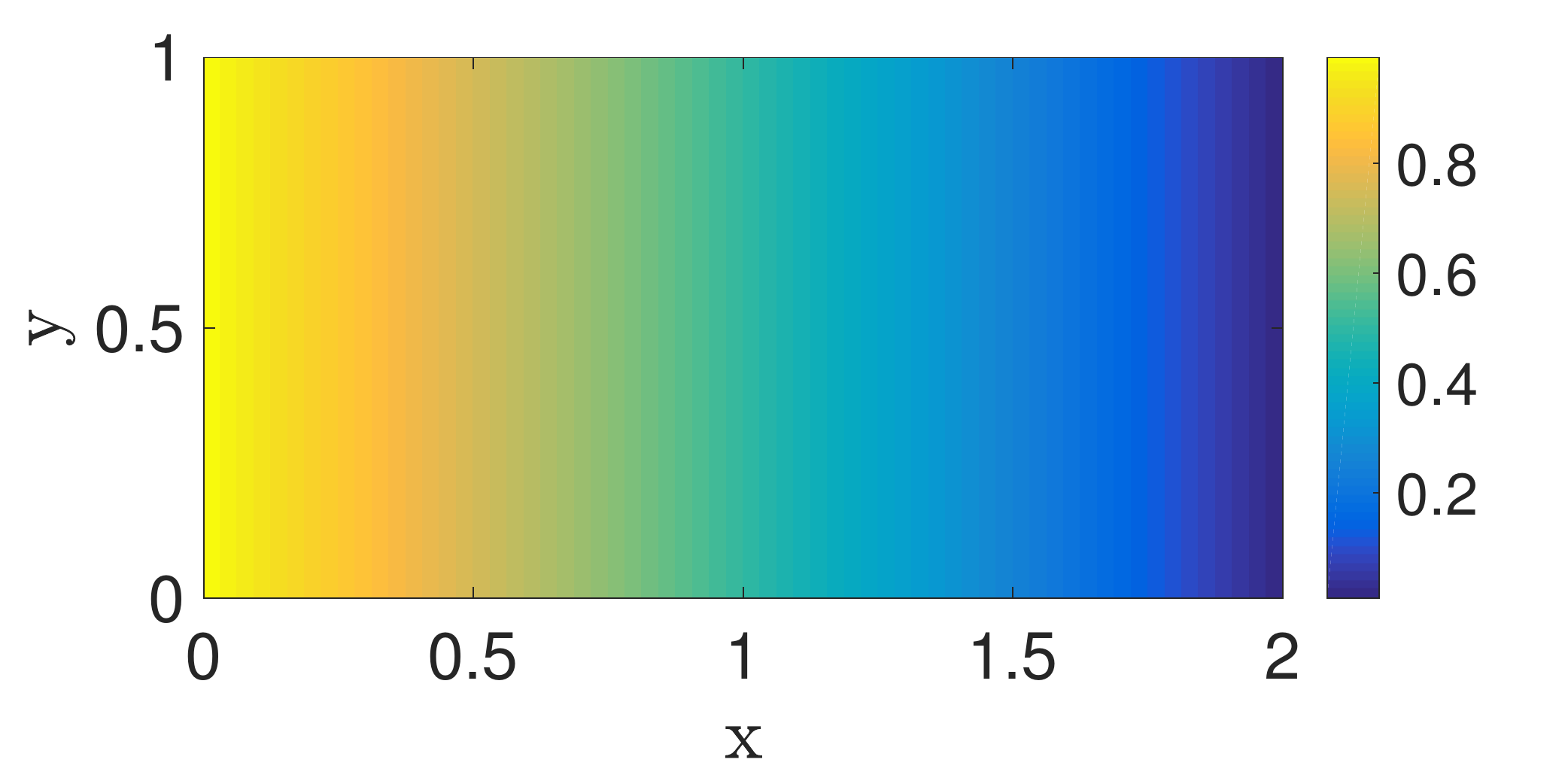}
    \caption{mean $u(\mathbf{x},\omega)$}
    \label{RK:fig:u_mean_mc_001_uncond_mrst_no_obs_plot}
  \end{subfigure}
  \begin{subfigure}[t]{0.48\textwidth}
    \centering%
    \includegraphics[scale=.32]{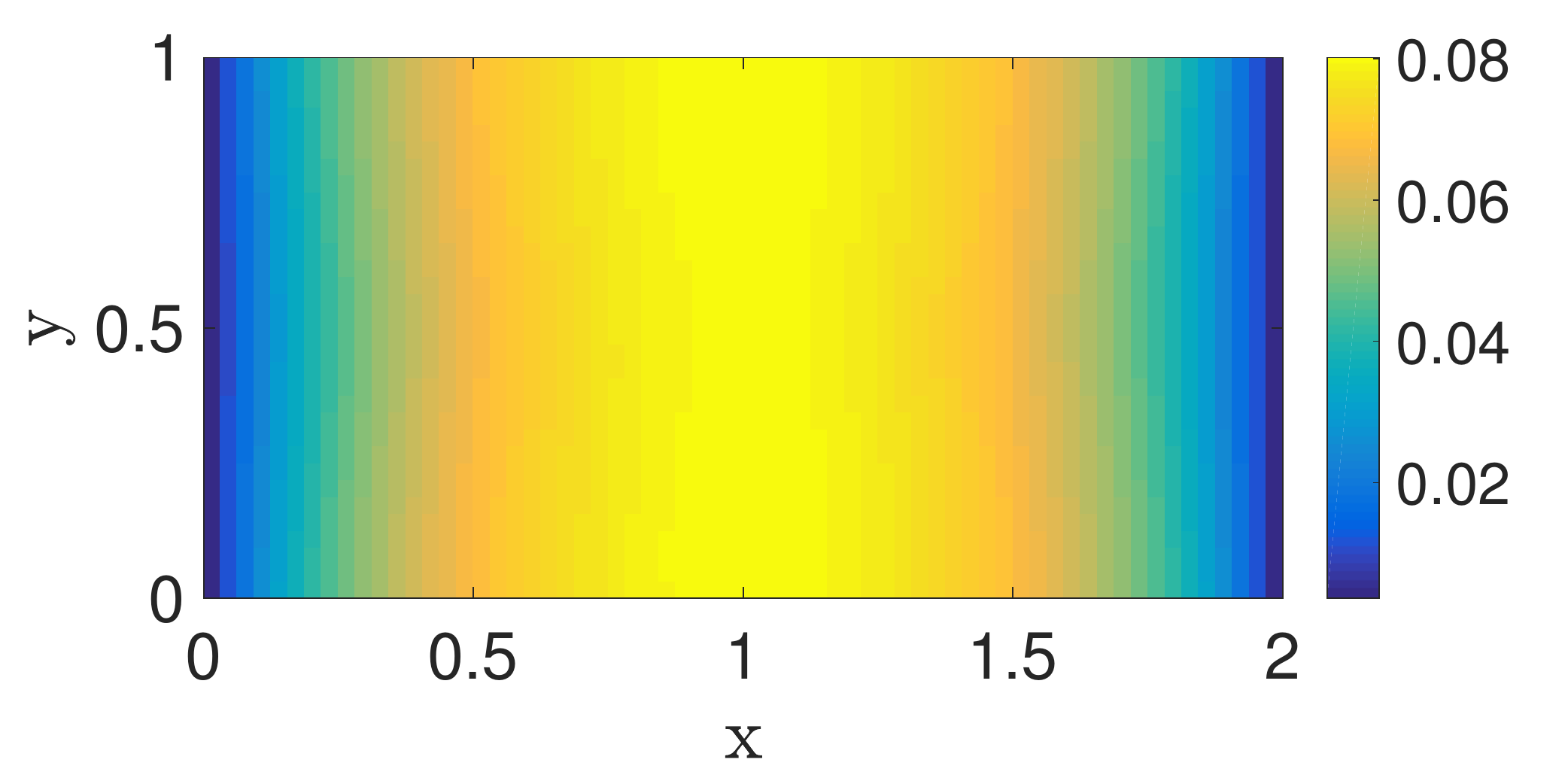}
    \caption{st. deviation of $u(\mathbf{x},\omega)$}
    \label{RK:fig:u_sdev_mc_001_uncond_mrst_no_obs_plot}
  \end{subfigure}
  \caption{Unconditional (a) mean of $g(\mathbf{x},\omega)$, (b) standard deviation of of $g(\mathbf{x},\omega)$, (c) mean of $u(\mathbf{x},\omega)$, and (d) standard deviation of $u(\mathbf{x},\omega)$ computed via MC simulation with \num{15000} realizations.}
  \label{fig:u_mc_001_uncond_mrst_no_obs_plot}
\end{figure}

\begin{figure}[ht!]
  \centering
  \begin{subfigure}[t]{0.95\textwidth}
    \centering%
    \includegraphics[scale=.4]{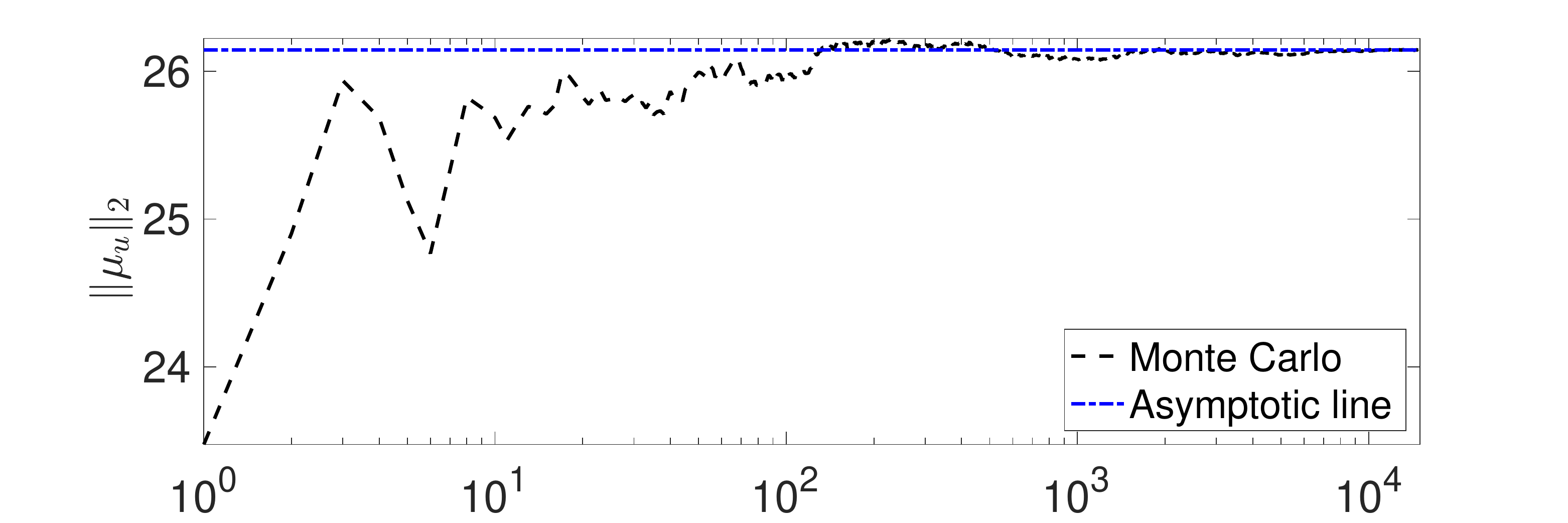}
    \caption{}
    \label{RK:fig:u_xi_mc_uncond_mean_norm_15k}
  \end{subfigure}
  \begin{subfigure}[t]{0.95\textwidth}
    \centering%
    \includegraphics[scale=.4]{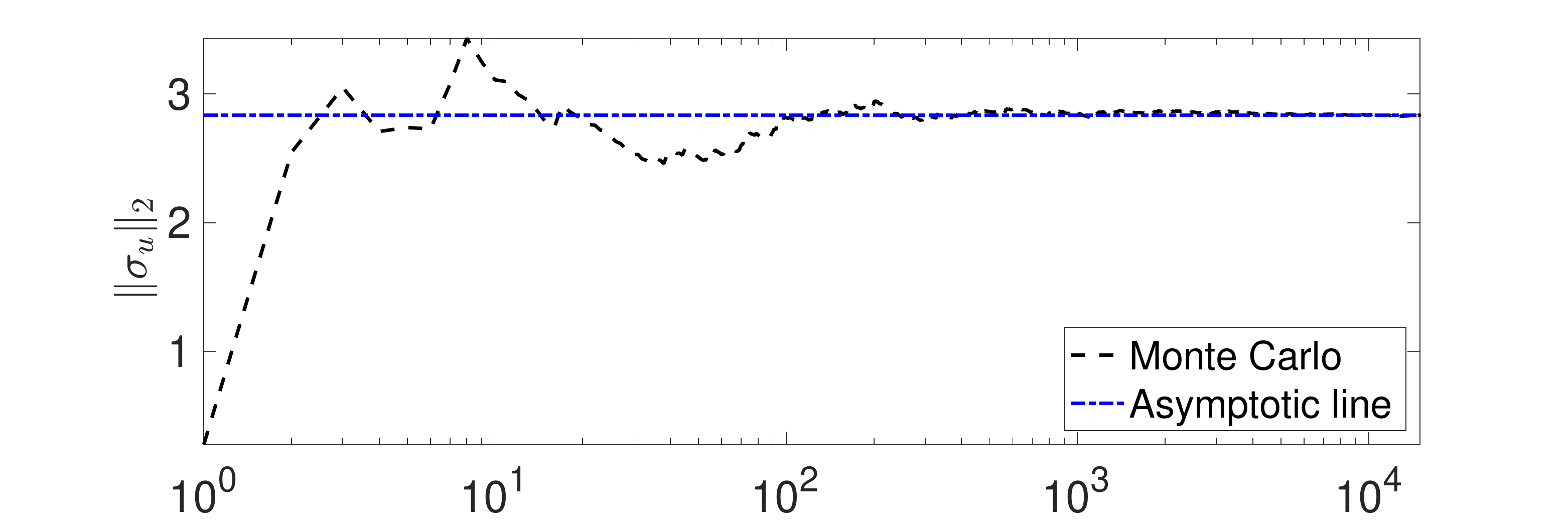}
    \caption{}
    \label{RK:fig:u_xi_mc_uncond_sdev_norm_15k}
  \end{subfigure}
  \caption{$L_2$ norm of unconditional (a) mean and (b) standard deviation of the unconditional solution $u(\mathbf{x},\omega)$, versus the number of MC realizations.}
  \label{fig:u_xi_mc_uncond_mean_sdev_norm_15k}
\end{figure}

To sample the reference unconditional $g$ field, we employ the unconditional KL expansion \eqref{eq:g-KL} truncated to $d = 60$ terms, where $d$ is chosen such as to retain $99\%$ of the unconditional variance, i.e.,  $\sum_{i=1}^d \lambda_i \geq 0.99 \operatorname{Tr}(C_g)$, where $\operatorname{Tr}(C_g)$ is the trace of $C_g$.
The unconditional mean and covariance of $u$ are computed from MC with \num{15000} samples. 
The resulting unconditional mean and standard deviation of $g$ and $u$ are presented in Figure~\ref{fig:u_mc_001_uncond_mrst_no_obs_plot}. In Figure~\ref{fig:u_xi_mc_uncond_mean_sdev_norm_15k} we present a convergence study of the MC estimators of the unconditional mean and standard deviation of $u$, where the $L_2$ norm of the estimators is plotted against the number of MC samples. This figure demonstrates that \num{15000} samples are sufficient to compute these estimators.

\begin{figure}[ht!]
  \centering
  \begin{subfigure}[t]{0.48\textwidth}
    \centering%
    \includegraphics[scale=.32]{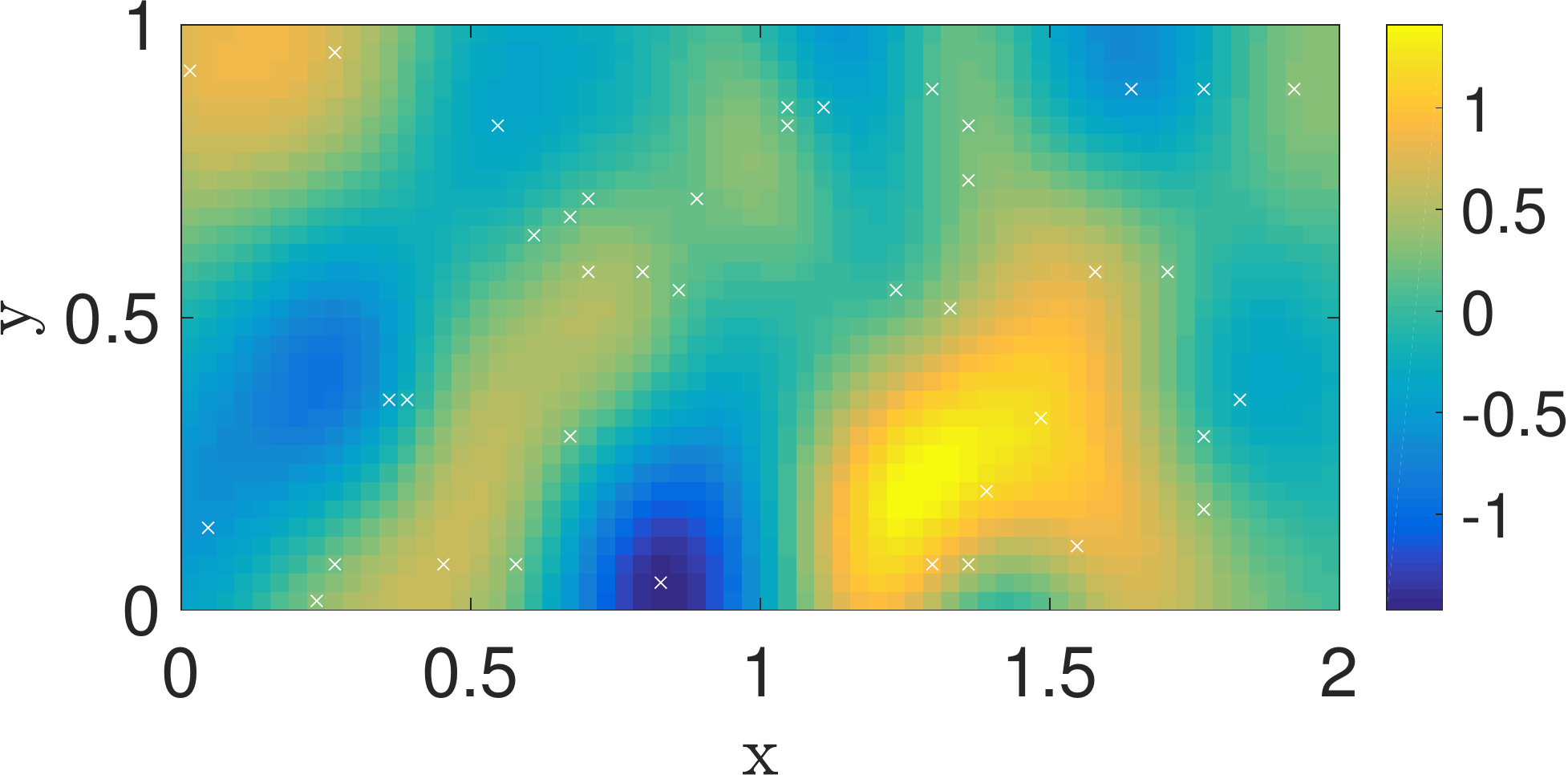}
    \caption{mean $g^c(\mathbf{x},\omega)$}
    \label{RK:fig:g_mean_exact_cond_mc_mrst}
  \end{subfigure}
  \begin{subfigure}[t]{0.48\textwidth}
    \centering%
    \includegraphics[scale=.32]{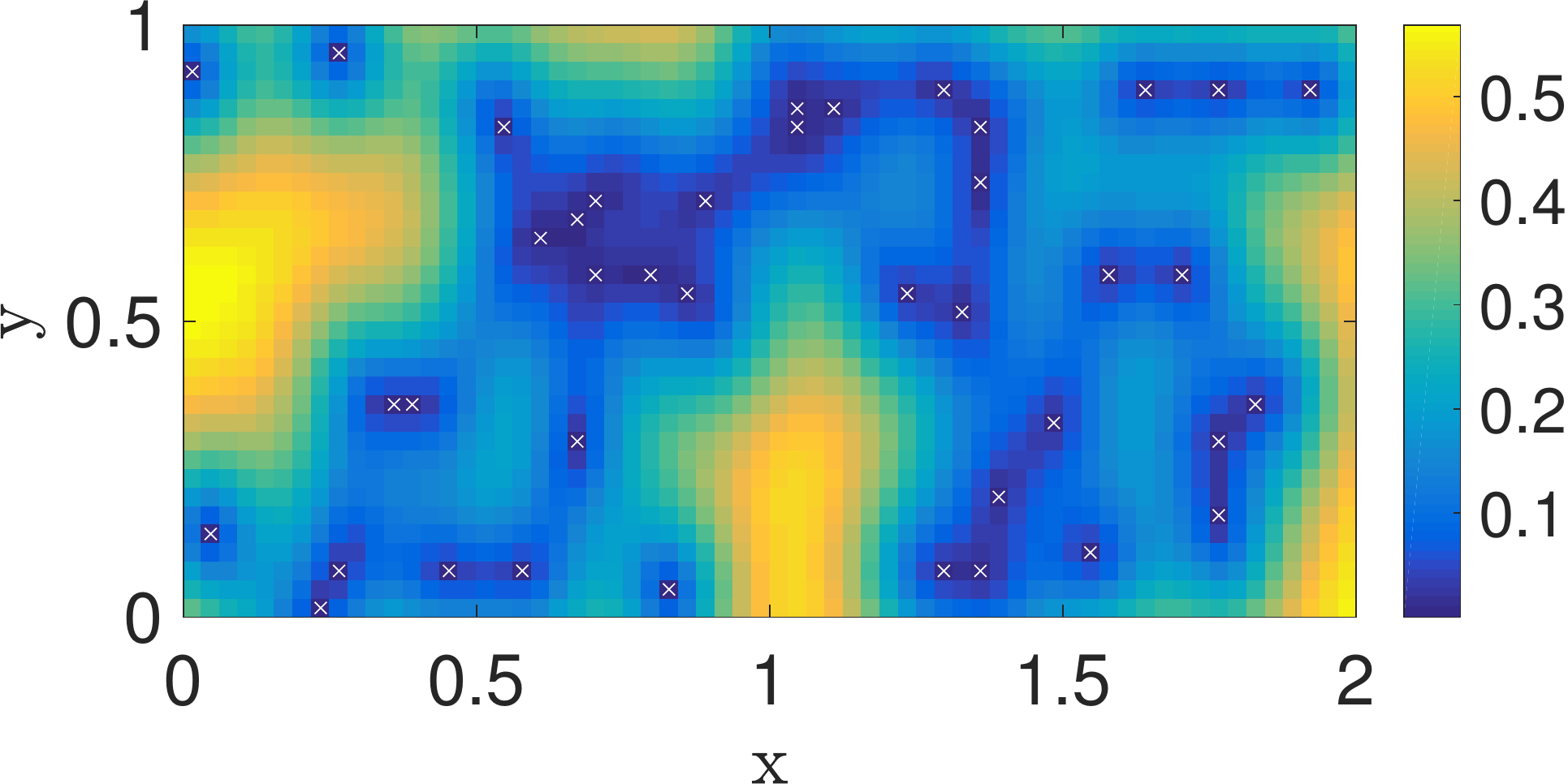}
    \caption{st. deviation of $g^c(\mathbf{x},\omega)$}
    \label{RK:fig:g_sdev_exact_cond_mc_mrst}
  \end{subfigure}
  \begin{subfigure}[t]{0.48\textwidth}
    \centering%
    \includegraphics[scale=.32]{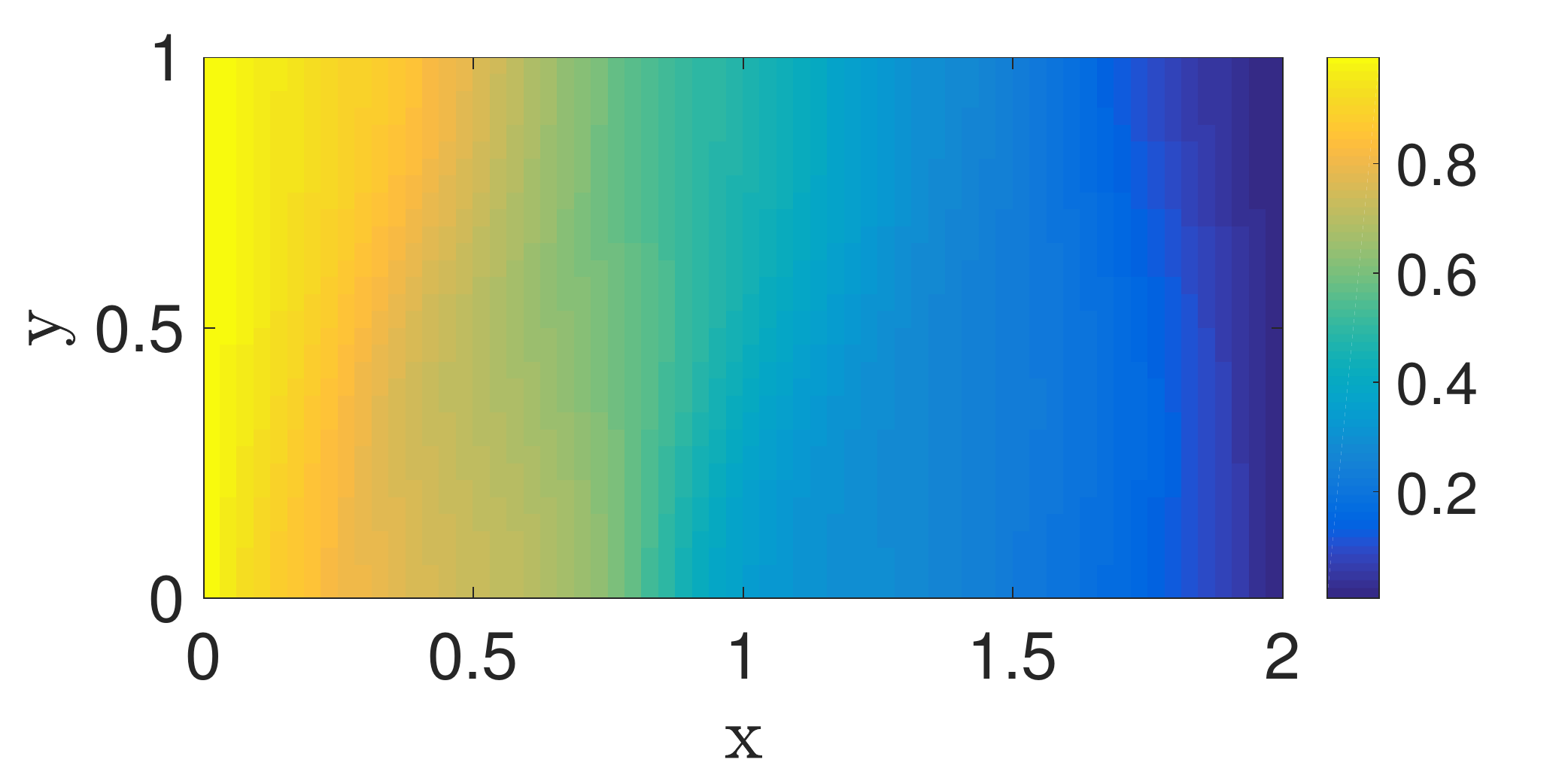}
    \caption{mean $u^c(\mathbf{x},\omega)$}
    \label{RK:fig:u_mean_exact_cond_mc_mrst_no_obs_plot}
  \end{subfigure}
  \begin{subfigure}[t]{0.48\textwidth}
    \centering%
    \includegraphics[scale=.32]{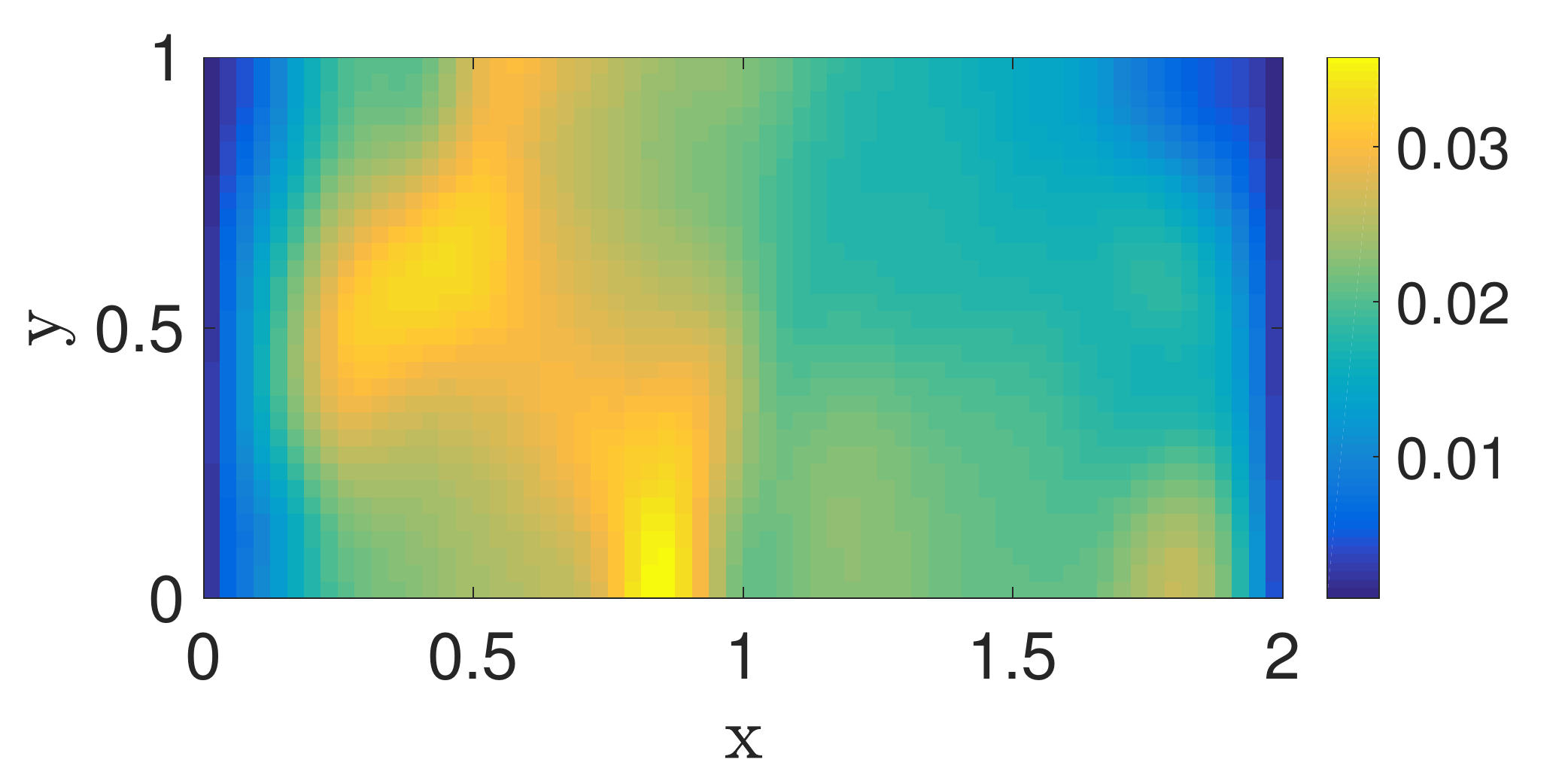}
    \caption{st. deviation of $u^c(\mathbf{x},\omega)$}
    \label{RK:fig:u_sdev_exact_cond_mc_mrst_no_obs_plot}
  \end{subfigure}
  \caption{(a) Mean and (b) standard deviation of $g^c(\mathbf{x},\omega)$ obtained from KL expansion with $53$ terms, and the corresponding (c) mean and (d) standard deviation of $u^c(\mathbf{x},\omega)$ computed via MC simulation with \num{15000} realizations}
  \label{fig:g_u_mean_exact_cond_mc_mrst}
\end{figure}

\begin{figure}[ht!]
  \centering
  \begin{subfigure}[t]{0.95\textwidth}
    \centering%
    \includegraphics[scale=.4]{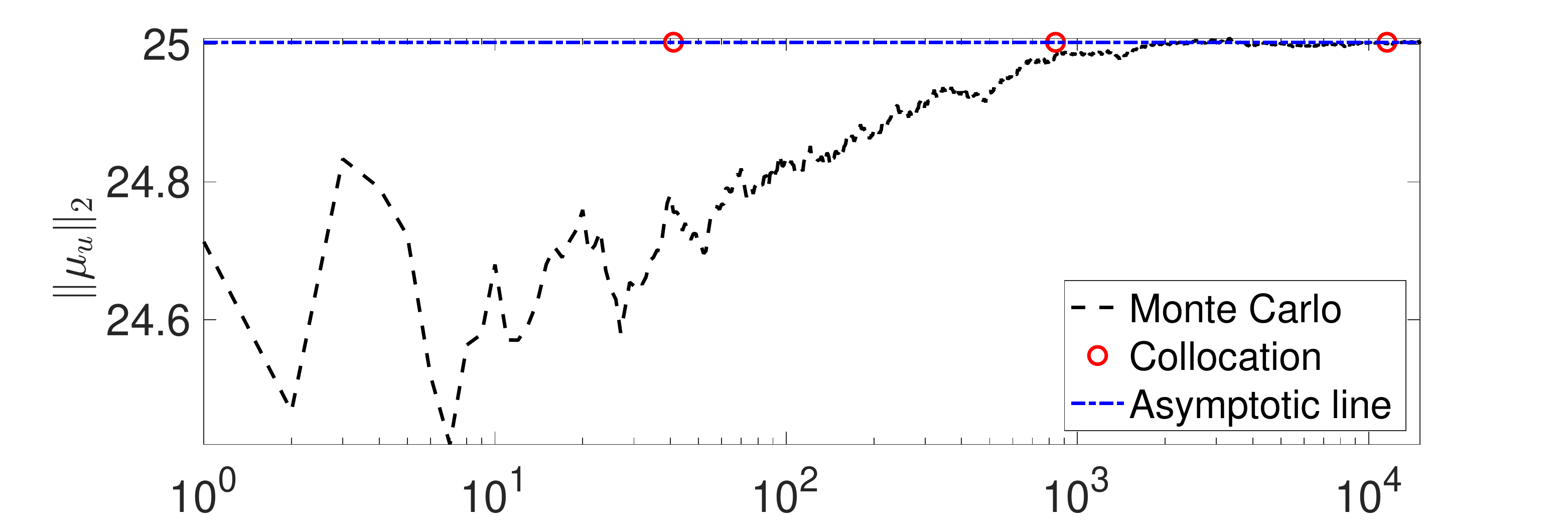}
    \caption{}
    \label{RK:fig:u_xi_mc_cond_mean_norm_15k}
  \end{subfigure}
  \begin{subfigure}[t]{0.95\textwidth}
    \centering%
    \includegraphics[scale=.4]{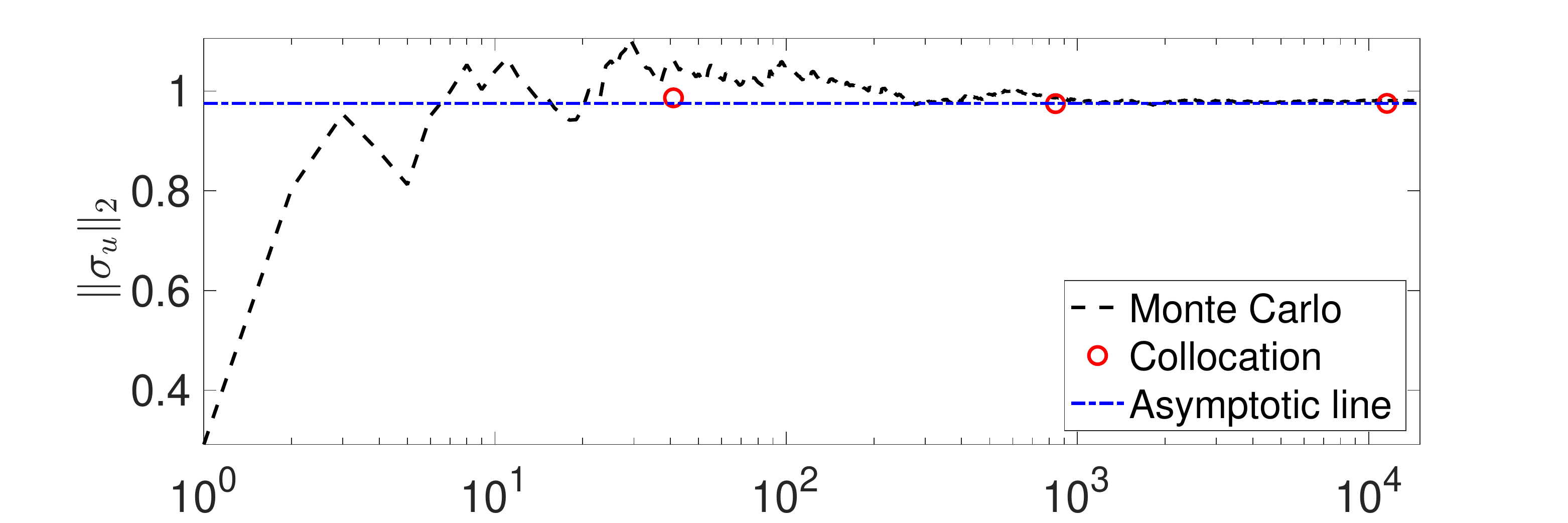}
    \caption{}
    \label{RK:fig:u_xi_mc_cond_sdev_norm_15k}
  \end{subfigure}
  \caption{$L_2$ norm of (a) mean and (b) standard deviation of the conditional solution $u(\mathbf{x},\omega)$, versus the number of realizations in the MC solution (black line) and stochastic collocation points (red circles).
    The number of collocation points are \num{41}, \num{841} and \num{11561} corresponding to the sparse grid levels 2, 3, and 4, respectively.}
  \label{fig:u_xi_mc_cond_mean_sdev_norm_15k}
\end{figure}

Similarly, to sample the reference conditional field, we compute the mean and covariance using Eqs \eqref{eq:cond_mean} and \eqref{eq:cond_cov}, compute the KL expansion of the conditional covariance, and truncate the resulting conditional KL expansion to $d = 53$ terms to retain $99\%$ of the conditional variance.
The conditional mean and covariance of $g$ and $u$ are computed using MC with \num{15000} samples. The results are presented in Figure~\ref{fig:g_u_mean_exact_cond_mc_mrst}. The convergence study of the MC estimators of the conditional mean and standard deviation of $u$ is shown in Figure~\ref{fig:u_xi_mc_cond_mean_sdev_norm_15k} and shows that 15000 samples are sufficient to compute these estimators.

\begin{figure}[ht!]
  \centering
  \begin{subfigure}[t]{0.48\textwidth}
    \centering%
    \includegraphics[scale=.32]{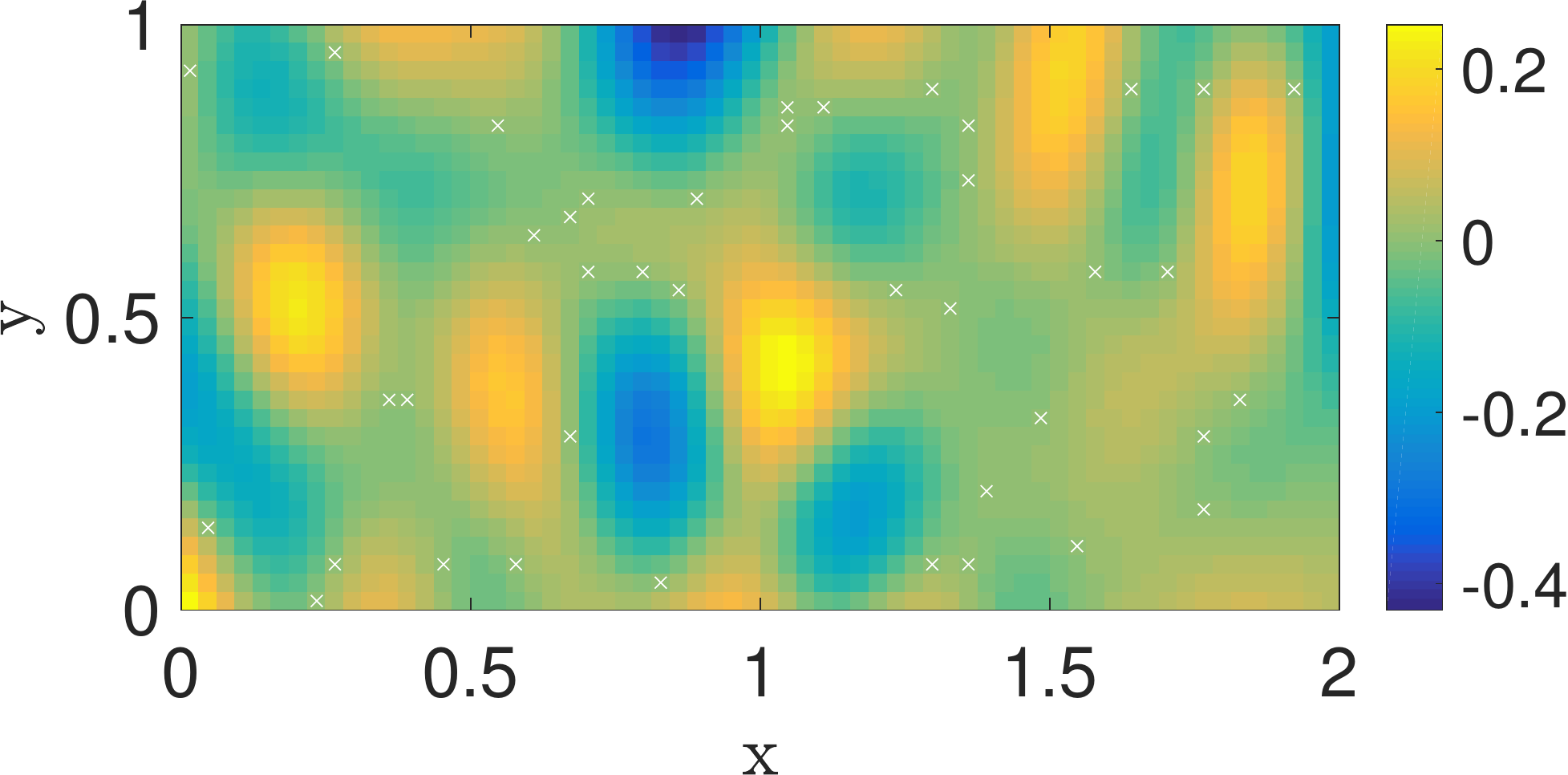}
    \caption{mean $g^c(\mathbf{x},\omega)$}%
    \label{RK:fig:gpm_est_exact-g_mu_gauss_xi_error}
  \end{subfigure}
  \begin{subfigure}[t]{0.48\textwidth}
    \centering%
    \includegraphics[scale=.32]{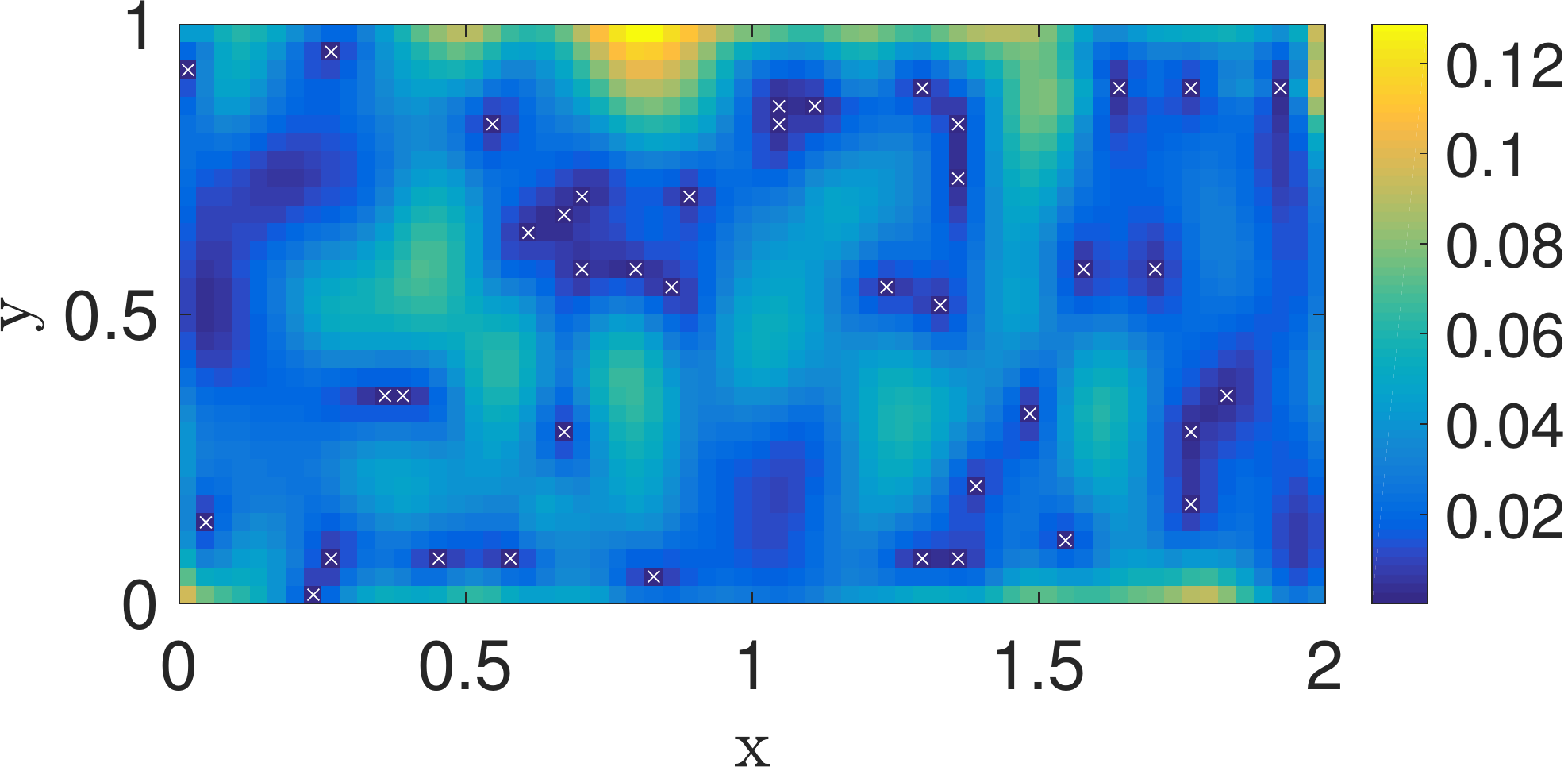}
    \caption{st. deviation of $g^c(\mathbf{x},\omega)$}%
    \label{RK:fig:C_pg_est_exact-g_sdev_gauss_xi_error}
  \end{subfigure}
  \begin{subfigure}[t]{0.48\textwidth}
    \centering%
    \includegraphics[scale=.32]{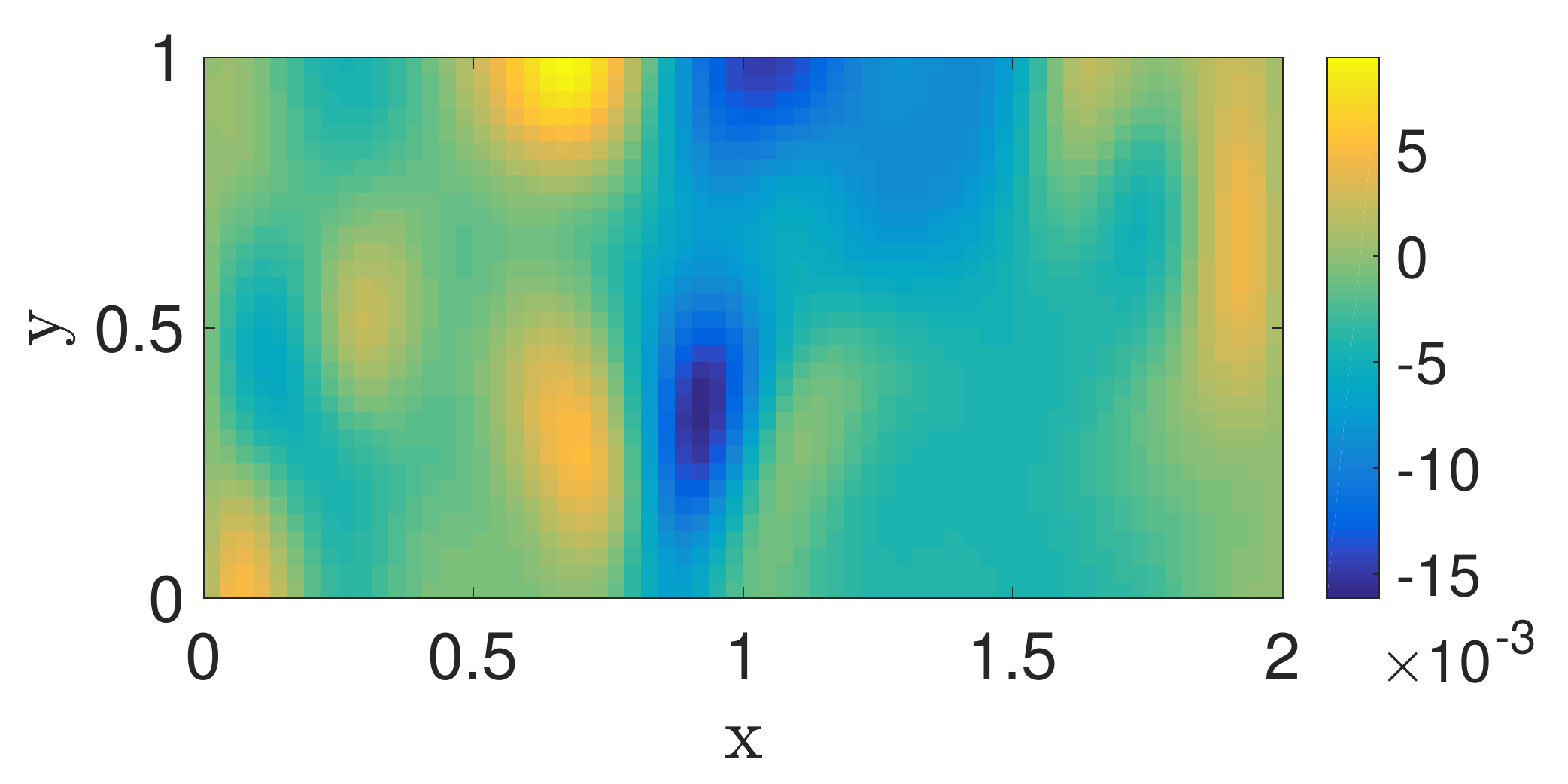}
    \caption{mean $u^c(\mathbf{x},\omega)$}%
    \label{RK:fig:u_mu_exact_cond_mc-u_mu_gauss_xi_error_no_obs_plot}
  \end{subfigure}
  \begin{subfigure}[t]{0.48\textwidth}
    \centering%
    \includegraphics[scale=.21]{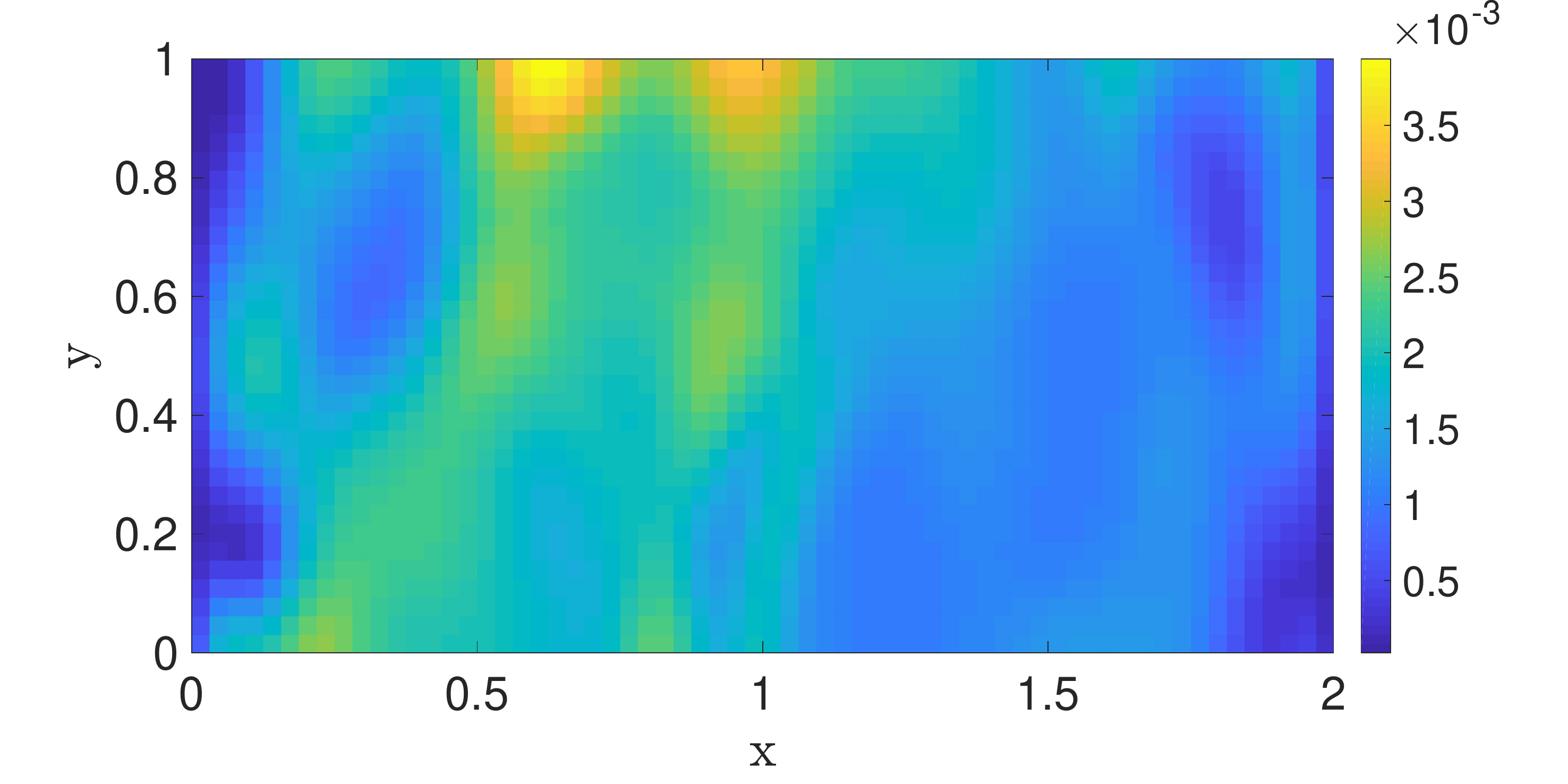}
    \caption{st. deviation of $u^c(\mathbf{x},\omega)$}%
    \label{RK:fig:u_sdev_exact_cond_mc-u_sdev_gauss_xi_error_no_obs_plot}
  \end{subfigure}
  \caption{Absolute point error in  (a) mean of $g^c(\mathbf{x},\omega)$, (b) standard deviation of $g^c(\mathbf{x},\omega)$, (c). mean of $u^c(\mathbf{x},\omega)$, and (d) standard deviation of $u^c(\mathbf{x},\omega)$ obtained with KL expansion truncated first and then conditioned (Approach 2) and the sparse grid collocation method.} \label{g_gauss_xi_d20_red_cond_mrst}
\end{figure}

The comparison of Figures~\ref{RK:fig:g_sdev_mc_001_uncond_mrst} and \ref{RK:fig:g_sdev_exact_cond_mc_mrst} show that the standard deviation of $g$ is reduced after conditioning on $g$ measurements. 
Similarly, the standard deviation of $u$, shown in Figures~\ref{RK:fig:u_sdev_mc_001_uncond_mrst_no_obs_plot} and \ref{RK:fig:u_sdev_exact_cond_mc_mrst_no_obs_plot}, also is reduced after conditioning on $g$ measurements. 

\subsubsection{Conditional solution computed using Approach 2: Conditioning truncated KL expansion of the unconditioned field}
\label{sec:numerical-approach-2}

In this section we employ Approach 2, presented in Section~\ref{sec:cond-trunc-kl}, for constructing a finite-dimensional conditional KL model of the random coefficient of the SPDE problem~\eqref{RK:eq:spde_nr}.
In Section~\ref{sec:numerical-reference}, we determined that 60 dimensions are required to represent the unconditional $g$ field.
By virtue of~\eqref{eq:r-dim}, it follows that conditioning these $d = 60$ dimensions on $N_s = 40$ observations reduces the dimensionality of the $g^c$ KL representation to $r = 20$.
As this approach for constructing a discretized conditional model disregards information provided by the higher eigenpairs of the unconditional expansion, we expect the resulting conditional moments (mean and standard deviation) of $g^c$ to differ from the reference moments computed in Section~\ref{sec:numerical-reference}.
The absolute point-wise errors in the mean and standard deviation of $g^c$ are shown in Figures~\ref{RK:fig:gpm_est_exact-g_mu_gauss_xi_error} and~\ref{RK:fig:C_pg_est_exact-g_sdev_gauss_xi_error}, respectively. 

Next, we employ the conditional KL expansion to estimate the mean and variance of $u^c$ using the sparse grid collocation method~\cite{RK:Xiu2002,Babuka2010} with \num{41}, \num{841}, and \num{11561} quadrature points.
Figure~\ref{fig:u_xi_mc_cond_mean_sdev_norm_15k} shows the $L_2$ norm of the estimators. It can be seen that \num{841} points are sufficient to obtain a convergent solution for both mean and standard deviation of $u^c$.
Absolute point-wise errors of the mean and standard deviation of $u^c$ with respect to the reference conditional solution are shown in Figures~\ref{RK:fig:u_mu_exact_cond_mc-u_mu_gauss_xi_error_no_obs_plot} and~\ref{RK:fig:u_sdev_exact_cond_mc-u_sdev_gauss_xi_error_no_obs_plot}, respectively.

\subsubsection{Conditioned solution computed using Approach 1: Truncated KL expansion of the conditional field}
\label{sec:numerical-approach-1}

\begin{figure}[ht!]
  \centering
  \begin{subfigure}[t]{0.48\textwidth}
    \centering%
    \includegraphics[scale=.3]{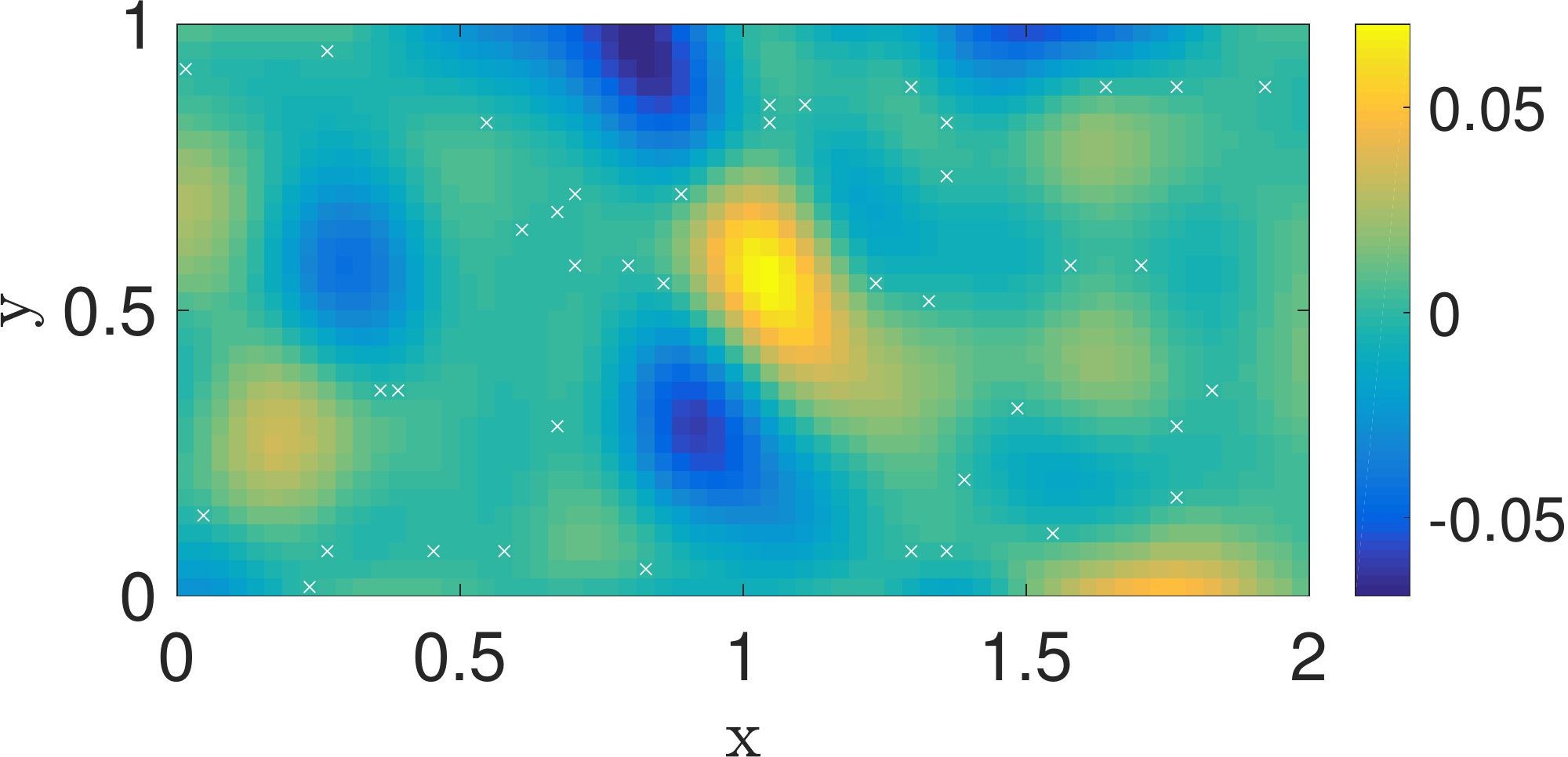}
    \caption{mean $g^c(\mathbf{x},\omega)$}%
    \label{RK:fig:gpm_est_exact-g_mu_gauss_xi_cond_trunc_error}
  \end{subfigure}
  \begin{subfigure}[t]{0.48\textwidth}
    \centering%
    \includegraphics[scale=.3]{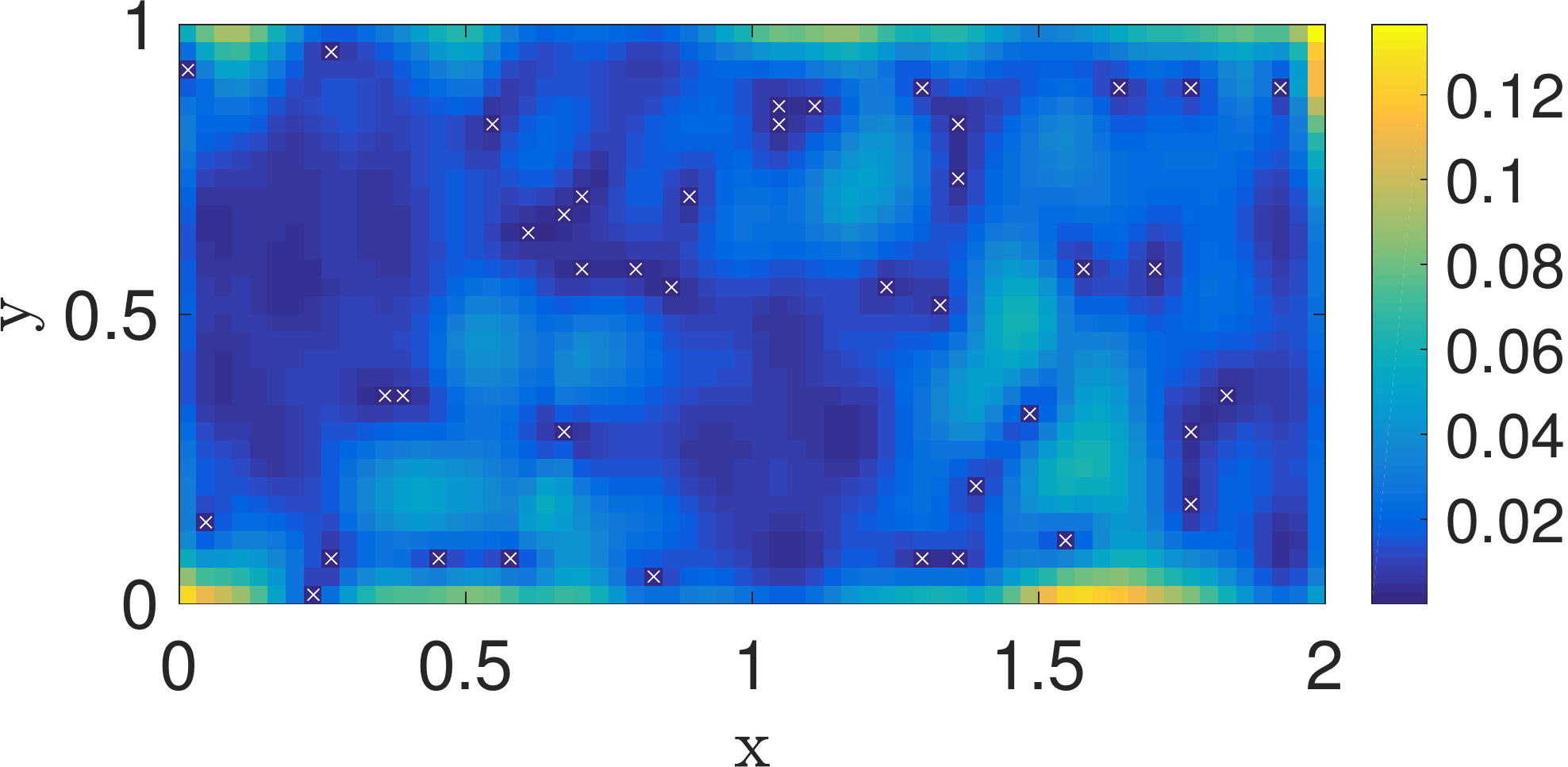}
    \caption{st. deviation of $g^c(\mathbf{x},\omega)$}%
    \label{RK:fig:C_pg_est_exact-g_sdev_gauss_xi_cond_trunc_error}
  \end{subfigure}
  \begin{subfigure}[t]{0.48\textwidth}
    \centering%
    \includegraphics[scale=.32]{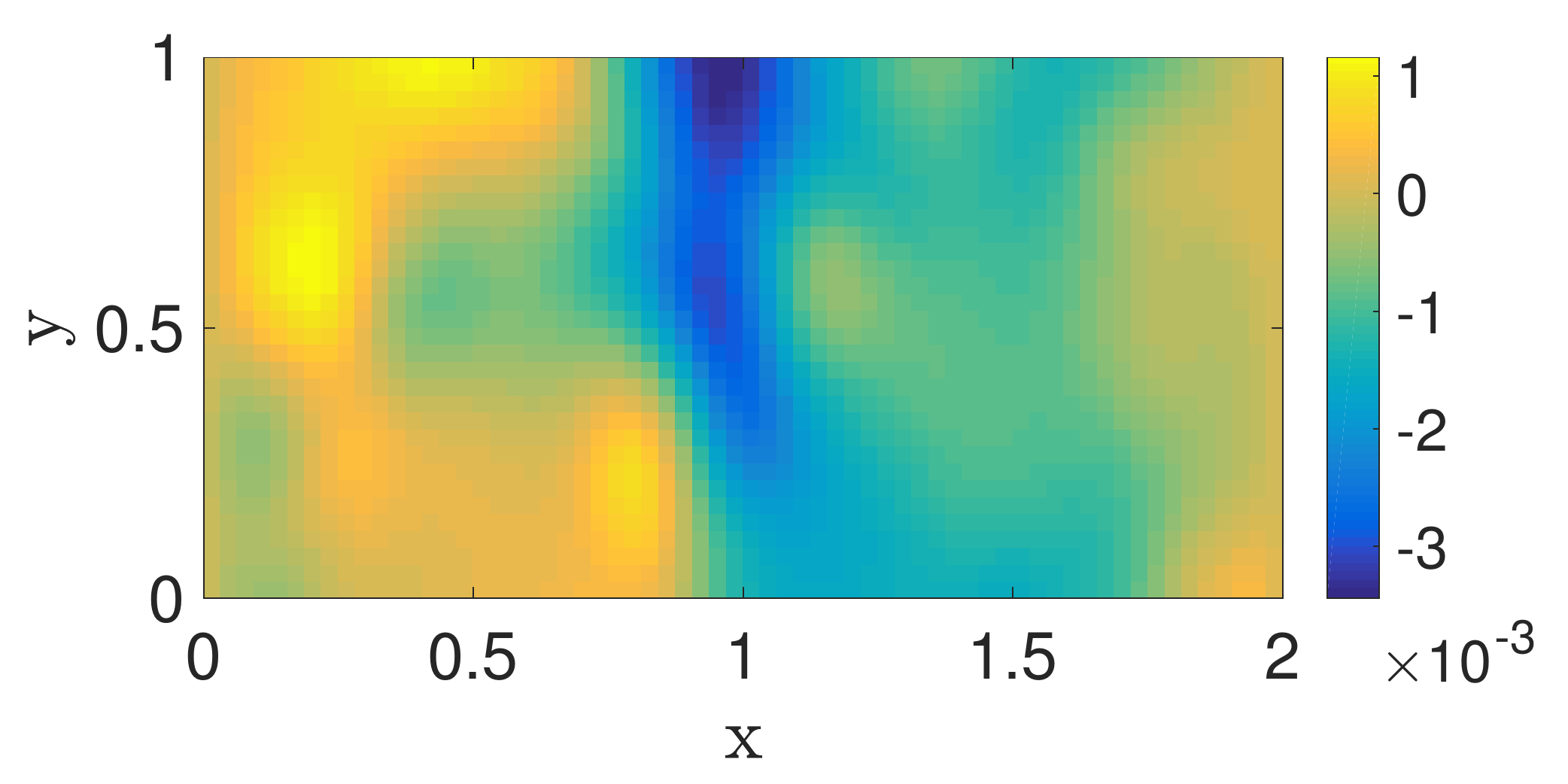}
    \caption{mean $u^c(\mathbf{x},\omega)$}%
    \label{RK:fig:u_mu_exact_cond_mc-u_mu_gauss_xi_cond_trunc_error_no_obs_plot}
  \end{subfigure}
  \begin{subfigure}[t]{0.48\textwidth}
    \centering%
    \includegraphics[scale=.21]{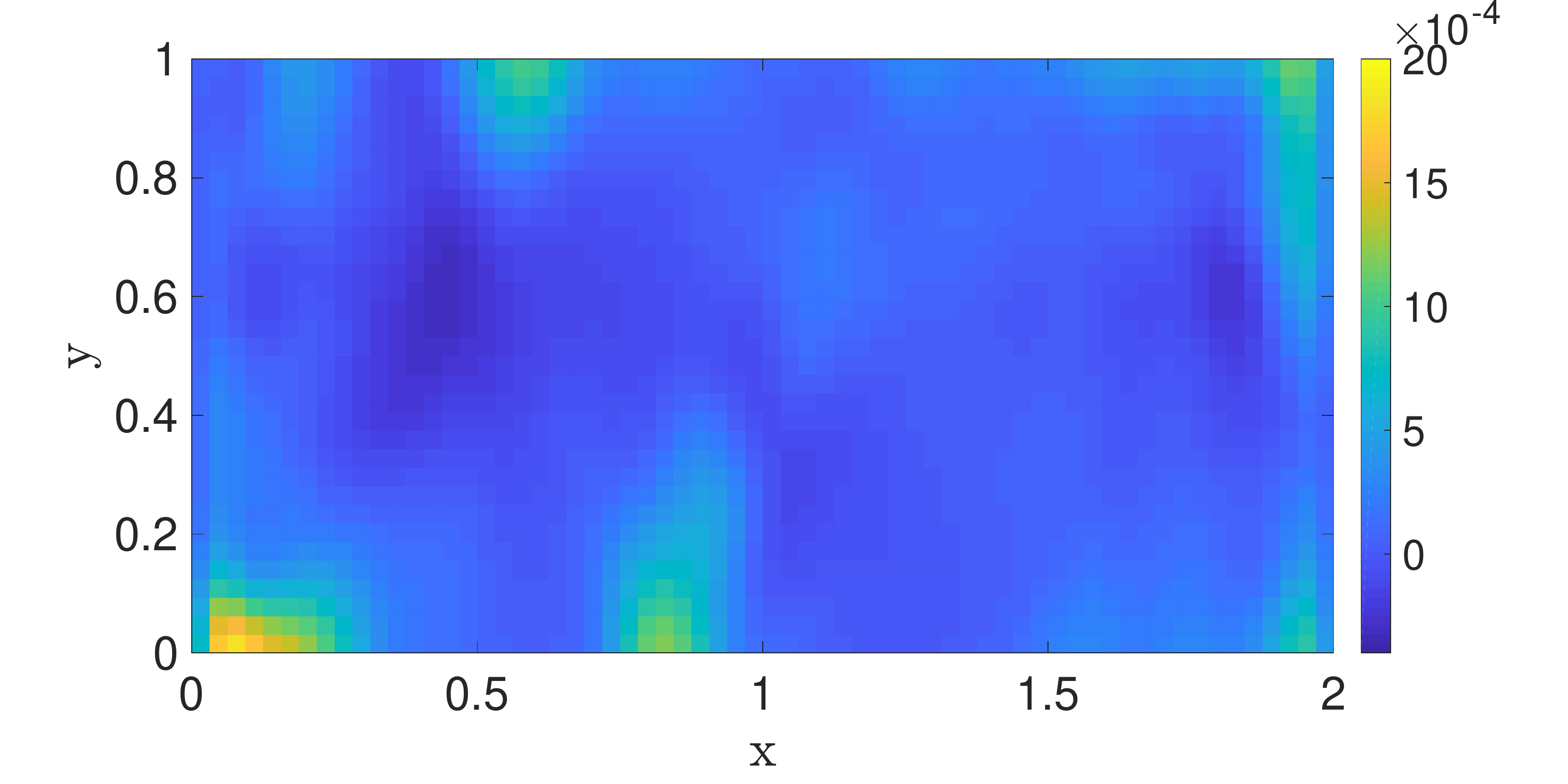}
    \caption{st. deviation of $u^c(\mathbf{x},\omega)$}%
    \label{RK:fig:u_sdev_exact_cond_mc-u_sdev_gauss_xi_cond_trunc_error_no_obs_plot}
  \end{subfigure}
  \caption{Error in (a) mean of $g^c(\mathbf{x},\omega)$, (b) standard deviation of $g^c(\mathbf{x},\omega)$, (c) mean of $u^c(\mathbf{x},\omega)$, and (d) standard deviation of $u^c(\mathbf{x},\omega)$ obtained using conditioning first and then truncated KL expansion (Approach 1) and the sparse grid collocation method.} \label{g_u_gauss_xi_d20_cond_trunc_mrst.}
\end{figure}

Next, we apply Approach 1 for constructing a finite-dimensional conditional KL model for the random coefficient of the SPDE problem \eqref{RK:eq:spde_nr}.
In the previous section, it was determined that the conditional KL expansion obtained with Approach 2 has 20 dimensions. 
To compare the accuracy of the two conditional KL approaches, here we truncate the conditional KL expansion of $g^c$ to $20$ dimensions.  
Note that in Section~\ref{sec:numerical-reference} we demonstrated that $53$ dimensions are necessary to retain $99\%$ of the variance of the exact conditional field $g^c$.
Therefore, by truncating the KL expansion to $20$ dimensions, we incur in the absolute point-wise errors in the conditional mean and standard deviation of $g^c$ shown in Figs.~\ref{RK:fig:gpm_est_exact-g_mu_gauss_xi_cond_trunc_error} and~\ref{RK:fig:C_pg_est_exact-g_sdev_gauss_xi_cond_trunc_error}.

As in the previous section, we employ the conditional KL expansion to estimate the mean and variance of $u^c$ using the sparse grid collocation method~\cite{RK:Xiu2002,Babuka2010}. % with \num{11561} quadrature points corresponding to $20$ stochastic dimensions and sparse grid level of 4.
Figures~\ref{RK:fig:u_mu_exact_cond_mc-u_mu_gauss_xi_cond_trunc_error_no_obs_plot} and \ref{RK:fig:u_sdev_exact_cond_mc-u_sdev_gauss_xi_cond_trunc_error_no_obs_plot} show the absolute point-wise error in the mean and standard deviation of $u^c$ with respect to the reference moments computed in Section~\ref{sec:numerical-reference}.

Comparing Figures~\ref{g_gauss_xi_d20_red_cond_mrst} and~\ref{g_u_gauss_xi_d20_cond_trunc_mrst.}, it can be seen that, for a set number of stochastic dimensions of the $g^c$ KL expansion, Approach 1 provides a more accurate approximation of the moments of $g^c$ and $u^c$.
On the other hand, through Eq.~\eqref{eq:r-dim}, Approach 2 provides a \emph{priori} estimate of the dimensionality of the conditional KL model sufficient to obtain an accurate solution. 

In order to study why Approach 1 provides a more accurate approximation of conditional moments, we compare the eigendecompositions of $g^c(x,\omega)$ provided by Approaches 1 and 2.
Here, we note that Approach 2 does not provide an explicit eigendecomposition (cf. Eq.~\eqref{eq:trunc-KL-g-cond-Psi-r}).
Therefore, we compute the corresponding eigendecomposition by first computing the covariance matrix induced by Approach 2, given by $(\tilde{\Psi}^r)^{\top} (\tilde{\Psi}^r)$, and then compute its eigendecomposition.

\begin{figure}[ht!]
  \centering
  \begin{subfigure}[t]{0.48\textwidth}
    \centering%
    \includegraphics[scale=.32]{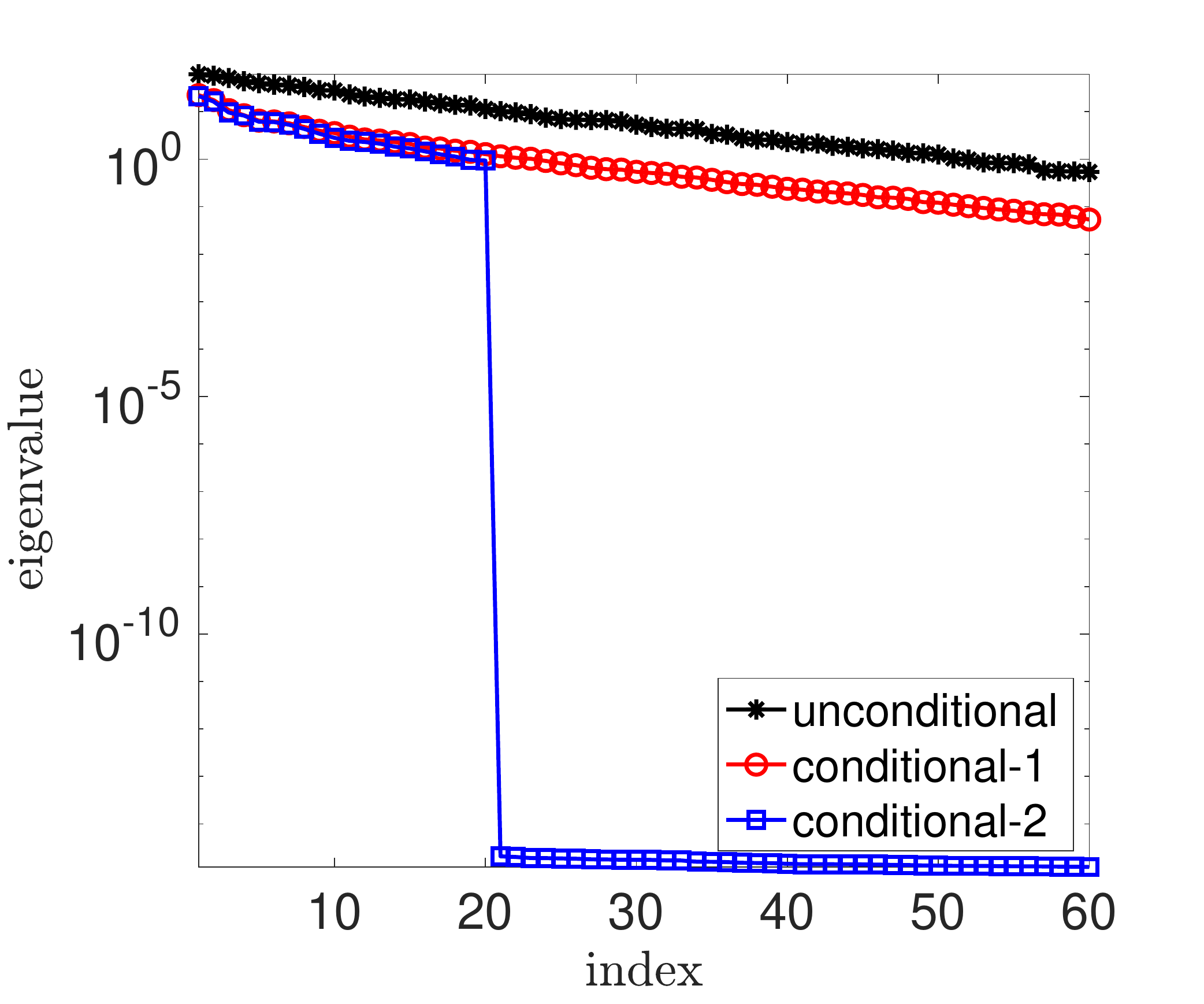}
    \caption{}
    \label{RK:fig:u_xi_mc_eigen_uncod_cond}
  \end{subfigure}
  \begin{subfigure}[t]{0.48\textwidth}
    \centering%
    \includegraphics[scale=.32]{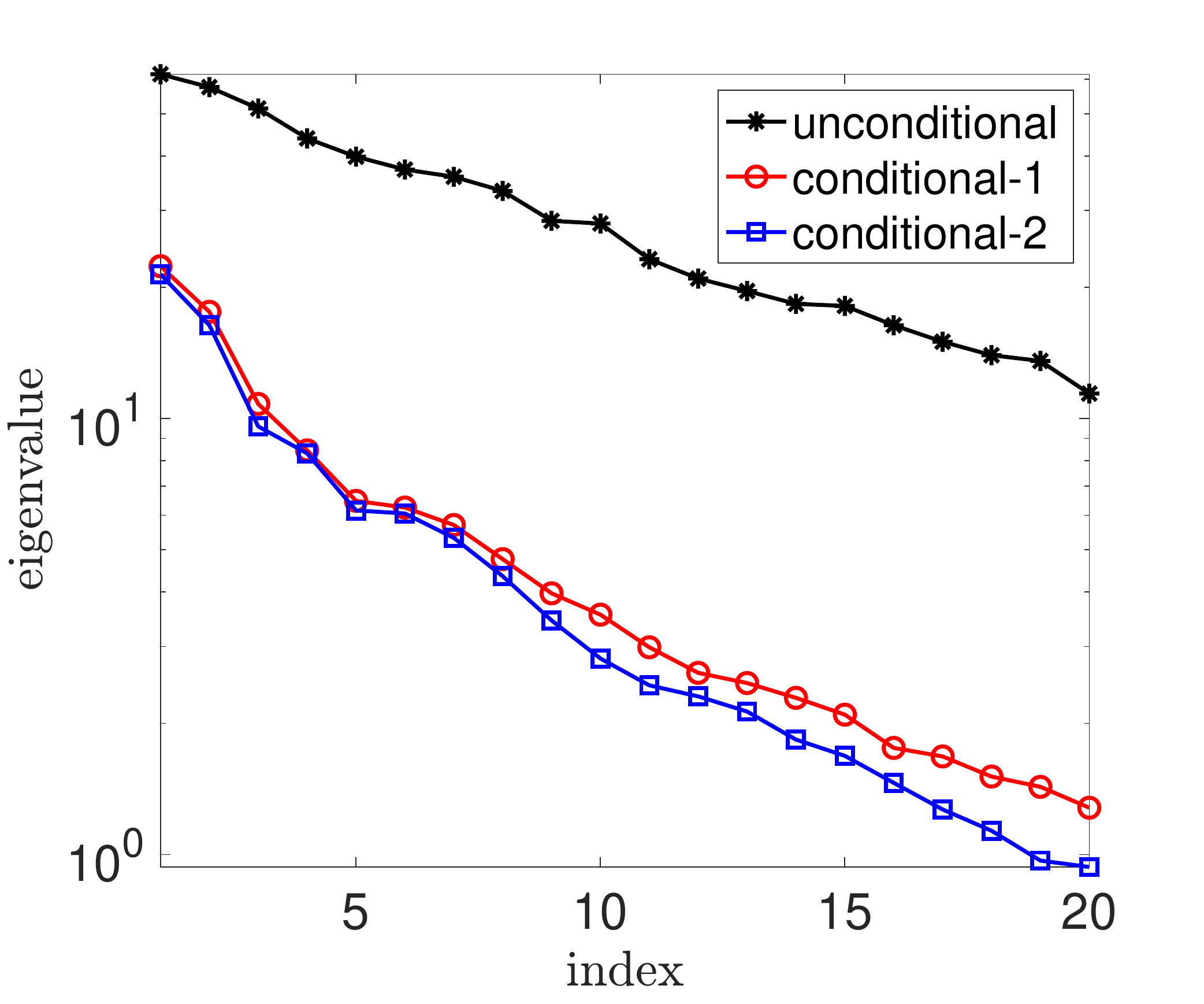}
    \caption{}
    \label{RK:fig:u_xi_mc_eigen_uncod_cond_red}
  \end{subfigure}
  \caption{Eigenvalue decay of $g(\mathbf{x}, \omega)$ (black) and $g^c(\mathbf{x}, \omega)$ (Approach 1, red, Approach 2, blue) random fields. (a) First 60 eigenvalues and (b) first 20 eigenvalues.}
  \label{u_xi_mc_eigen_uncod_cond}
\end{figure}

Figure~\ref{u_xi_mc_eigen_uncod_cond} shows the eigenvalues of $g^c(x,\omega)$ resulting from Approaches 1 and 2, together with the eigenvalues of $g(x,\omega)$.
It can be seen that the magnitude of the eigenvalues of $g^c(x,\omega)$ is smaller than those of $g(x,\omega)$, which follows from the fact that the variance of the conditional KL expansions is smaller than the variance of the unconditional KL expansion. 
It can also be seen that the conditional eigenvalues decay faster than the unconditional eigenvalues, especially for the first ten eigenvalues.
% \textcolor{red}{Rama will add eigen values for Approach 2.}

\begin{figure}[ht!]
  \centering
  \begin{subfigure}[t]{0.48\textwidth}
    \centering%
    \includegraphics[scale=.32]{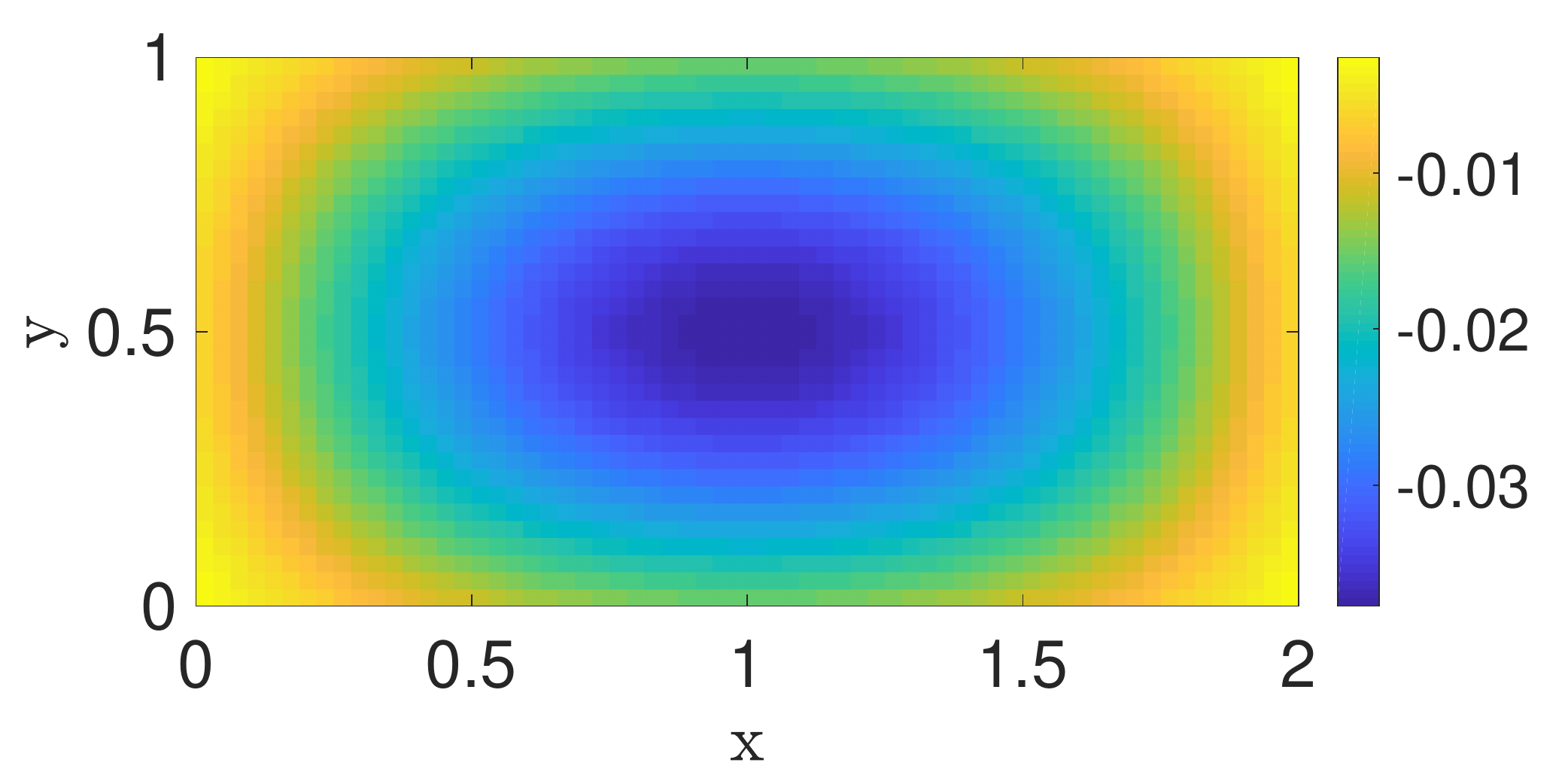}
    \caption{$g(\mathbf{x},\omega)$}
    \label{RK:fig:u_xi_mc_eigen_uncond_eig_func_1}
  \end{subfigure}
  \begin{subfigure}[t]{0.48\textwidth}
    \centering%
    \includegraphics[scale=.32]{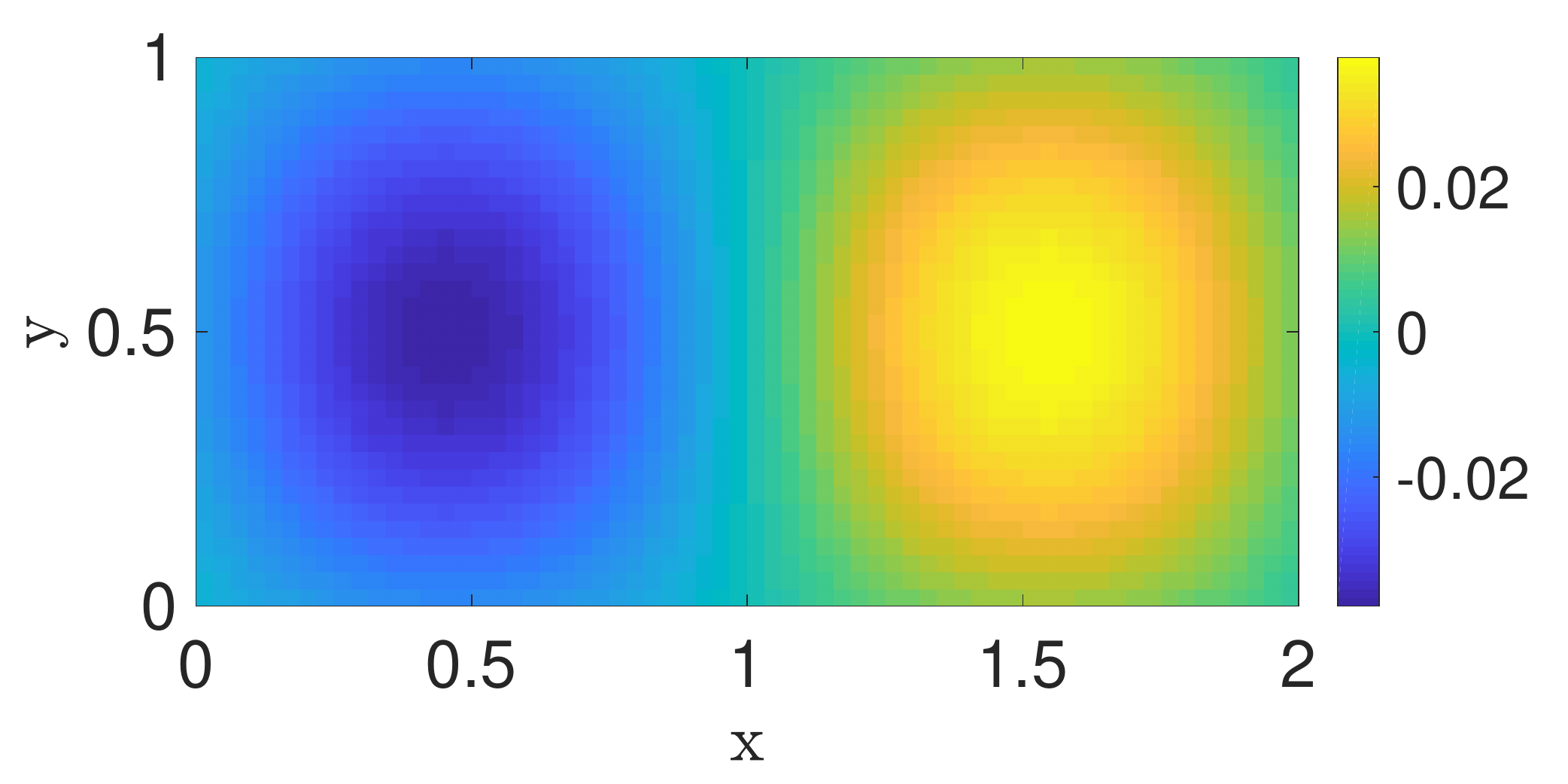}
    \caption{$g(\mathbf{x},\omega)$}
    \label{RK:fig:u_xi_mc_eigen_uncond_eig_func_2}
  \end{subfigure}
  \begin{subfigure}[t]{0.48\textwidth}
    \centering%
    \includegraphics[scale=.32]{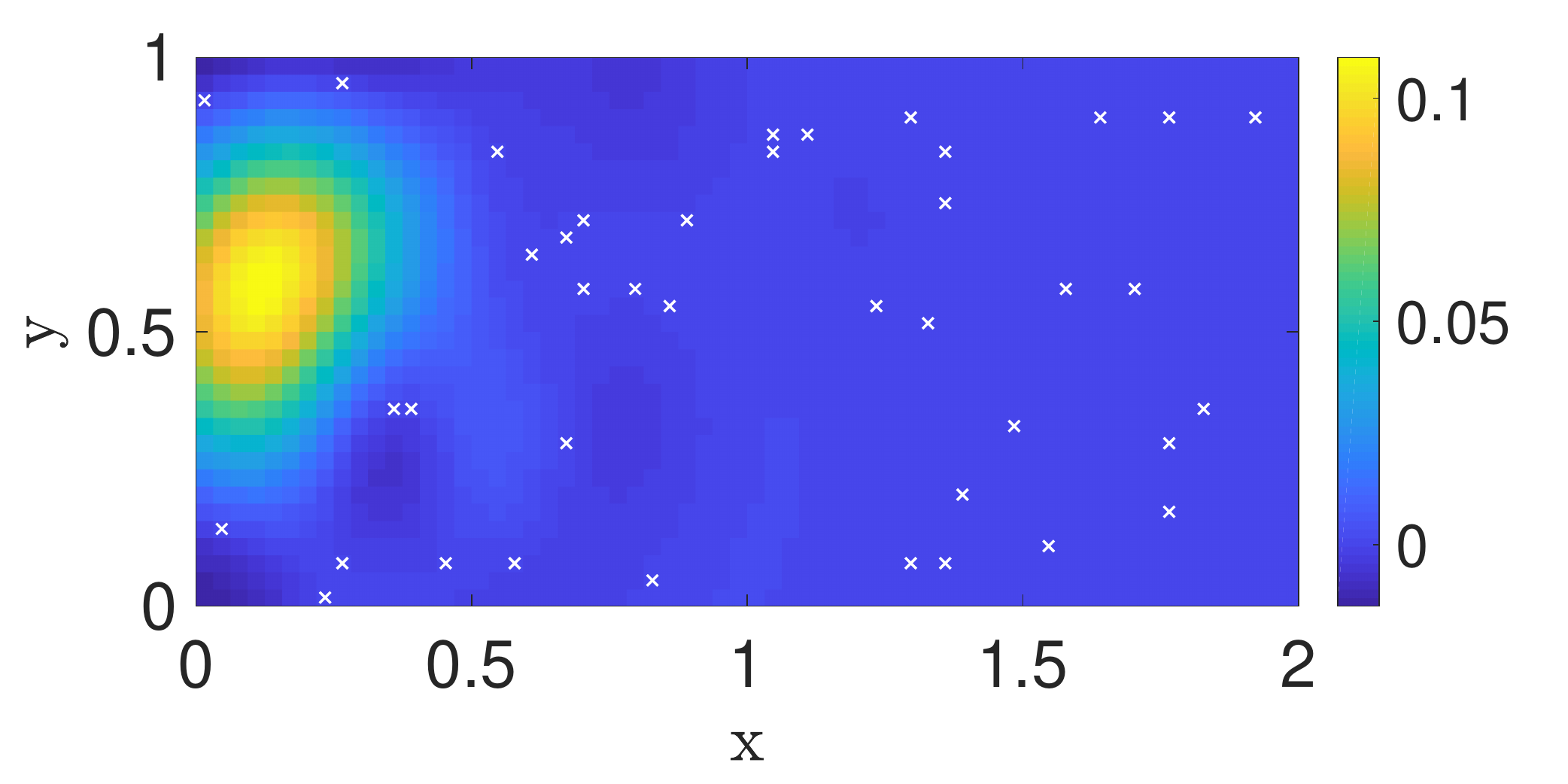}
    \caption{$g^c(\mathbf{x},\omega)$, Approach 1}
    \label{RK:fig:u_xi_mc_eigen_cond_eig_func_approach1_1}
  \end{subfigure}
  \begin{subfigure}[t]{0.48\textwidth}
    \centering%
    \includegraphics[scale=.32]{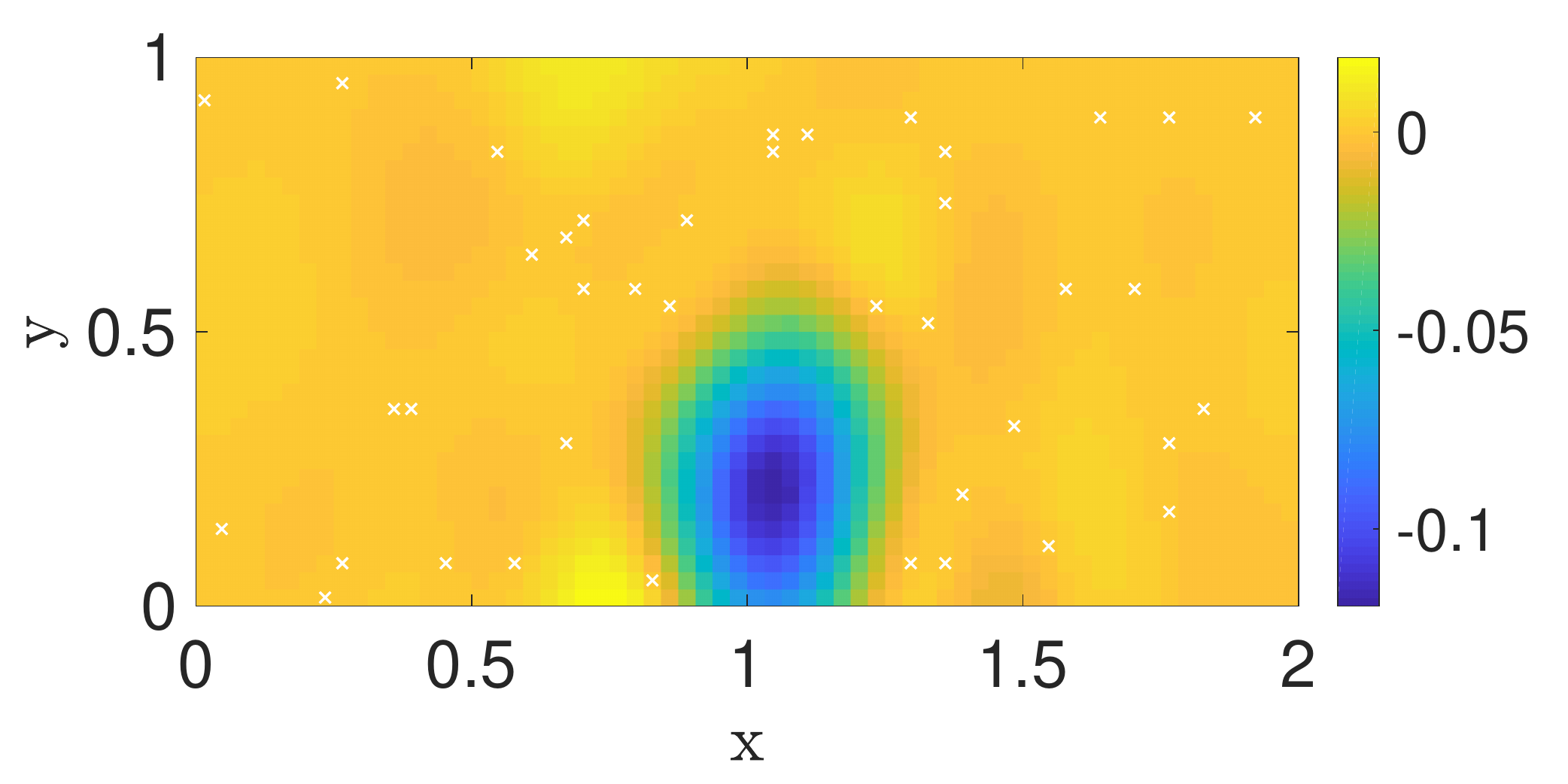}
    \caption{$g^c(\mathbf{x},\omega)$, Approach 1}
    \label{RK:fig:u_xi_mc_eigen_cond_eig_func_approach1_2}
  \end{subfigure}
  \begin{subfigure}[t]{0.48\textwidth}
    \centering%
    \includegraphics[scale=.32]{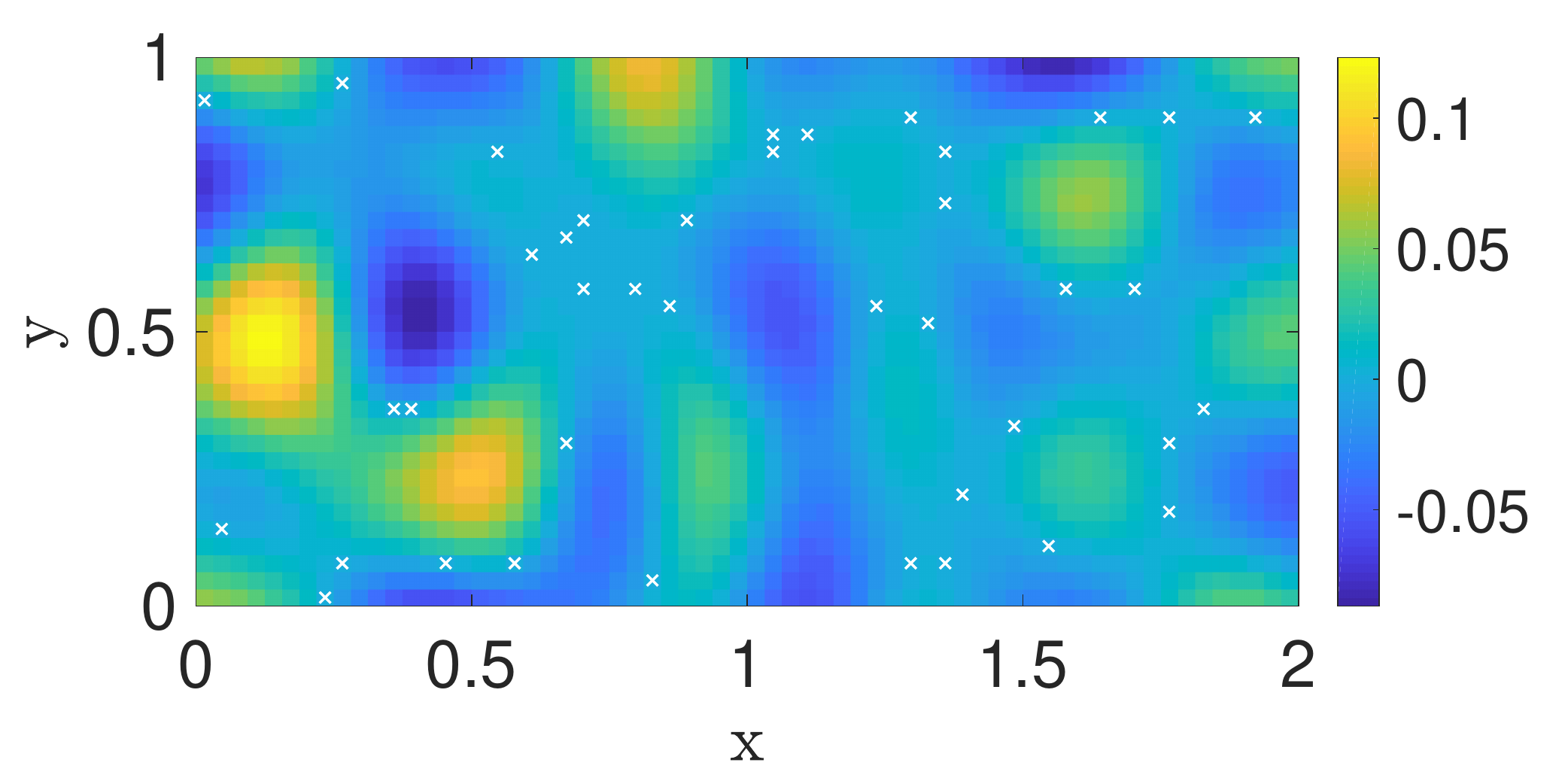}
    \caption{$g^c(\mathbf{x},\omega)$, Approach 2}
    \label{RK:fig:u_xi_mc_eigen_cond_eig_func_1}
  \end{subfigure}
  \begin{subfigure}[t]{0.48\textwidth}
    \centering%
    \includegraphics[scale=.32]{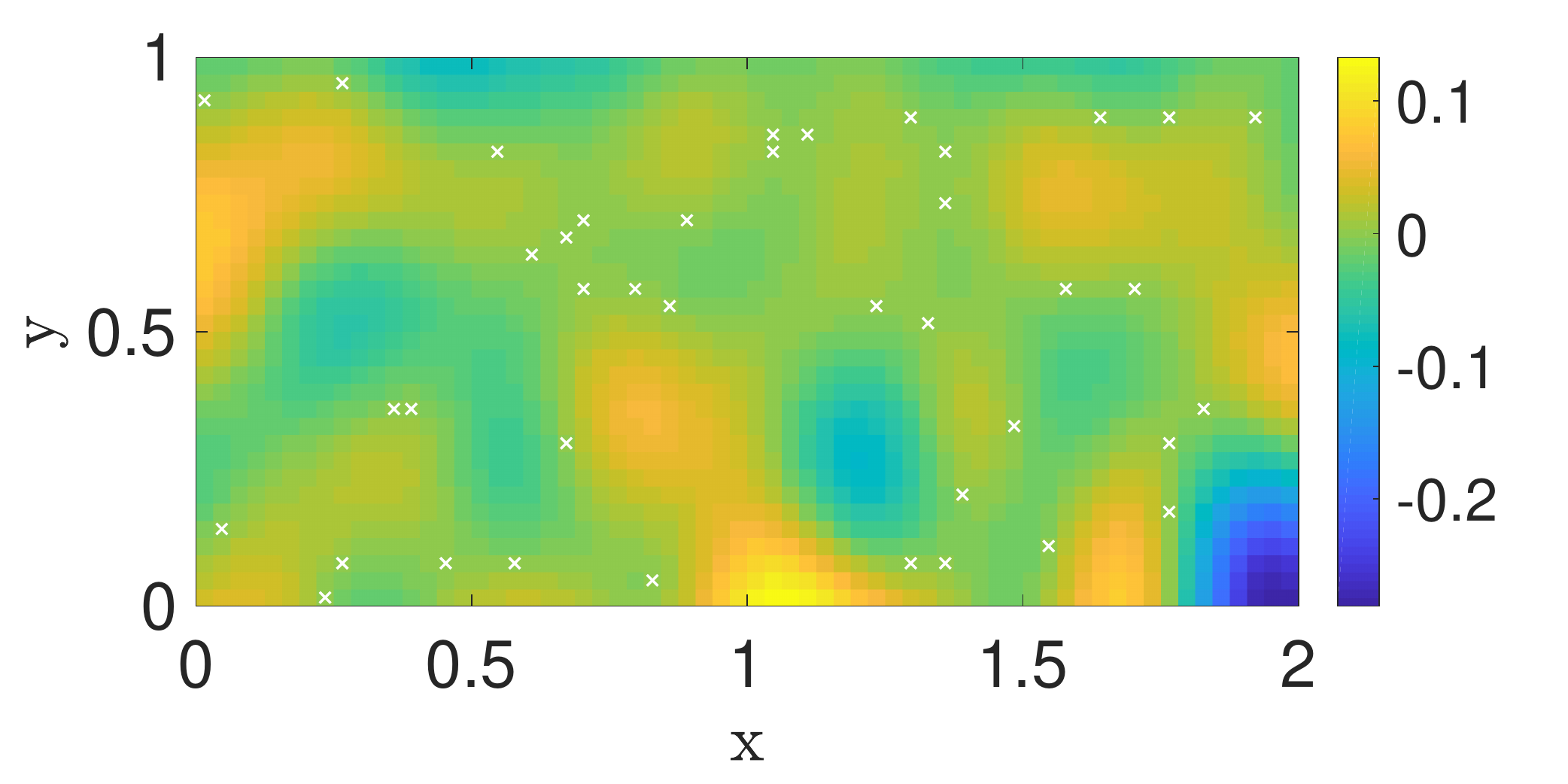}
    \caption{$g^c(\mathbf{x},\omega)$, Approach 2}
    \label{RK:fig:u_xi_mc_eigen_cond_eig_func_2}
  \end{subfigure}
  \caption{First (left) and second (right) eigenfunctions of $g(\mathbf{x},\omega)$ and $g^c(\mathbf{x},\omega)$}
  \label{fig:u_xi_mc_eigen_cond_eig_func}
\end{figure}

Figure~\ref{fig:u_xi_mc_eigen_cond_eig_func} shows the first and second eigenfunctions of $g^c(x,\omega)$ resulting from Approaches 1 and 2, together with the eigenvalues of $g(x,\omega)$.
Here, we note that, by construction, the eigenfunctions of Approach 2 are calculated from the first 20 eigenfunctions of $g(\mathbf{x}, \omega)$.
In contrast, in Approach 1, all the eigenfunctions of $g(\mathbf{x}, \omega)$ contribute to the first 20 eigenfunctions of $g^c(\mathbf{x}, \omega)$; therefore, the eigenfunctions of Approach 1 can resolve finer-scale features than the eigenfunctions of Approach 2.
We attribute the superior accuracy of Approach 1 for approximating conditional moments to its superior capacity for resolving fine-scale features of $g^c(\mathbf{x}, \omega)$.

\subsection{Active learning}
\label{sec:active-learning}

Here, we apply the two active learning data acquisition methods presented in Section~\ref{sec:active_learning} to identify additional measurement locations for the $g$ field.
We use these additional observations together with the previously available observations to construct the conditional KL expansion of $g$ using Approach 1, as it was shown in Section~\ref{sec:numerical-approach-1} to be more accurate than Approach 2 for the considered application.

\begin{figure}[ht!]
    \centering
    \begin{subfigure}[t]{0.48\textwidth}
        \centering
        \includegraphics[scale=.28]{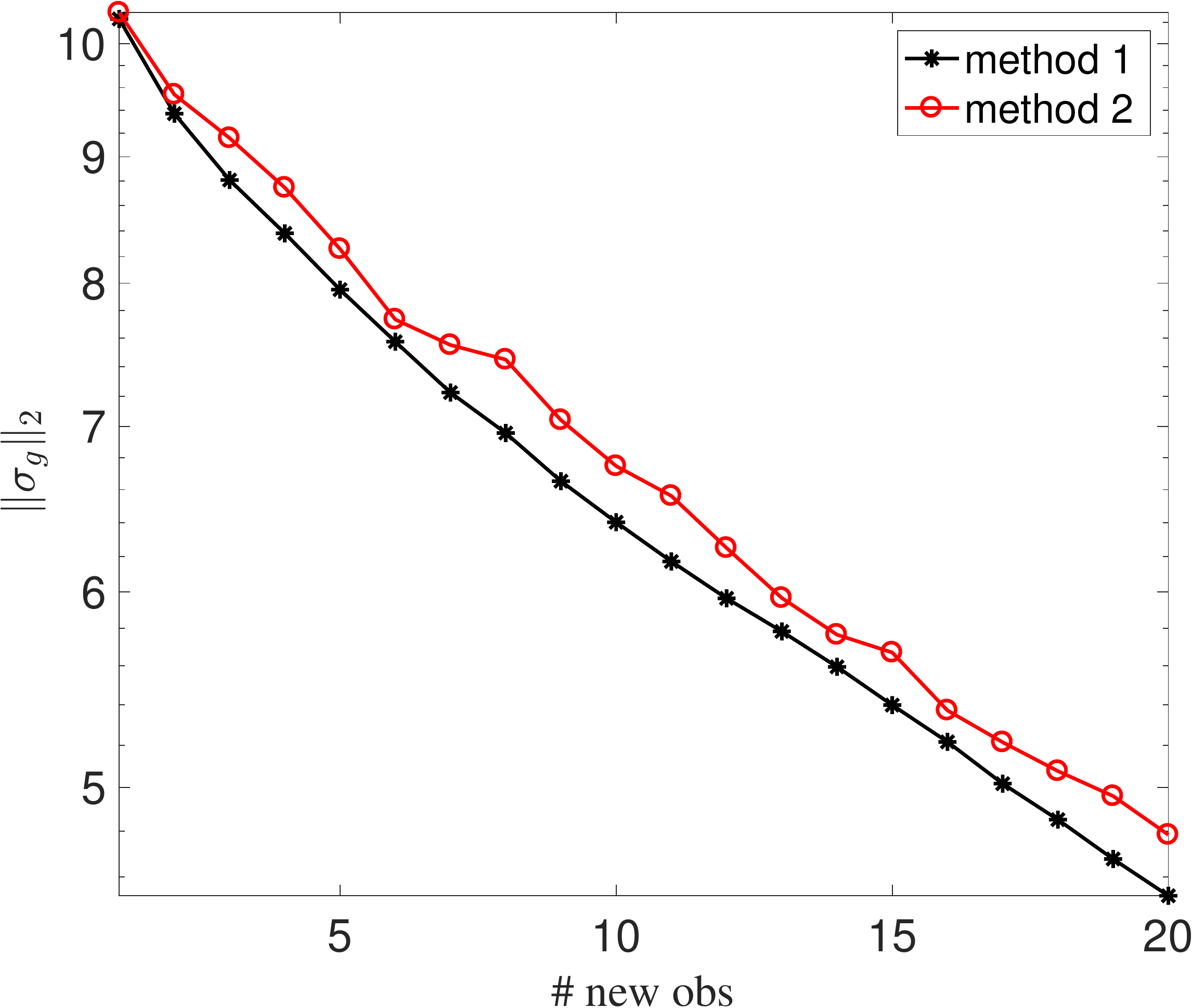}
        \caption{$|| \sigma_{g^c}||_2$} \label{RK:fig:norm_g_sdv_k_and_k_u}
    \end{subfigure}        
    \begin{subfigure}[t]{0.48\textwidth}
        \centering
        \includegraphics[scale=.28]{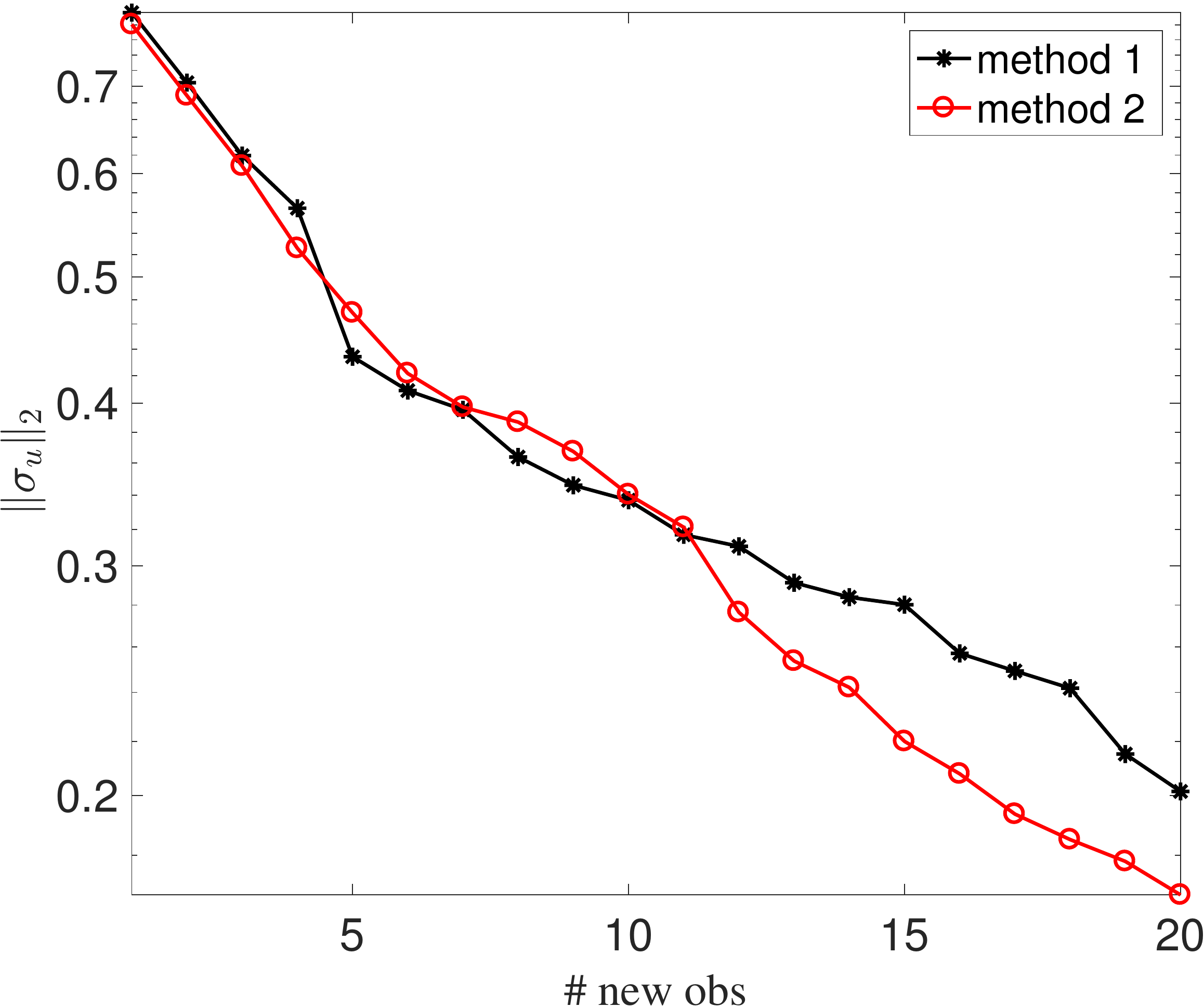}
        \caption{$|| \sigma_{u^c}||_2$} \label{RK:fig:norm_u_sdv_k_and_k_u}
    \end{subfigure}    
       \caption{$L_2$ norm of standard deviation of (a) $g^c(\mathbf{x},\omega)$ and (b) $u(\mathbf{x},\omega)$ as a function of the number of additional measurements. Standard deviation of unconditional $g(\mathbf{x},\omega)$ is $\sigma_g = 0.65$.}   \label{norm_g_u_sdv_k_and_k_u}
\end{figure}

Here, we consider again the SPDE problem~\eqref{RK:eq:spde_nr}, and we aim to explore the behavior of the presented data acquisition methods for two choices of $\sigma_g$, namely $0.65$ and $1.3$.
For $\sigma_g = 0.65$, Figure~\ref{norm_g_u_sdv_k_and_k_u} shows the $L^2$-norm of the standard deviation of $g^c$ and $u^c$ as a function of the number of additional measurements, $N_{am}$ identified with both active learning methods. 
In Method 2, the conditional covariances $\hat{C}^c_{u}$ and $\hat{C}^c_{ug}$ are computed using \num{200} MC realizations.
This small number of MC realizations is justified as empirical observation shows that very accurate estimates of the covariances are not necessary for obtaining good estimates of the new observation locations from the minimization problem~\eqref{minimization}.

Figure~\ref{RK:fig:norm_g_sdv_k_and_k_u} shows that the decay rate of the norm of the $g^c$ standard deviation is approximately the same for both methods. As expected, the variance reduction in Method 1 is larger (by approximately $5\%$) than in Method 2, because by construction, Method 1 reduces the variance of $g^c$ optimally.
On the other hand, Figure~\ref{RK:fig:norm_u_sdv_k_and_k_u} shows that Method 2 reduces the $L^2$-norm of $\sigma^c_u$ more than Method 1 for most $N_{am}$, and by more than $15\%$ for $N_{am}>11$.
For some $N_{am}<11$, Method 2 has larger norm of $\sigma^c_u$ than Method 1, which we attribute to the error in the $\sigma^c_u$ approximation \eqref{eq:GP_cov_u}. 

\begin{figure}[ht!]
    \centering
    \begin{subfigure}[t]{0.48\textwidth}
        \centering
        \includegraphics[scale=.28]{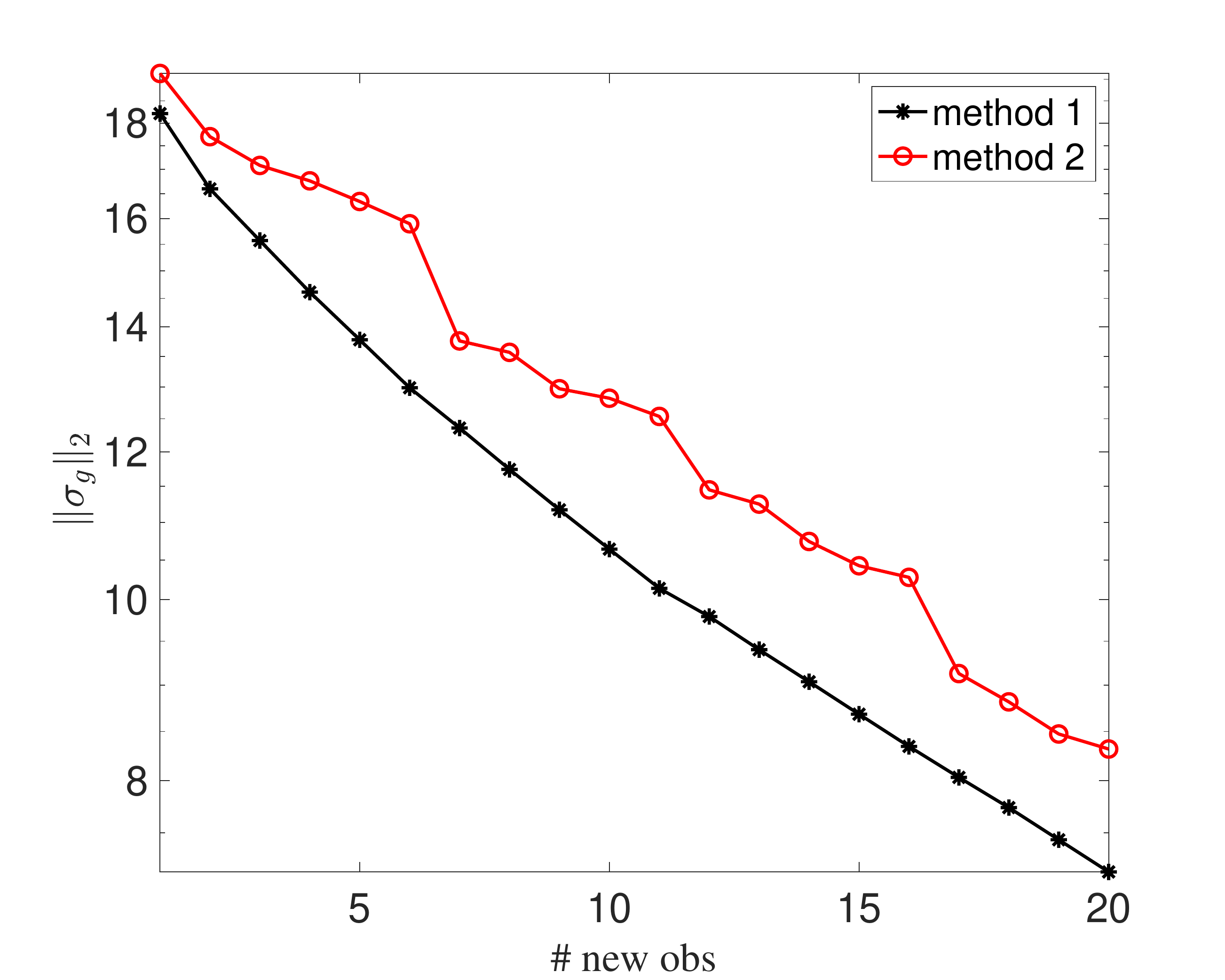}
        \caption{$|| \sigma_{g^c}||_2$} \label{RK:fig:norm_g_sdv_k_and_k_u_sdev_1_03}
    \end{subfigure}        
    \begin{subfigure}[t]{0.48\textwidth}
        \centering
        \includegraphics[scale=.28]{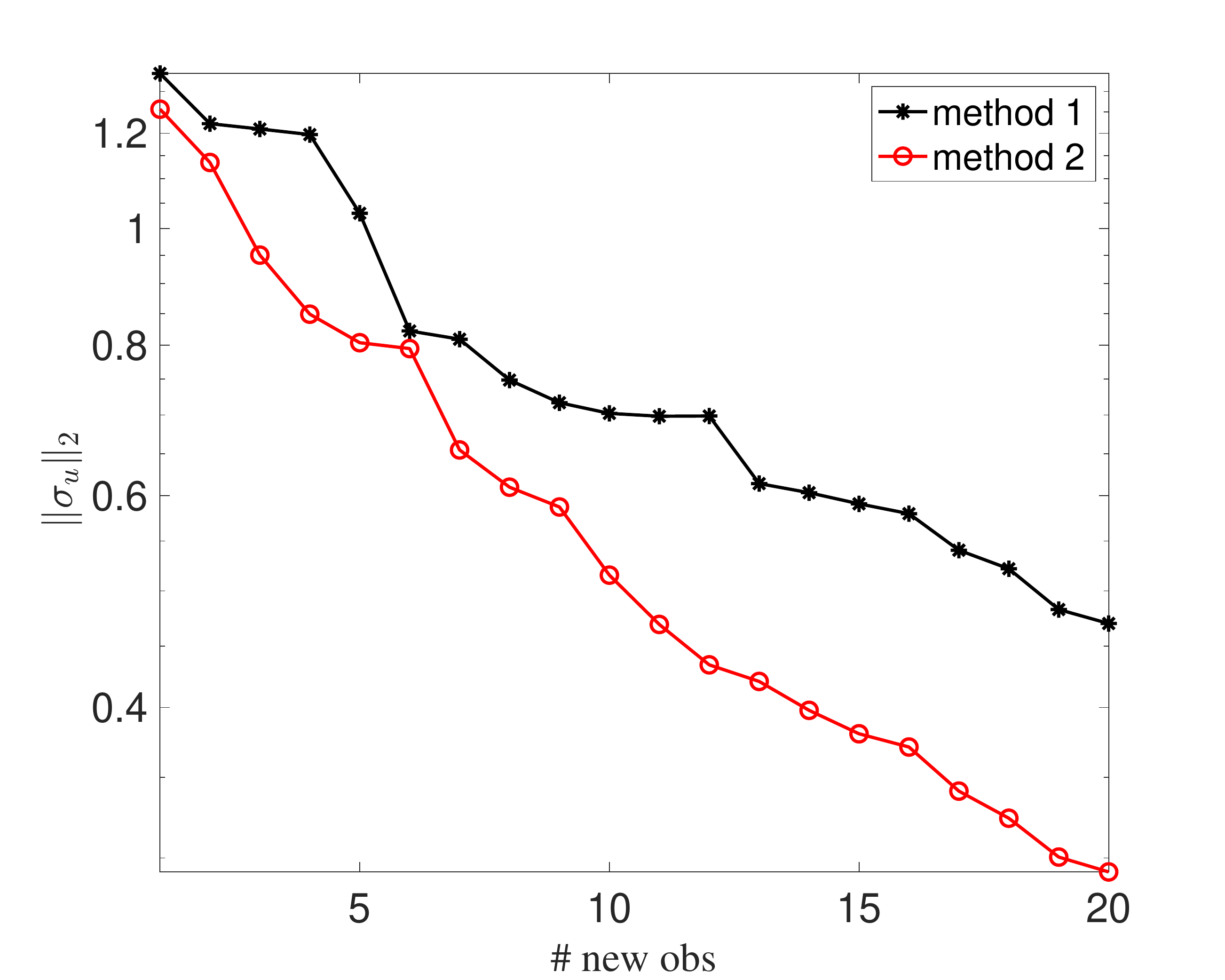}
        \caption{$|| \sigma_{u^c}||_2$} \label{RK:fig:norm_u_sdv_k_and_k_u_sdev_1_03}
    \end{subfigure}    
       \caption{$L_2$ norm of standard deviation of (a) $g^c(\mathbf{x},\omega)$ and (b) $u(\mathbf{x},\omega)$ as a function of the number of additional measurements. Standard deviation of unconditional $g(\mathbf{x},\omega)$ is $\sigma_g = 1.3$.}   \label{norm_g_u_sdv_k_and_k_u_sdev_1_03}
\end{figure}

We obtain qualitatively similar (but more pronounced) results for $\sigma_g = 1.3$, as shown in Figure~\ref{norm_g_u_sdv_k_and_k_u_sdev_1_03}.
Method 1 leads to a sharper decrease of the $\sigma^c_g$ norm and  Method 2 results in a sharper decrease of the $\sigma^c_u$ norm than what we observe for  $\sigma_g = 0.65$.
From Figures~\ref{norm_g_u_sdv_k_and_k_u} and \ref{norm_g_u_sdv_k_and_k_u_sdev_1_03}, we conclude that Method 2 is more efficient than Method 1 for reducing uncertainty in $u^c$, and that \eqref{eq:GP_cov_u} provides a sufficiently accurate approximation for solving the minimization problem \eqref{eq:active-learning-criteria-u}.

\begin{figure}[ht!]
    \centering
    \begin{subfigure}[t]{0.48\textwidth}
        \centering
        \includegraphics[scale=.3]{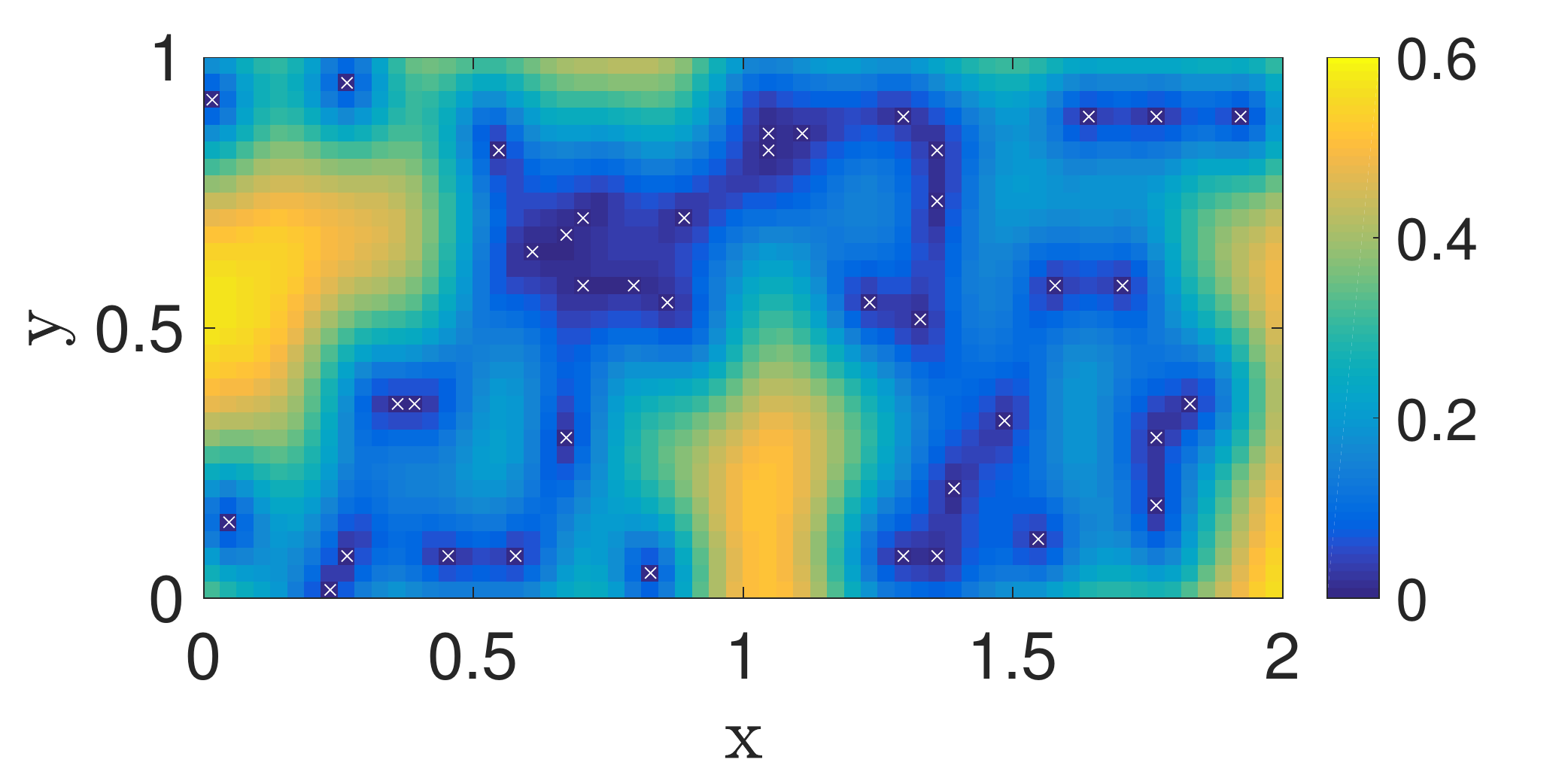}
        \caption{$\sigma_{g^c}$, $N_{am}=0$} \label{RK:fig:C_pg_est_exact_nobs40_mrst_al}
    \end{subfigure}        \\
    \begin{subfigure}[t]{0.48\textwidth}
        \centering
        \includegraphics[scale=.3]{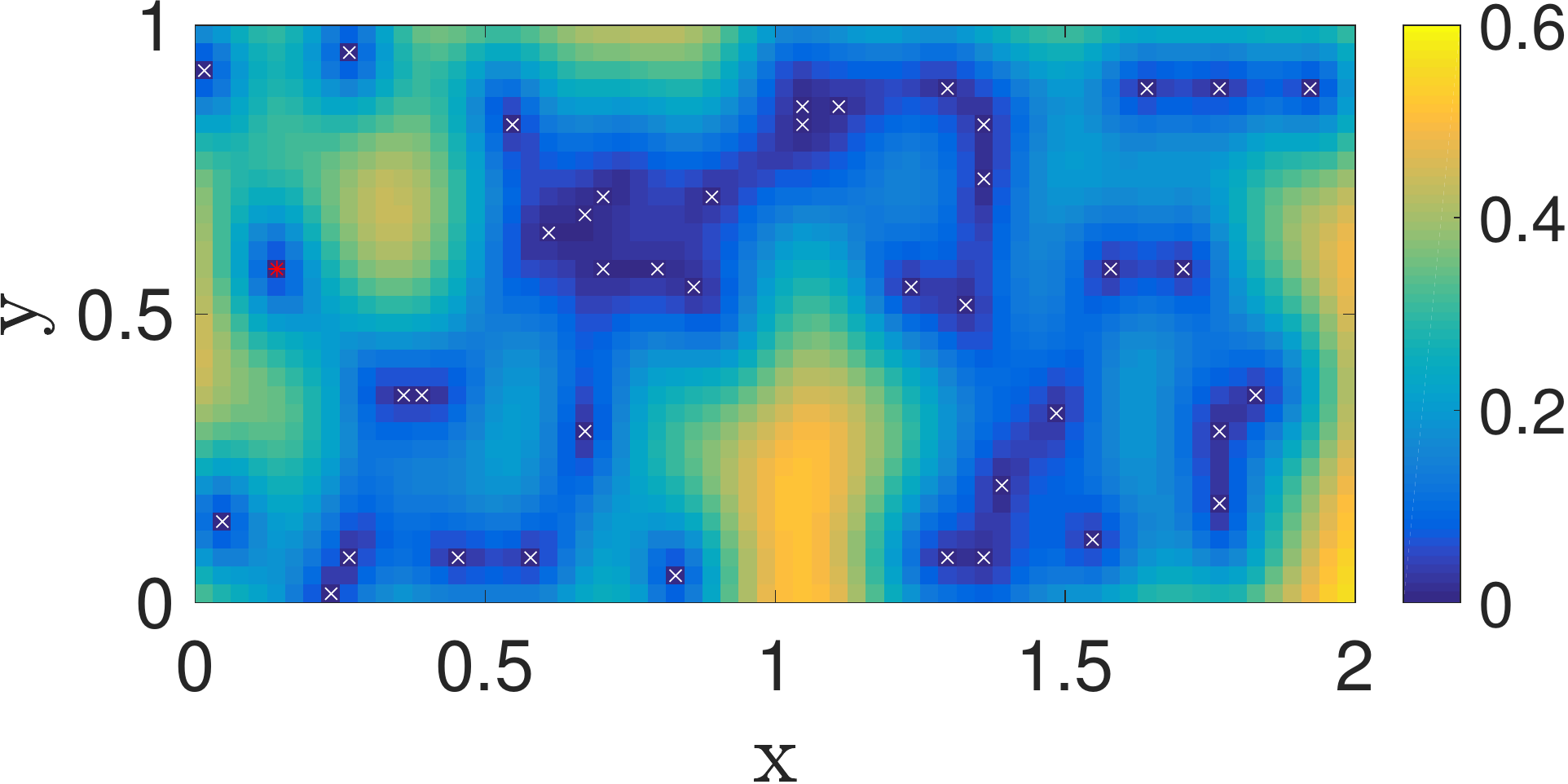}
        \caption{$N_{am}=1$, Method 1} \label{RK:fig:C_pg_exact_nobs40_mrst_x_star_1}
    \end{subfigure}    
     \begin{subfigure}[t]{0.48\textwidth}
        \centering
        \includegraphics[scale=.3]{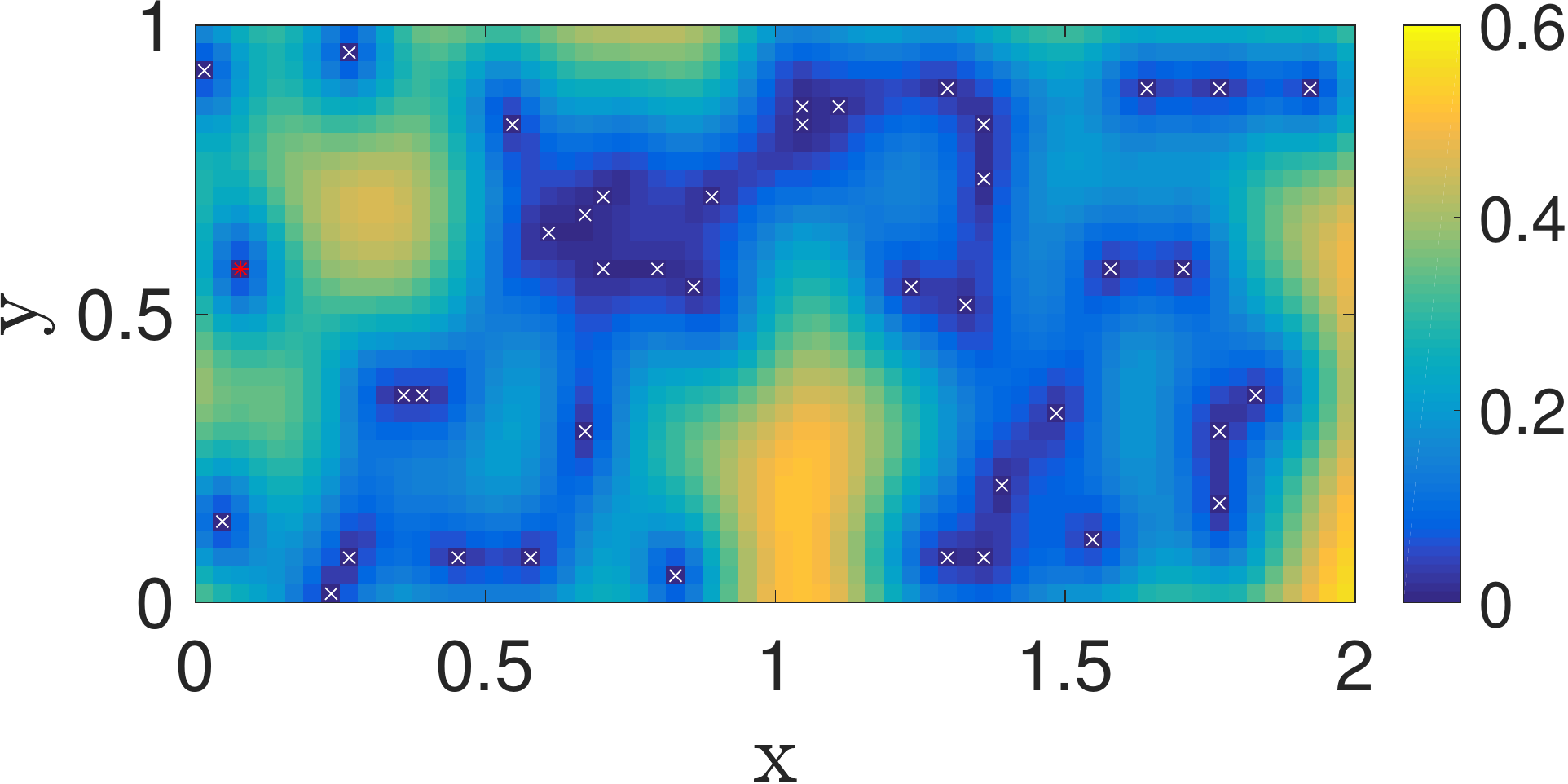}
        \caption{$N_{am}=1$, Method 2} \label{RK:fig:C_pg_exact_nobs40_mrst_x_star_k_u_1}
    \end{subfigure}    
    \begin{subfigure}[t]{0.48\textwidth}
        \centering
        \includegraphics[scale=.3]{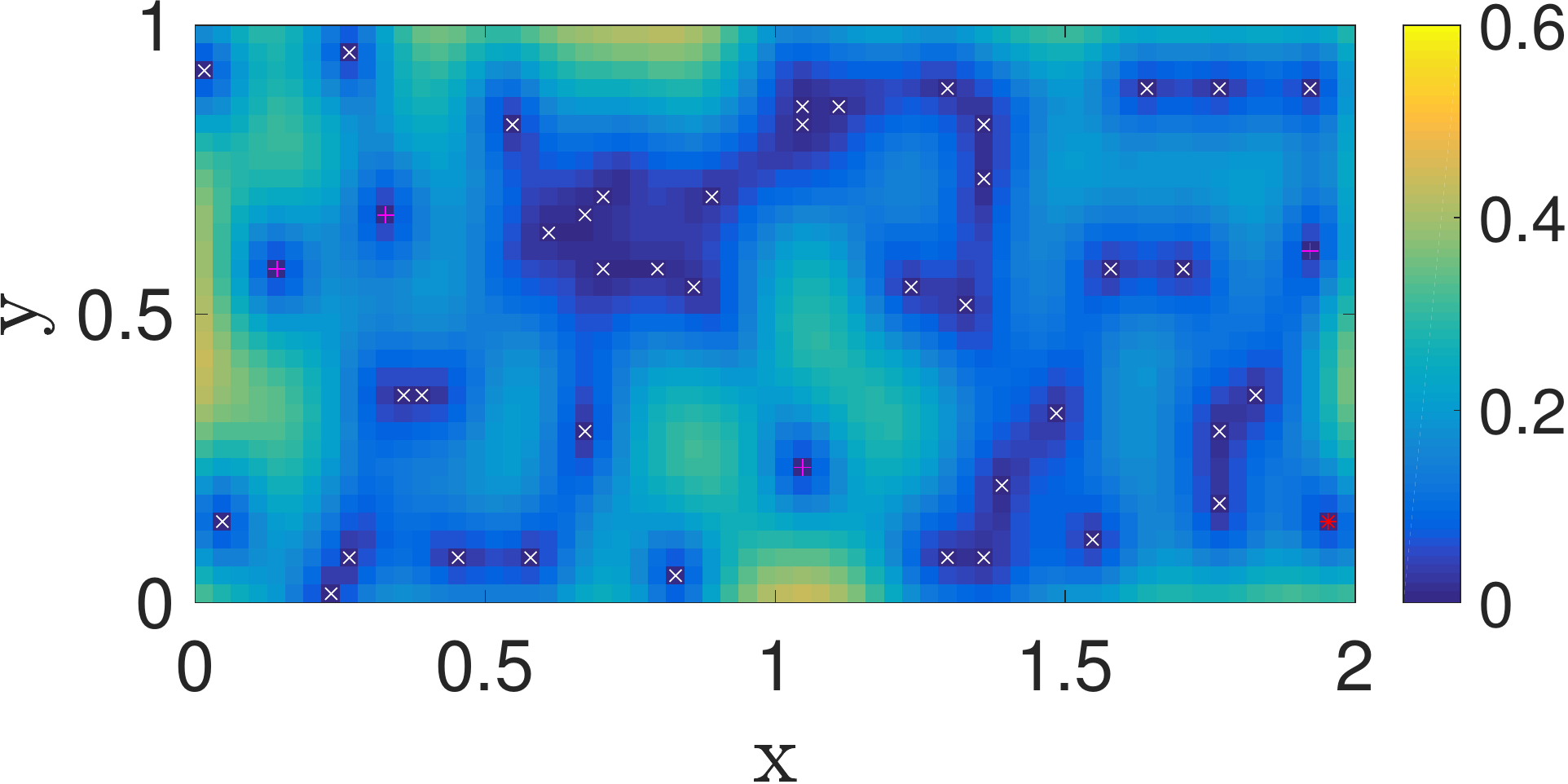}
        \caption{$N_{am}=5$, Method 1} \label{RK:fig:C_pg_exact_nobs40_mrst_x_star_5}
    \end{subfigure}        
    \begin{subfigure}[t]{0.48\textwidth}
        \centering
        \includegraphics[scale=.3]{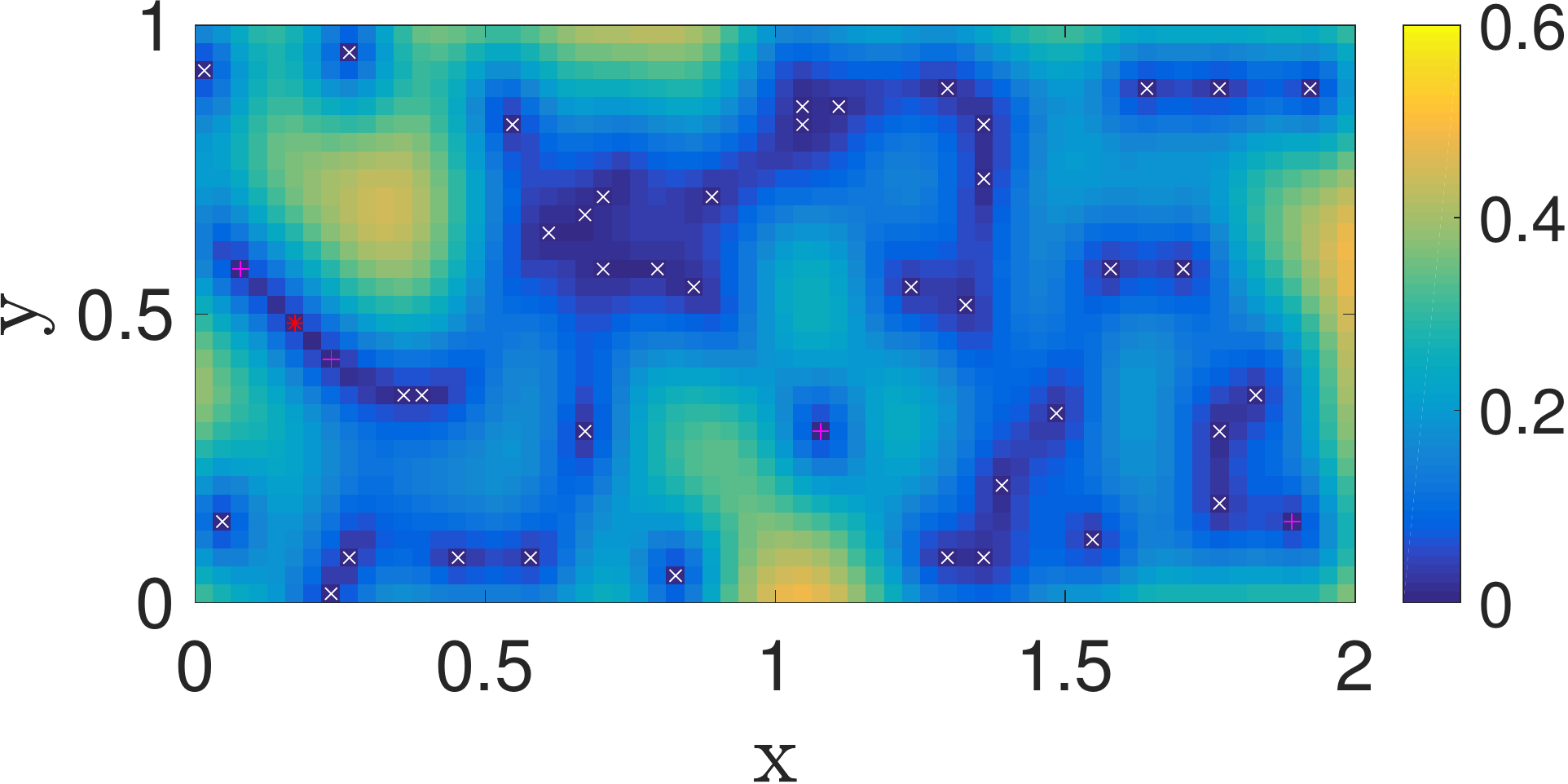}
        \caption{$N_{am}=5$, Method 2} \label{RK:fig:C_pg_exact_nobs40_mrst_x_star_k_u__5}
    \end{subfigure}  
    \begin{subfigure}[t]{0.48\textwidth}
        \centering
        \includegraphics[scale=.3]{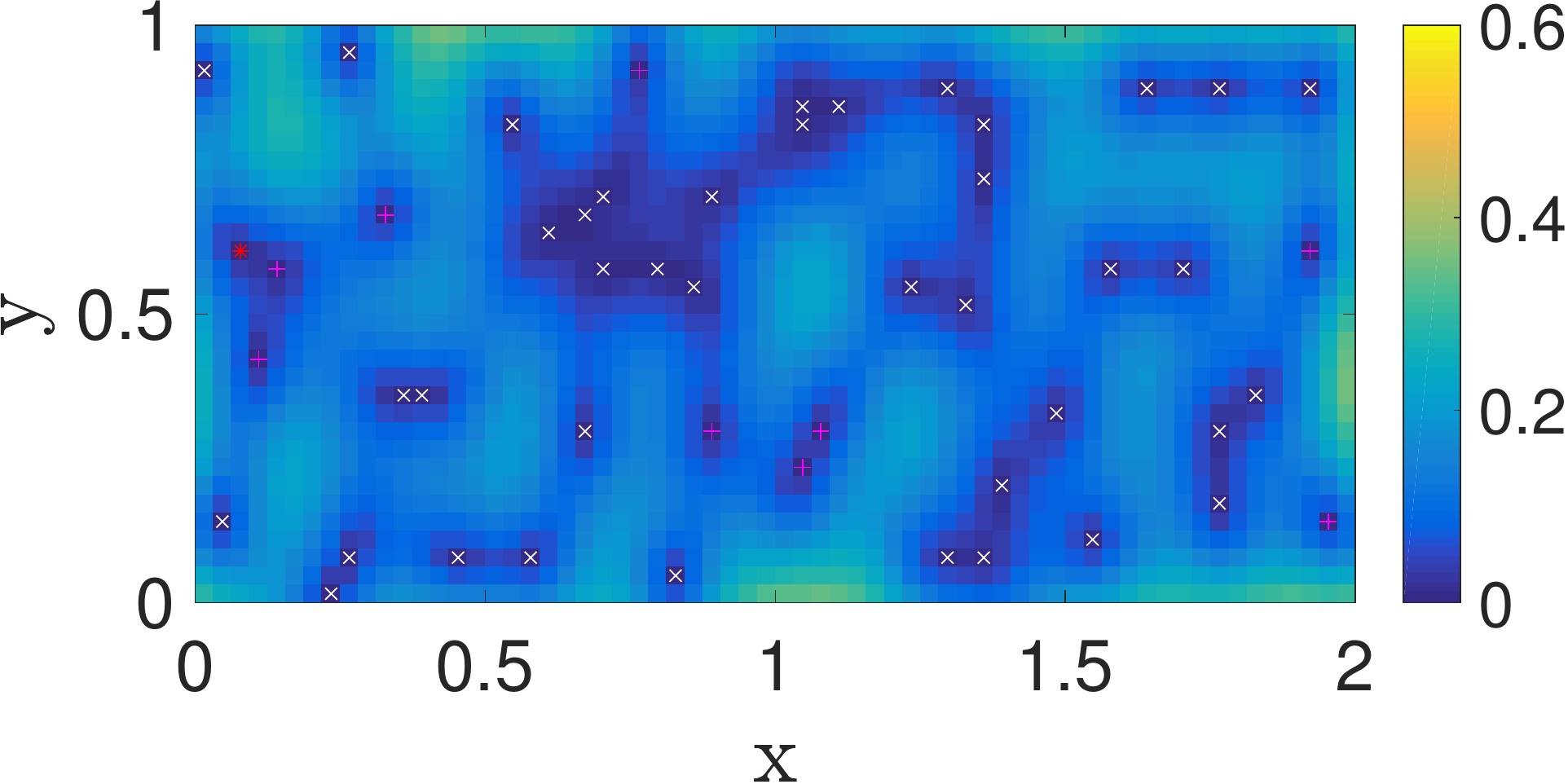}
        \caption{$N_{am}=10$, Method 1} \label{RK:fig:C_pg_exact_nobs40_mrst_x_star_10}
    \end{subfigure}    
    \begin{subfigure}[t]{0.48\textwidth}
        \centering
        \includegraphics[scale=.3]{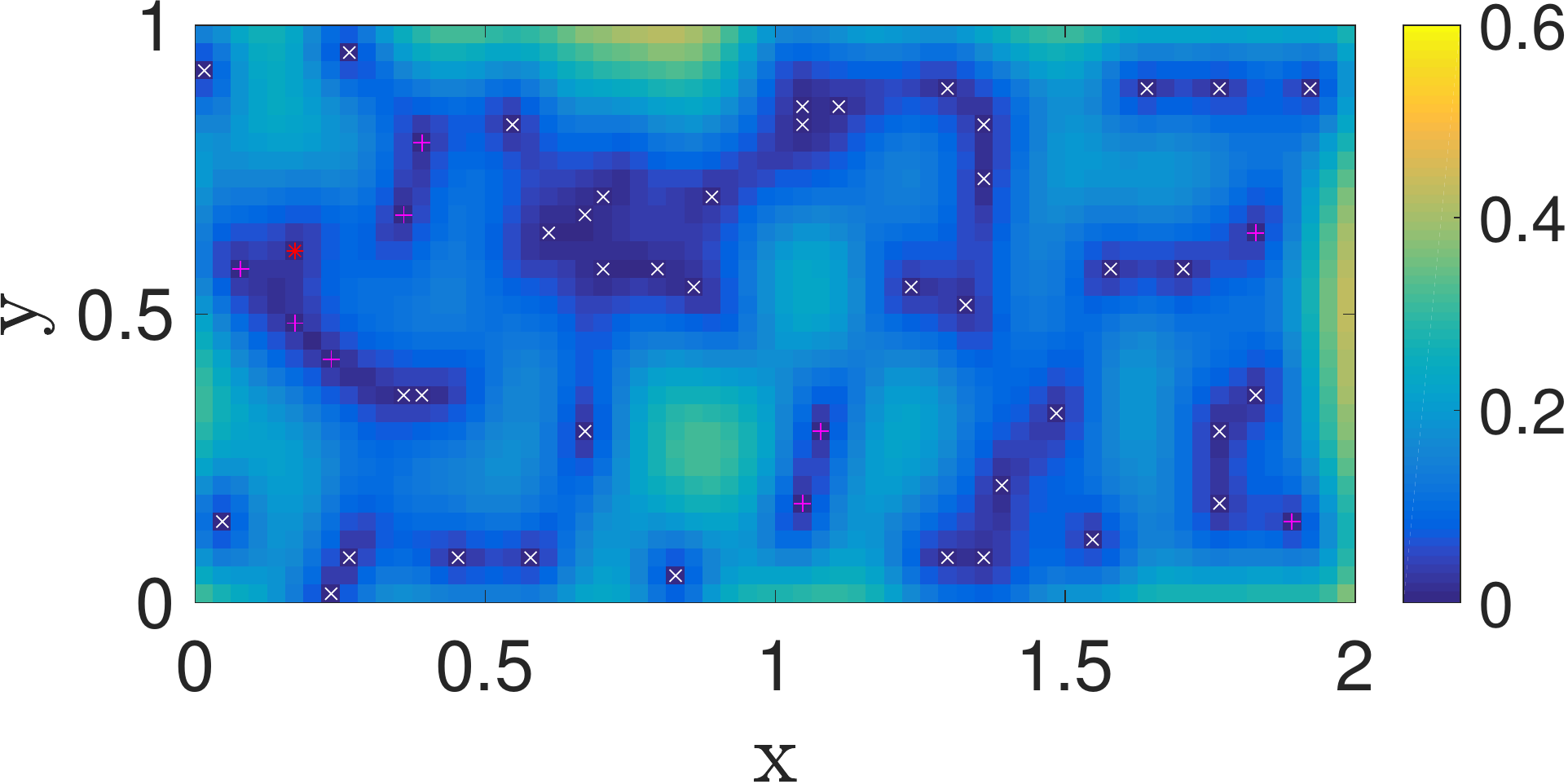}
        \caption{$N_{am}=10$, Method 2} \label{RK:fig:C_pg_exact_nobs40_mrst_x_star_k_u_10}
    \end{subfigure}    
    \caption{ (a) Standard deviation of $g^c(\mathbf{x},\omega)$ without additional measurements. Standard deviation of $g^c(\mathbf{x},\omega)$ with additional measurements: (b) $N_{am}=1$, Method 1; (c). $N_{am}=1$, Method 2; (d) $N_{am}=1$, Method 5; (e) $N_{am}=5$, Method 2; (f) $N_{am}=10$, Method 1; and (g) $N_{am}=10$, Method 2. Additional observation locations are shown in magenta and the original measurements are shown in red. Standard deviation of unconditional $g(\mathbf{x},\omega)$ is $\sigma_g = 0.65$.} \label{C_pg_exact_nobs40_mrst_x_star}
\end{figure}

\begin{figure}[ht!]
    \centering
    \begin{subfigure}[t]{0.48\textwidth}
        \centering
        \includegraphics[scale=.3]{figures/u_sdev_exact_cond_mc_mrst_no_obs_plot}
        \caption{$\sigma_{u^c}$, $N_{am}=0$} \label{RK:fig:u_sdev_exact_cond_mc_mrst_no_obs_plot_al}
    \end{subfigure}        \\
    \begin{subfigure}[t]{0.48\textwidth}
        \centering
        \includegraphics[scale=.3]{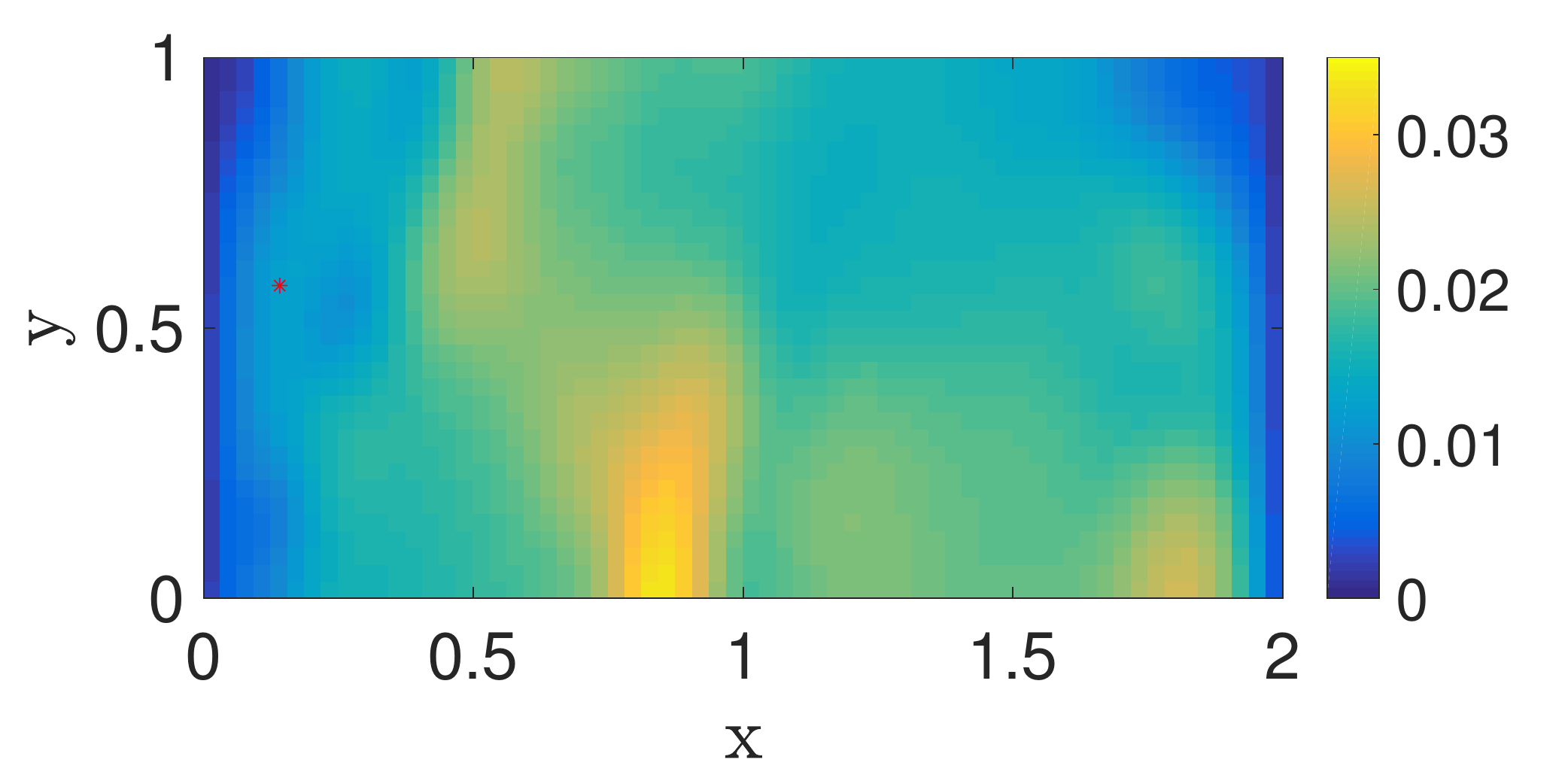}
        \caption{$N_{am}=1$, Method 1} \label{RK:fig:u_sdev_exact_nobs40_mrst_x_star_1}
    \end{subfigure}    
     \begin{subfigure}[t]{0.48\textwidth}
        \centering
        \includegraphics[scale=.3]{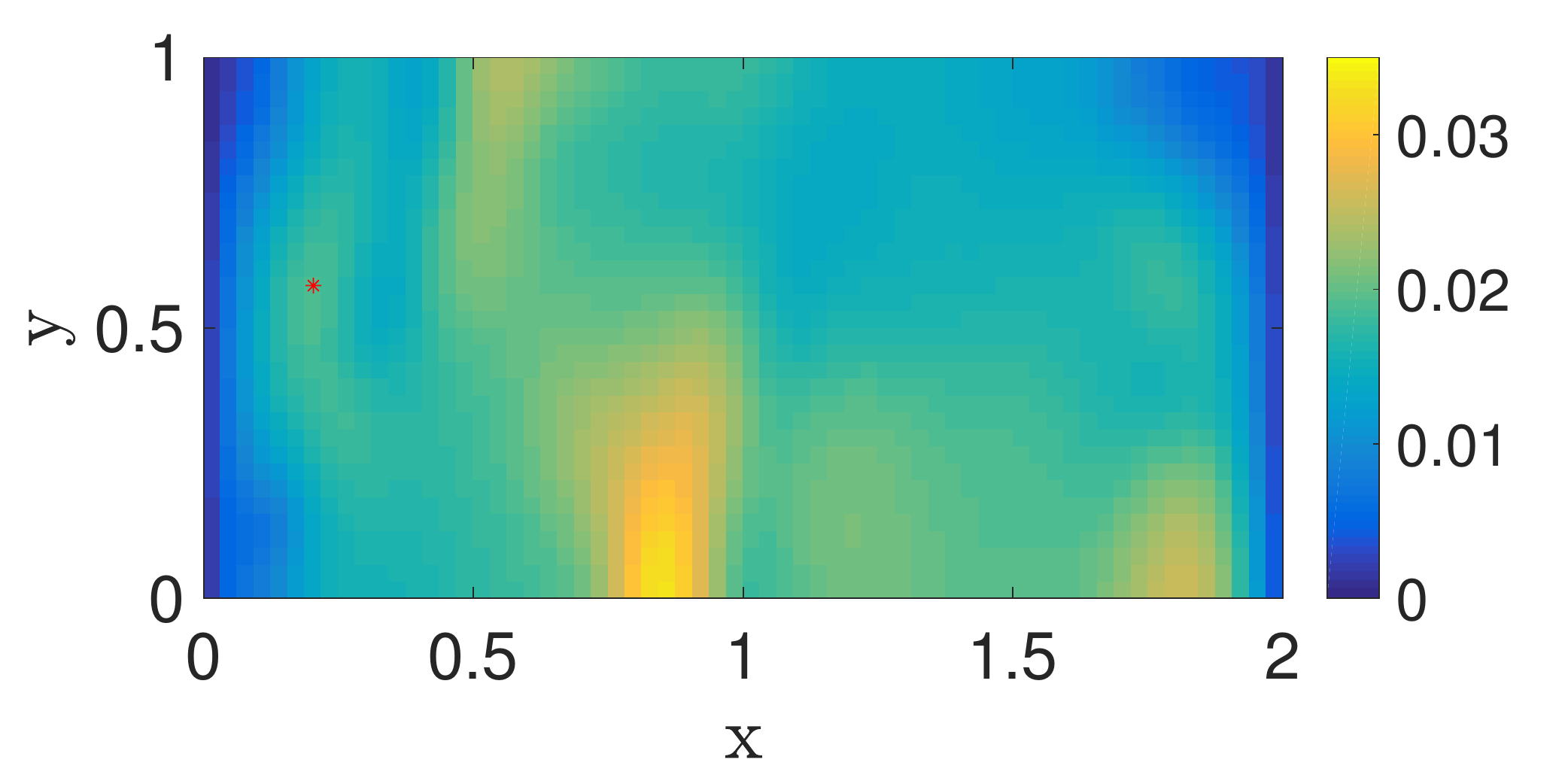}
        \caption{$N_{am}=1$, Method 2} \label{RK:fig:u_sdev_exact_nobs40_mrst_x_star_k_u_1}
    \end{subfigure}    
    \begin{subfigure}[t]{0.48\textwidth}
        \centering
        \includegraphics[scale=.3]{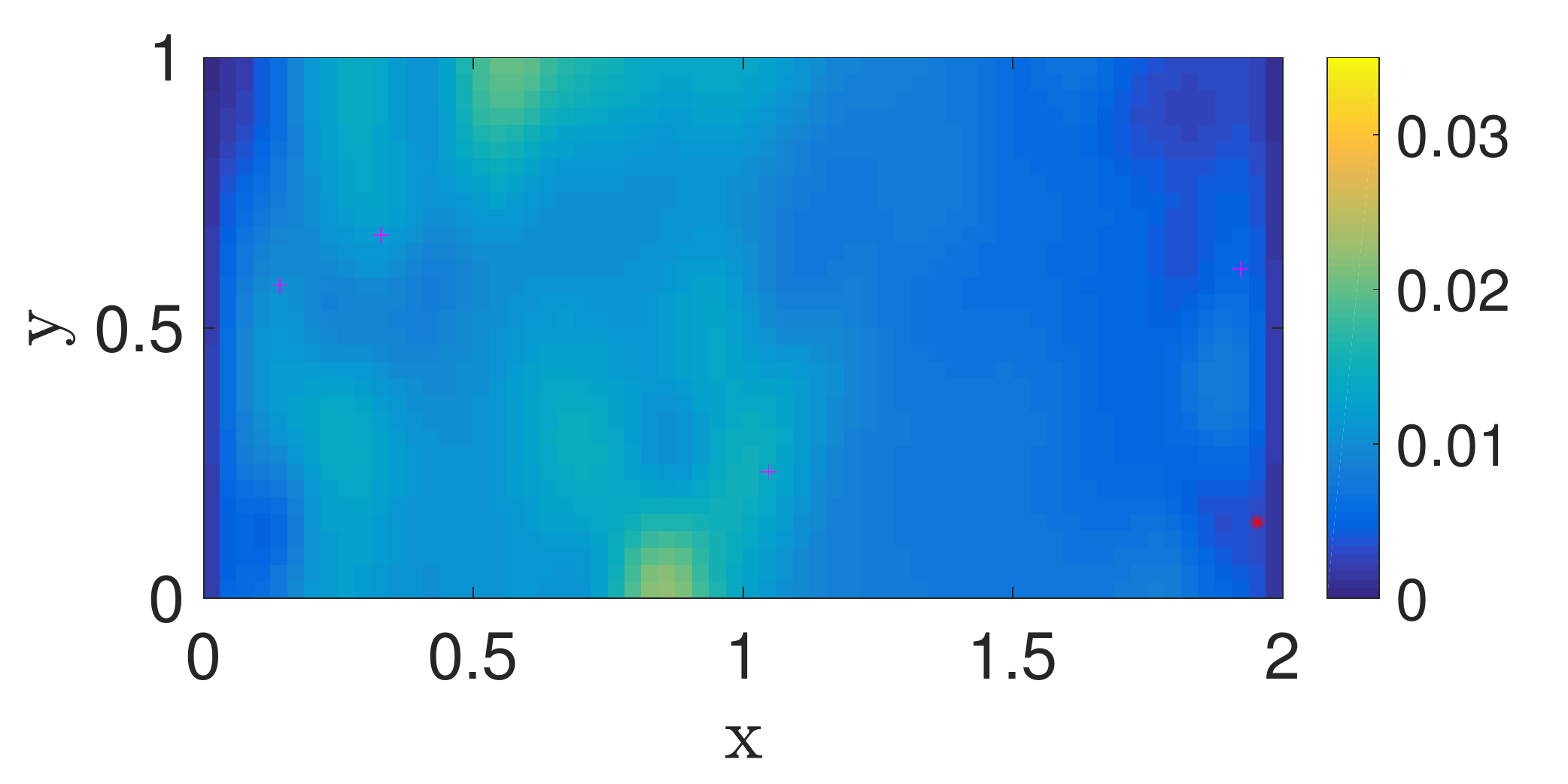}
        \caption{$N_{am}=5$, Method 1} \label{RK:fig:u_sdev_exact_nobs40_mrst_x_star_5}
    \end{subfigure}        
    \begin{subfigure}[t]{0.48\textwidth}
        \centering
        \includegraphics[scale=.3]{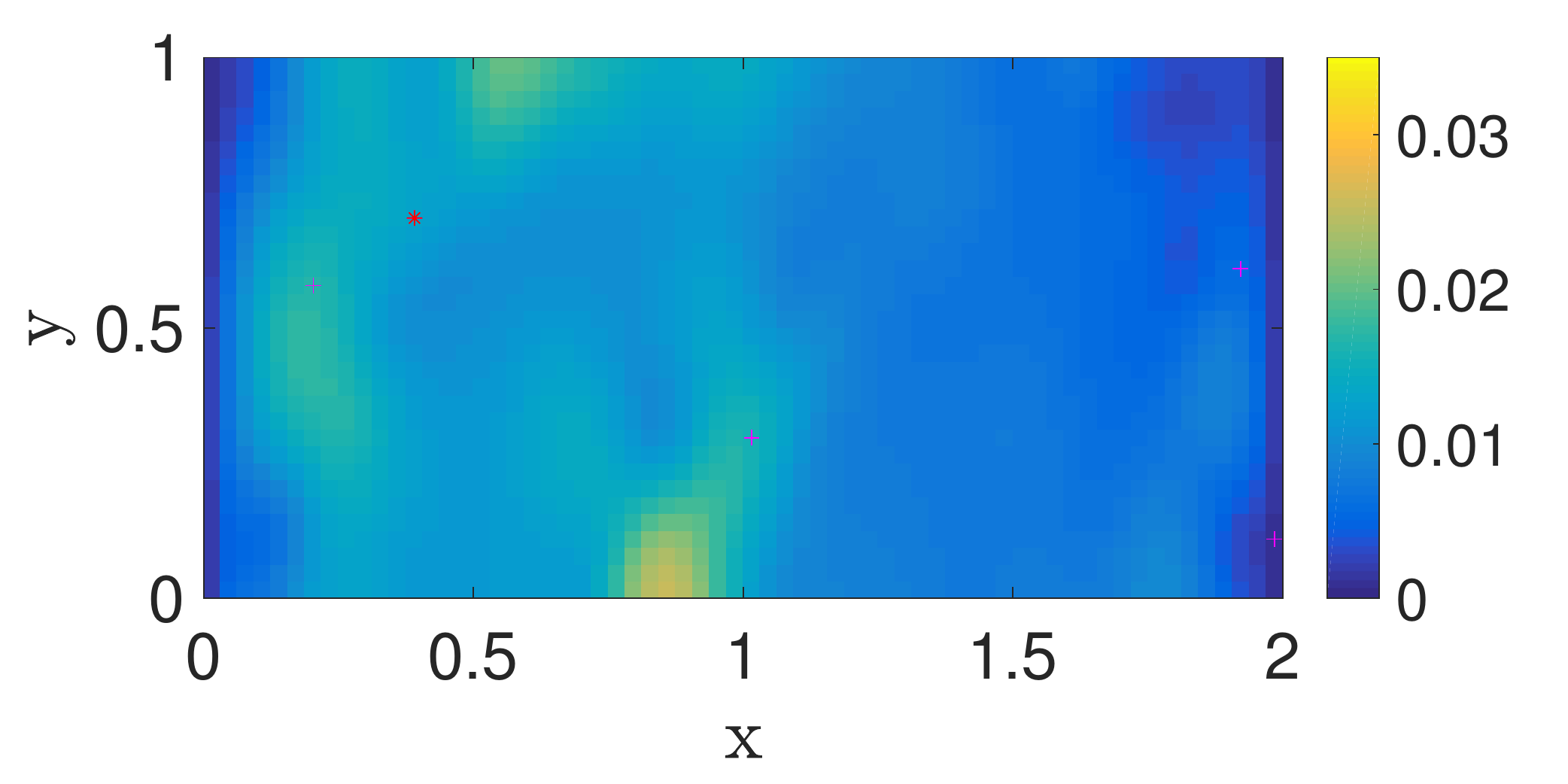}
        \caption{$N_{am}=5$, Method 2} \label{RK:fig:u_sdev_exact_nobs40_mrst_x_star_k_u__5}
    \end{subfigure}  
    \begin{subfigure}[t]{0.48\textwidth}
        \centering
        \includegraphics[scale=.3]{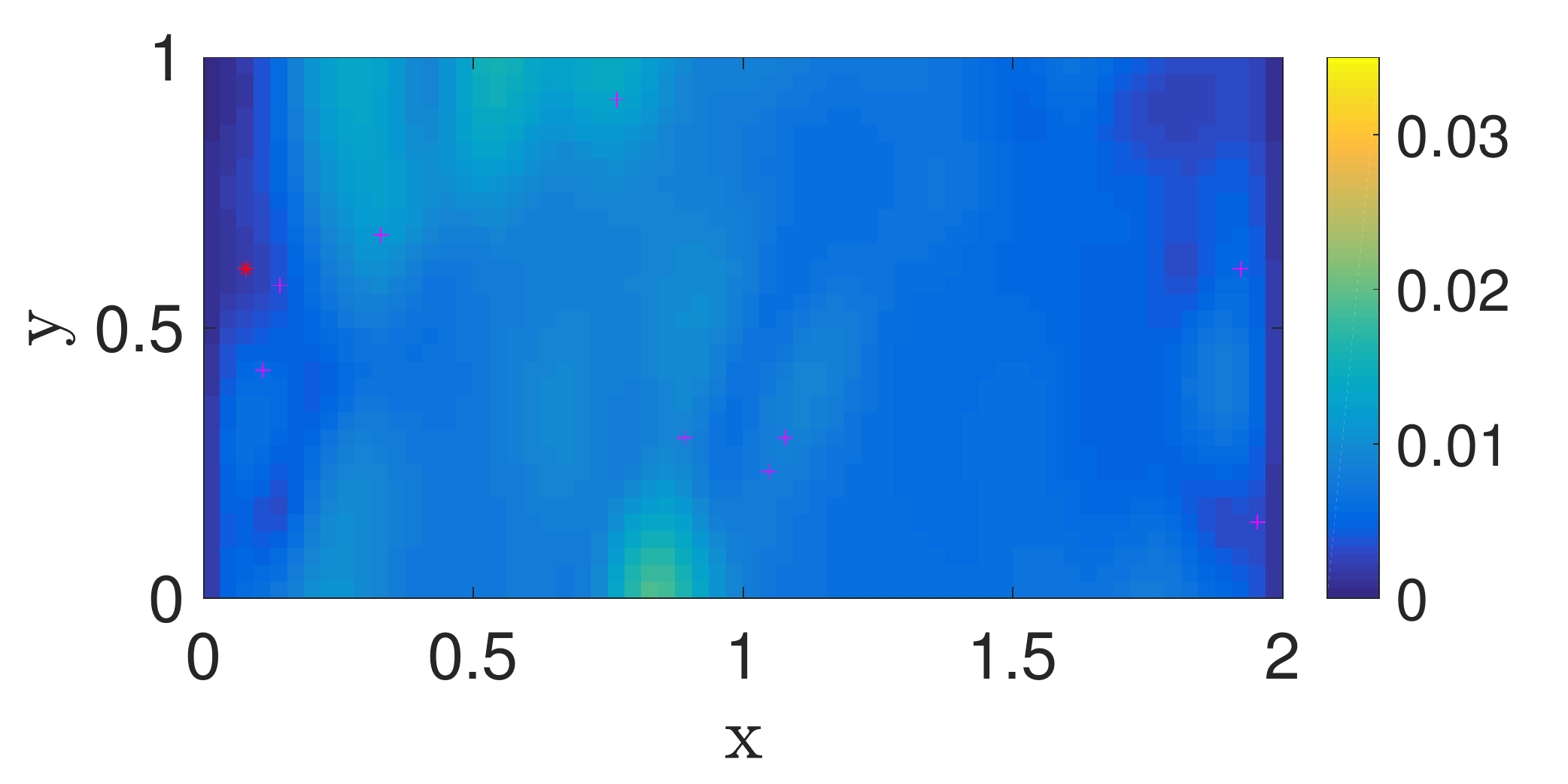}
        \caption{$N_{am}=10$, Method 1} \label{RK:fig:u_sdev_exact_nobs40_mrst_x_star_10}
    \end{subfigure}    
    \begin{subfigure}[t]{0.48\textwidth}
        \centering
        \includegraphics[scale=.3]{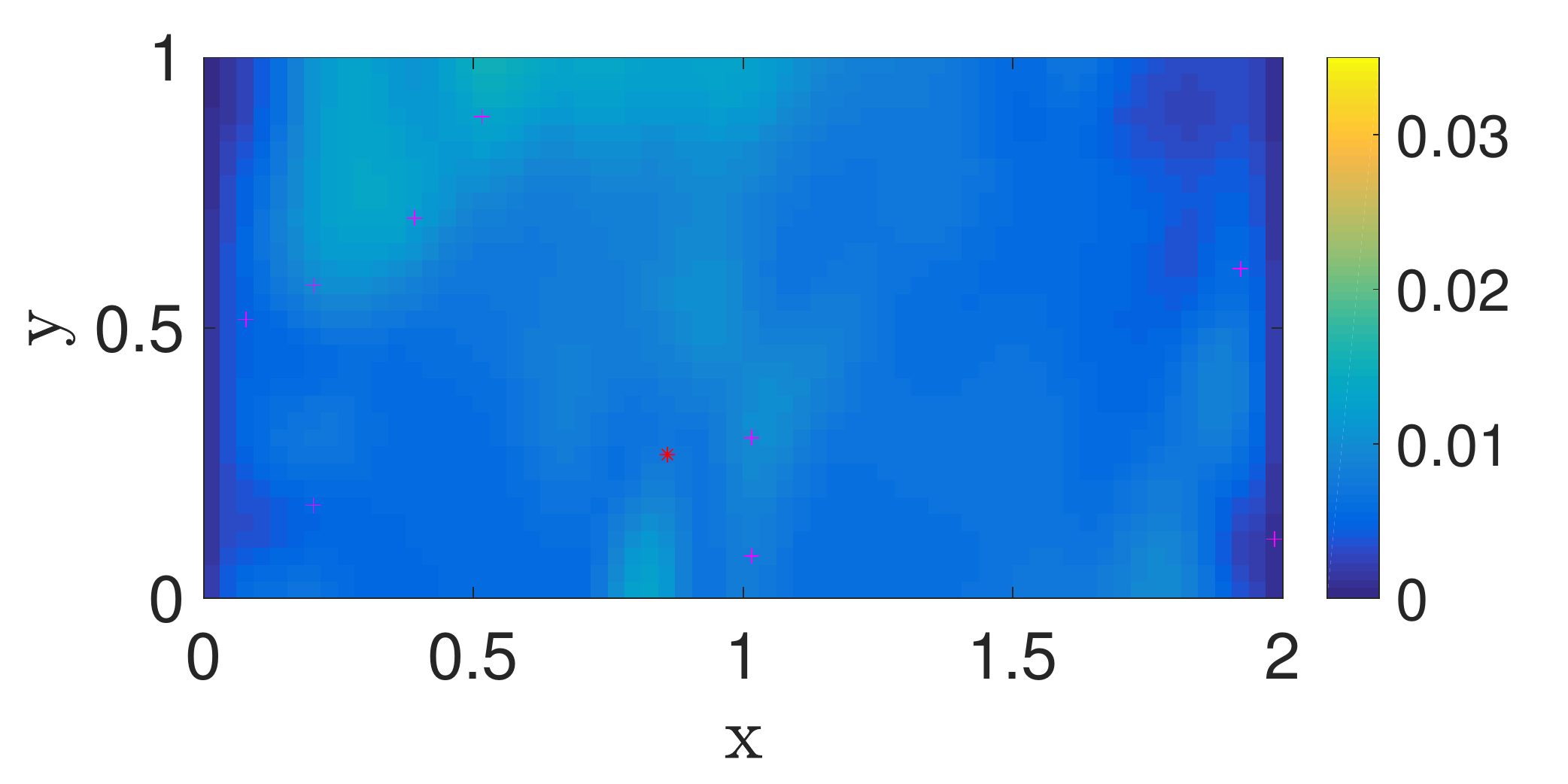}
        \caption{$N_{am}=10$, Method 2} \label{RK:fig:u_sdev_exact_nobs40_mrst_x_star_k_u_10}
    \end{subfigure}    
\caption{(a) Standard deviation of $u^c(\mathbf{x},\omega)$ without additional measurements. Standard deviation of $u^c(\mathbf{x},\omega)$ with additional measurements of $g$: (b) $N_{am}=1$, Method 1; (c). $N_{am}=1$, Method 2; (d) $N_{am}=1$, Method 5; (e) $N_{am}=5$, Method 2; (f) $N_{am}=10$, Method 1; and (g) $N_{am}=10$, Method 2. Additional observation locations of $g$ are shown in magenta and the original measurements are shown in red. Standard deviation of unconditional $g(\mathbf{x},\omega)$ is $\sigma_g = 0.65$.} \label{u_sdev_exact_nobs40_mrst_x_star_k_u}
\end{figure}

\begin{figure}[ht!]
    \centering
    \begin{subfigure}[t]{0.48\textwidth}
        \centering
        \includegraphics[scale=.3]{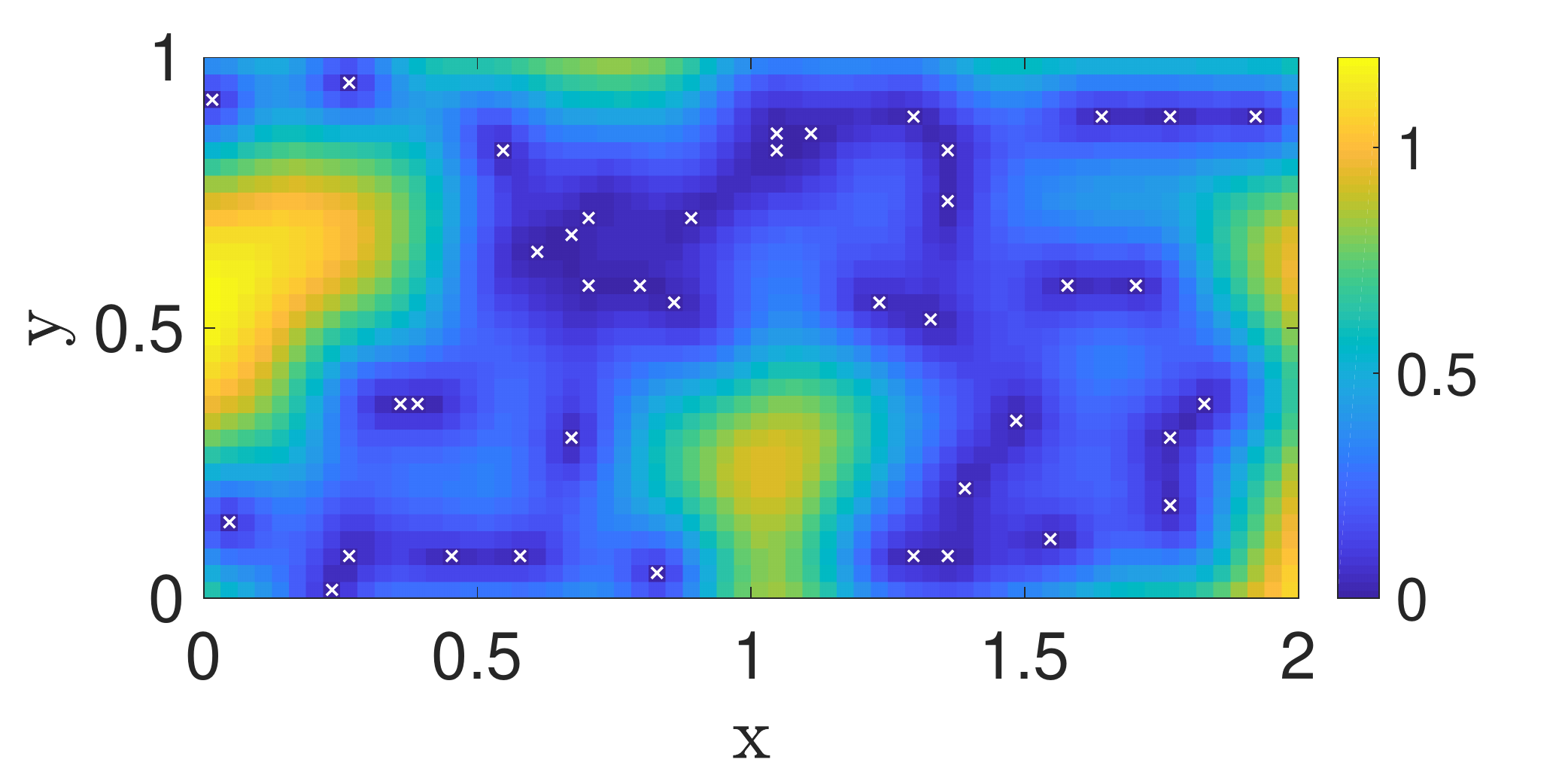}
        \caption{$\sigma_{g^c}$, $N_{am}=0$} \label{RK:fig:C_pg_est_exact_nobs40_mrst_1_03_al}
    \end{subfigure}        \\
    \begin{subfigure}[t]{0.48\textwidth}
        \centering
        \includegraphics[scale=.3]{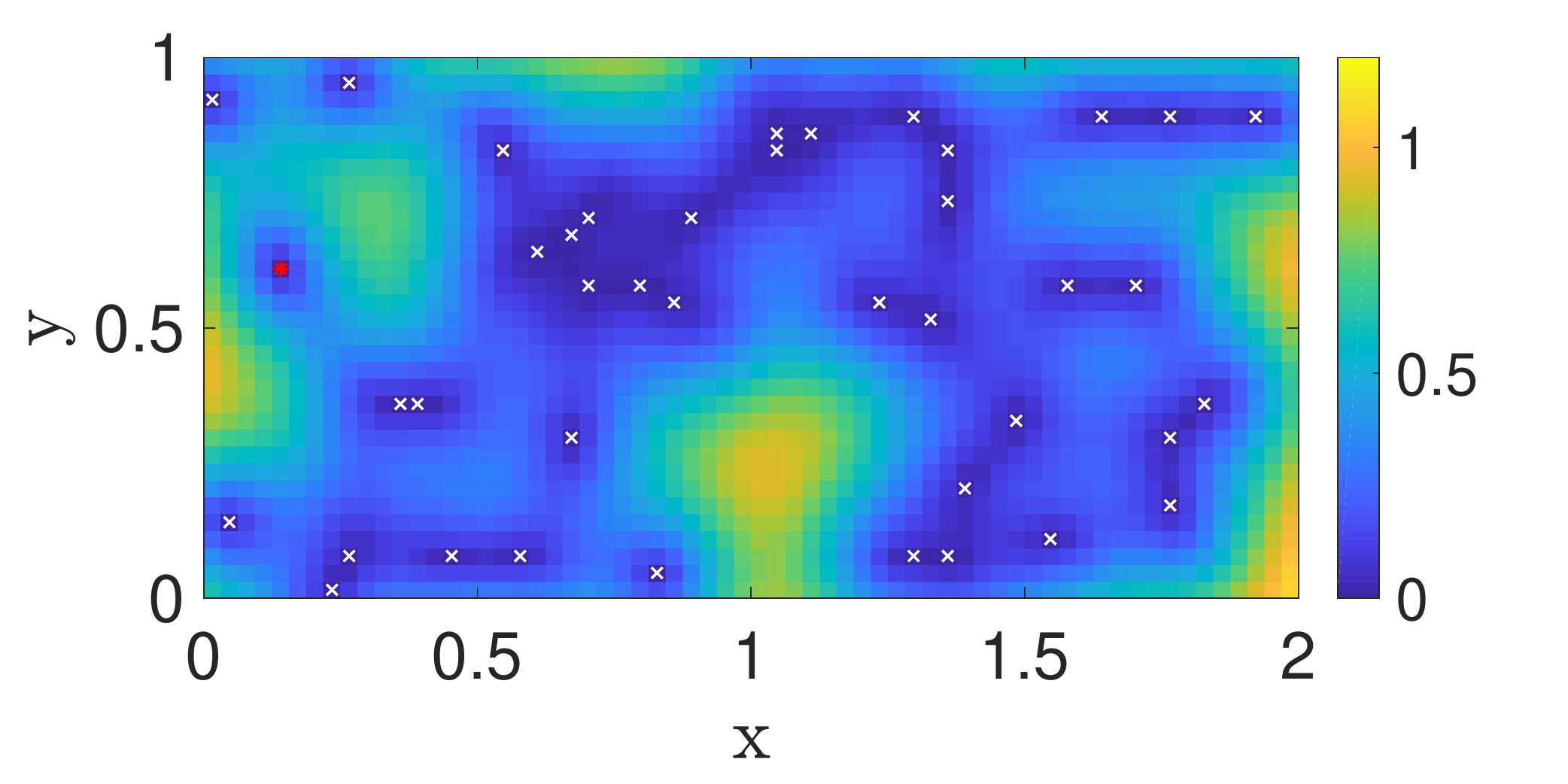}
        \caption{$N_{am}=1$, Method 1} \label{RK:fig:C_pg_exact_nobs40_mrst_x_star_sdev_1_03_1}
    \end{subfigure}    
     \begin{subfigure}[t]{0.48\textwidth}
        \centering
        \includegraphics[scale=.3]{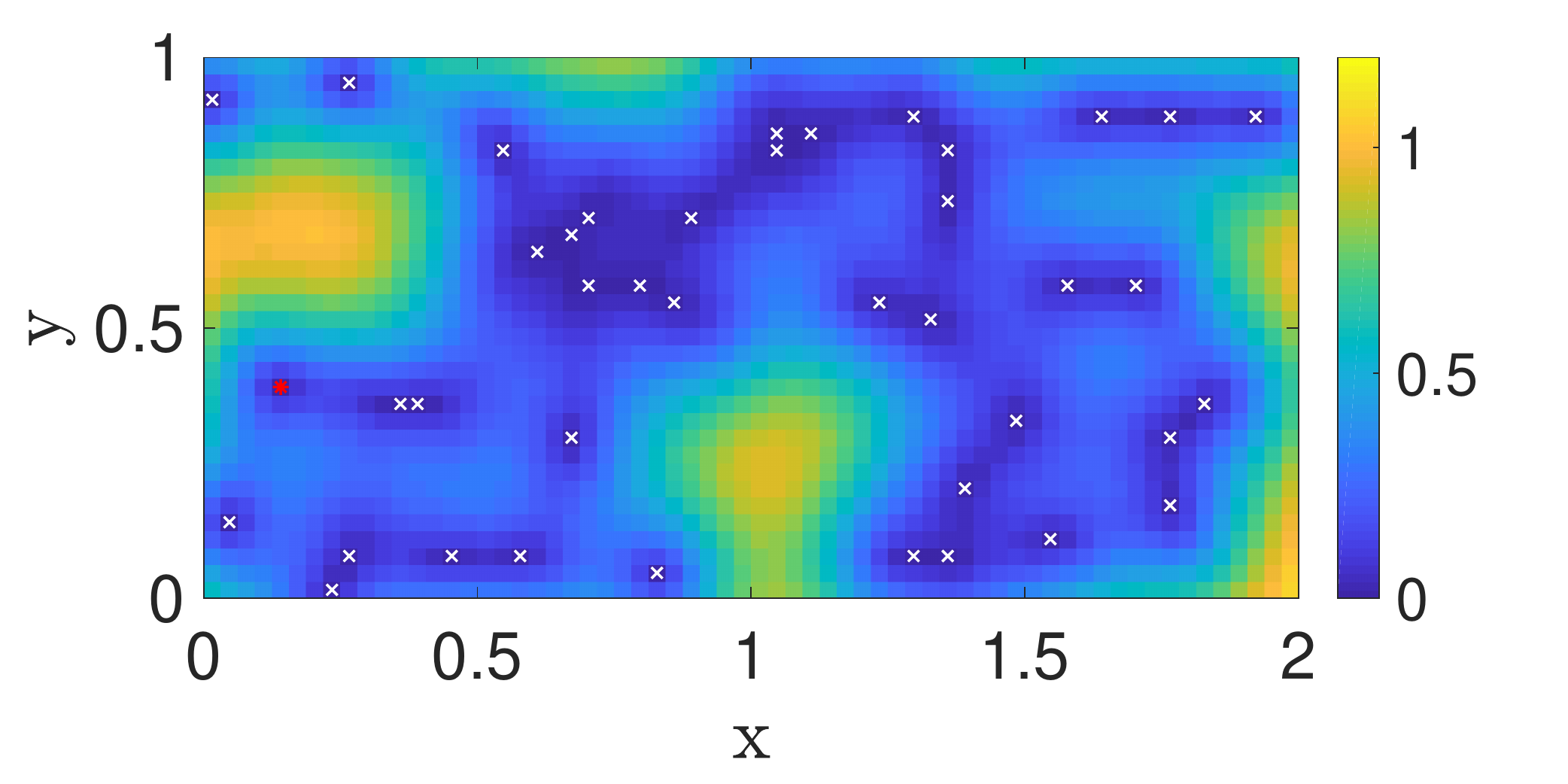}
        \caption{$N_{am}=1$, Method 2} \label{RK:fig:C_pg_exact_nobs40_mrst_x_star_k_u_sdev_1_03_1}
    \end{subfigure}    
    \begin{subfigure}[t]{0.48\textwidth}
        \centering
        \includegraphics[scale=.3]{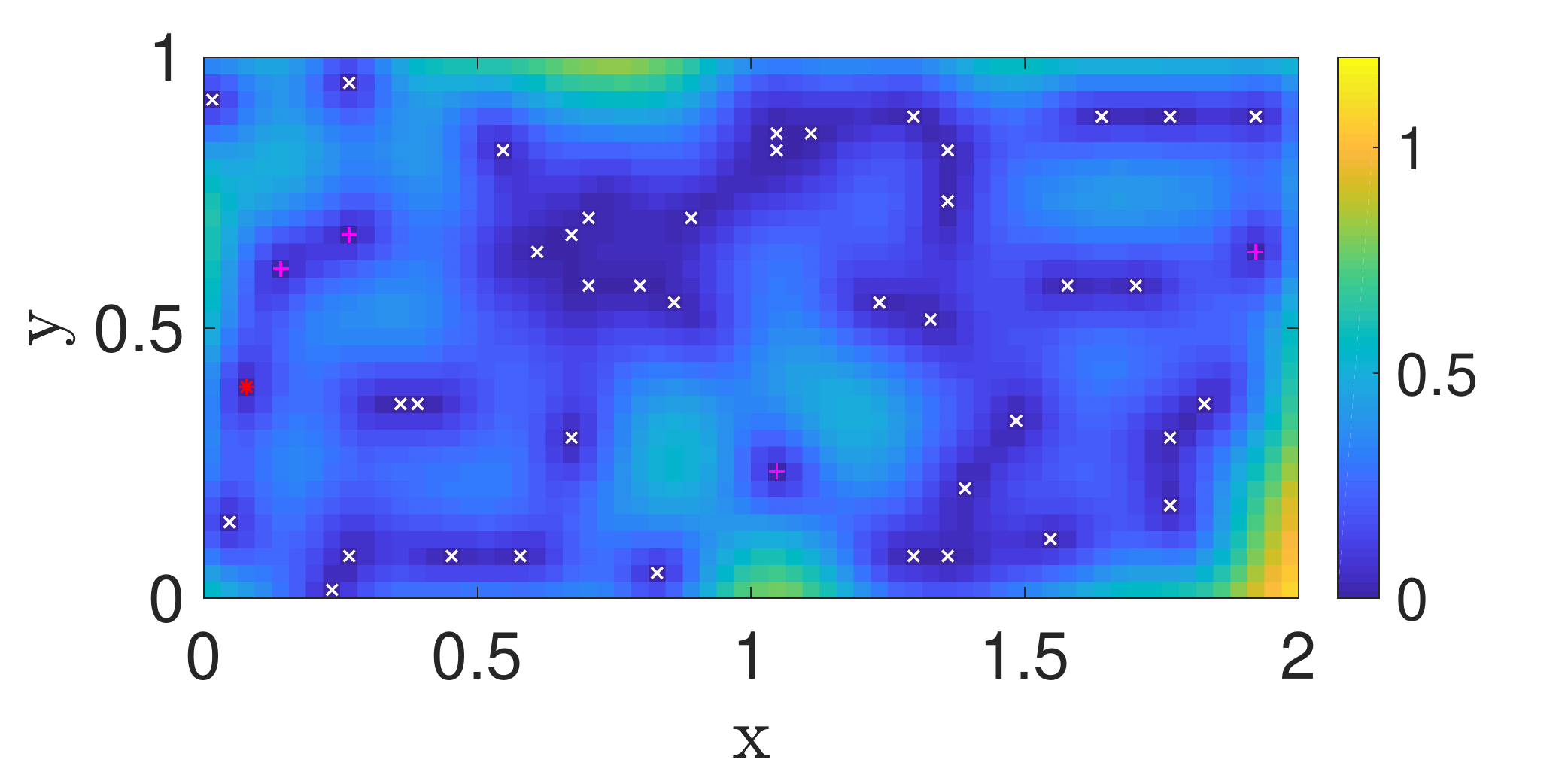}
        \caption{$N_{am}=5$, Method 1} \label{RK:fig:C_pg_exact_nobs40_mrst_x_star_sdev_1_03_5}
    \end{subfigure}        
    \begin{subfigure}[t]{0.48\textwidth}
        \centering
        \includegraphics[scale=.3]{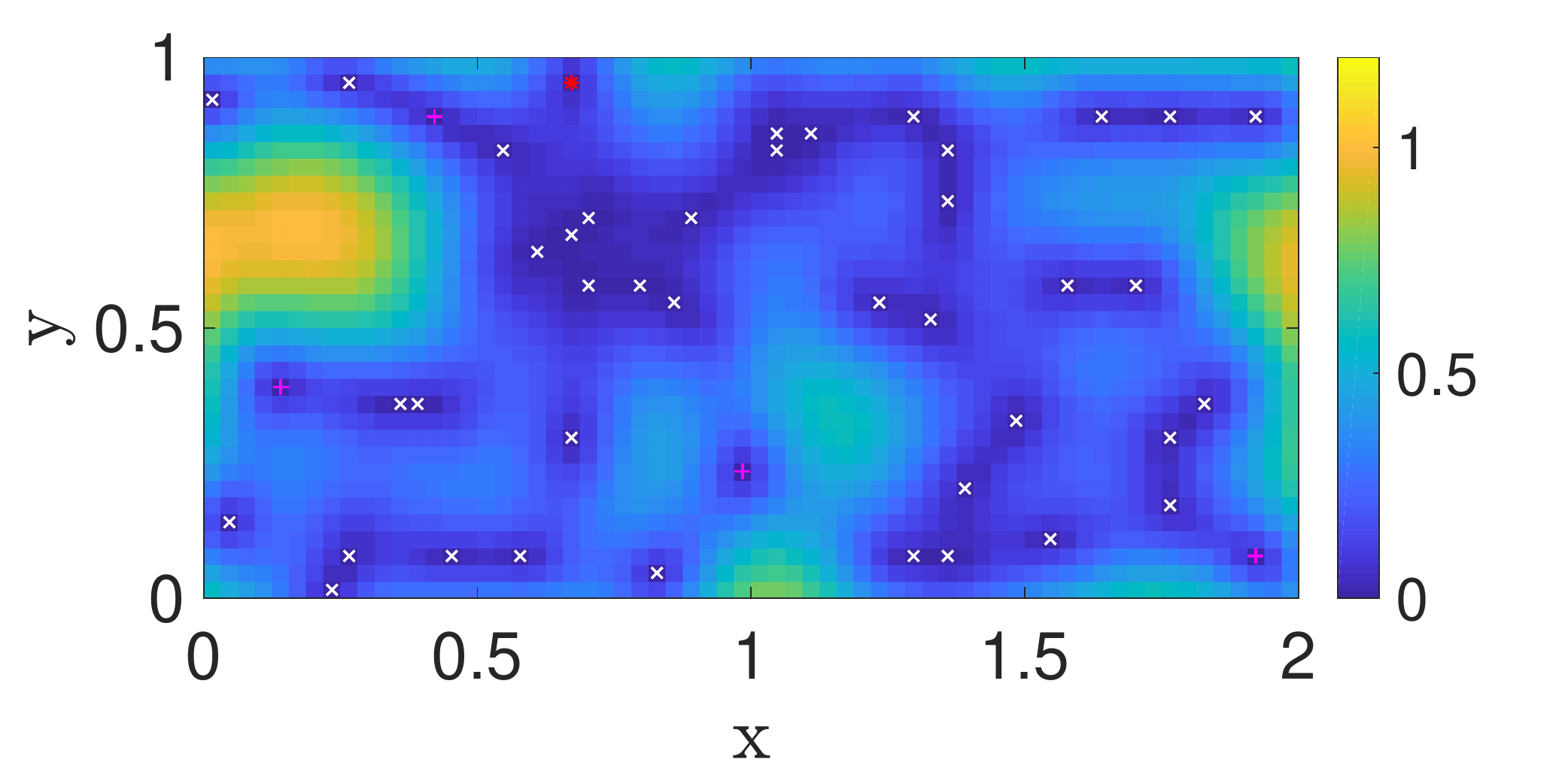}
        \caption{$N_{am}=5$, Method 2} \label{RK:fig:C_pg_exact_nobs40_mrst_x_star_k_u_sdev_1_03_5}
    \end{subfigure}  
    \begin{subfigure}[t]{0.48\textwidth}
        \centering
        \includegraphics[scale=.3]{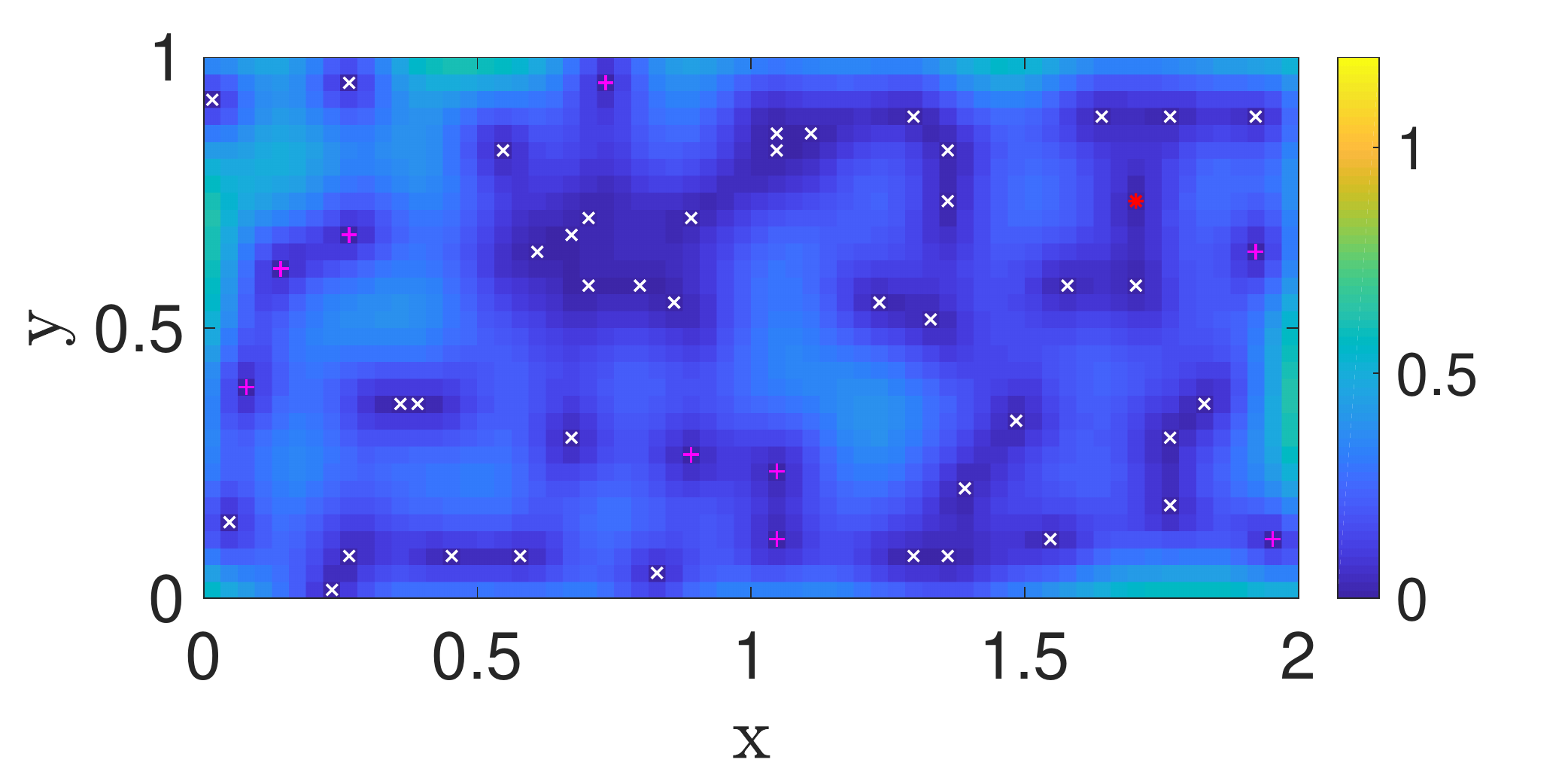}
        \caption{$N_{am}=10$, Method 1} \label{RK:fig:C_pg_exact_nobs40_mrst_x_star_sdev_1_03_10}
    \end{subfigure}    
    \begin{subfigure}[t]{0.48\textwidth}
        \centering
        \includegraphics[scale=.3]{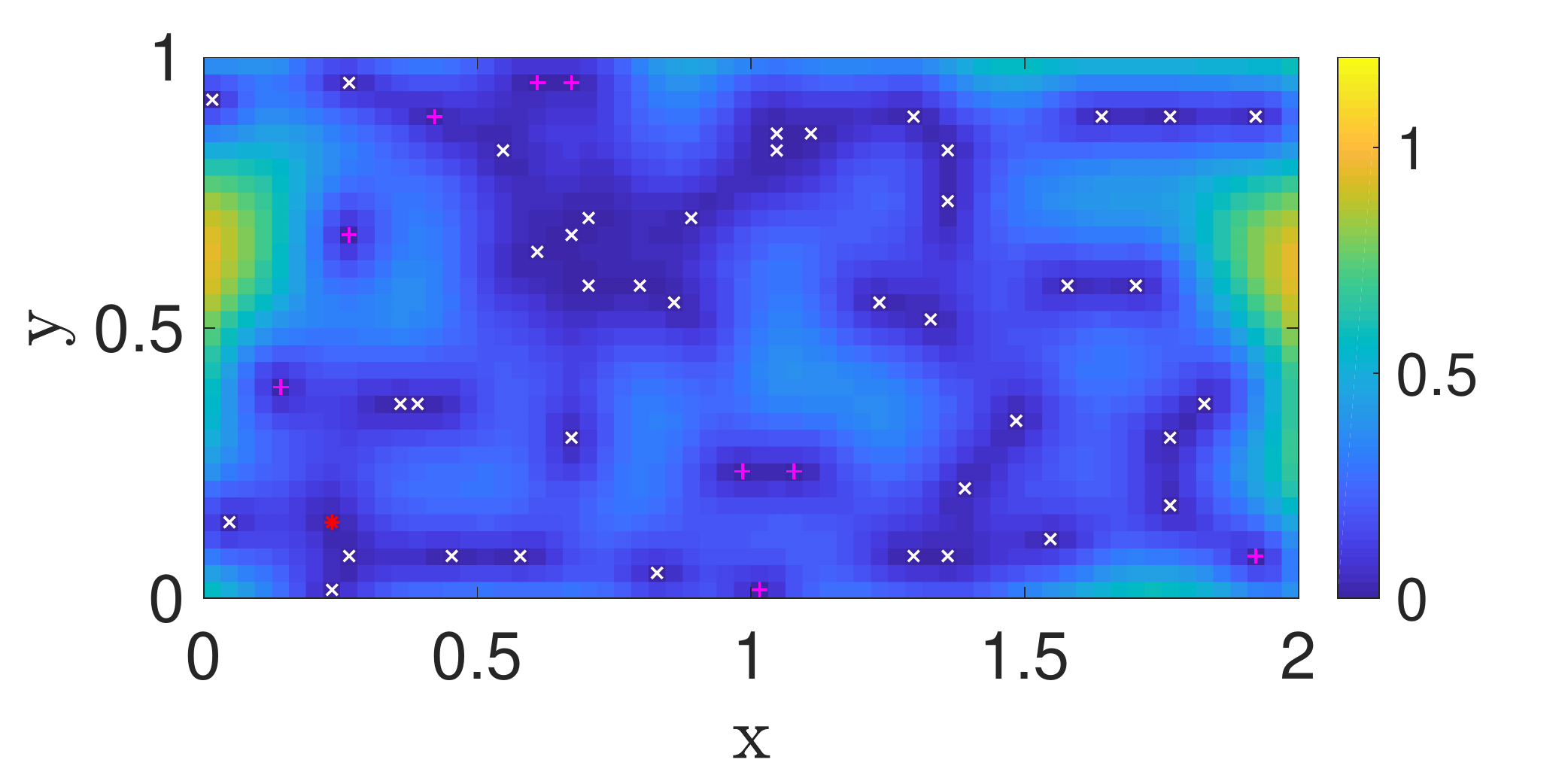}
        \caption{$N_{am}=10$, Method 2} \label{RK:fig:C_pg_exact_nobs40_mrst_x_star_k_u_sdev_1_03_10}
    \end{subfigure}    
    \caption{(a) Standard deviation of $g^c(\mathbf{x},\omega)$ without additional measurements. Standard deviation of $g^c(\mathbf{x},\omega)$ with additional measurements: (b) $N_{am}=1$, Method 1; (c). $N_{am}=1$, Method 2; (d) $N_{am}=1$, Method 5; (e) $N_{am}=5$, Method 2; (f) $N_{am}=10$, Method 1; and (g) $N_{am}=10$, Method 2. Additional observation locations are shown in magenta and the original measurements are shown in red. Standard deviation of unconditional $g(\mathbf{x},\omega)$ is $\sigma_g = 1.63$.}\label{C_pg_exact_nobs40_mrst_x_star_sdev_1_03}
\end{figure}

\begin{figure}[ht!]
    \centering
    \begin{subfigure}[t]{0.48\textwidth}
        \centering
        \includegraphics[scale=.3]{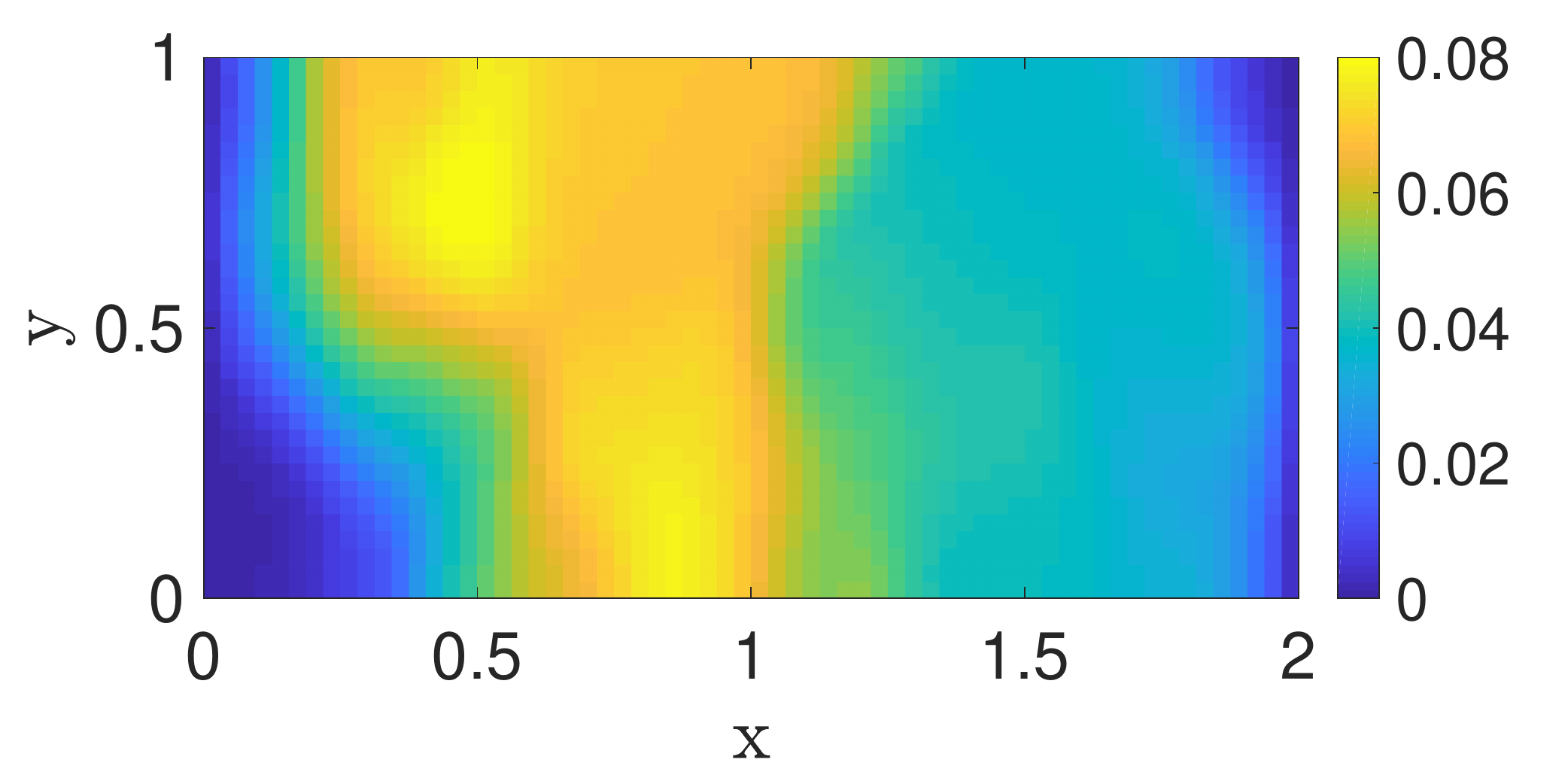}
        \caption{$\sigma_{u^c}$, $N_{am}=0$} \label{RK:fig:u_sdev_exact_cond_mc_mrst_sdev_1_03}
    \end{subfigure}        \\
    \begin{subfigure}[t]{0.48\textwidth}
        \centering
        \includegraphics[scale=.3]{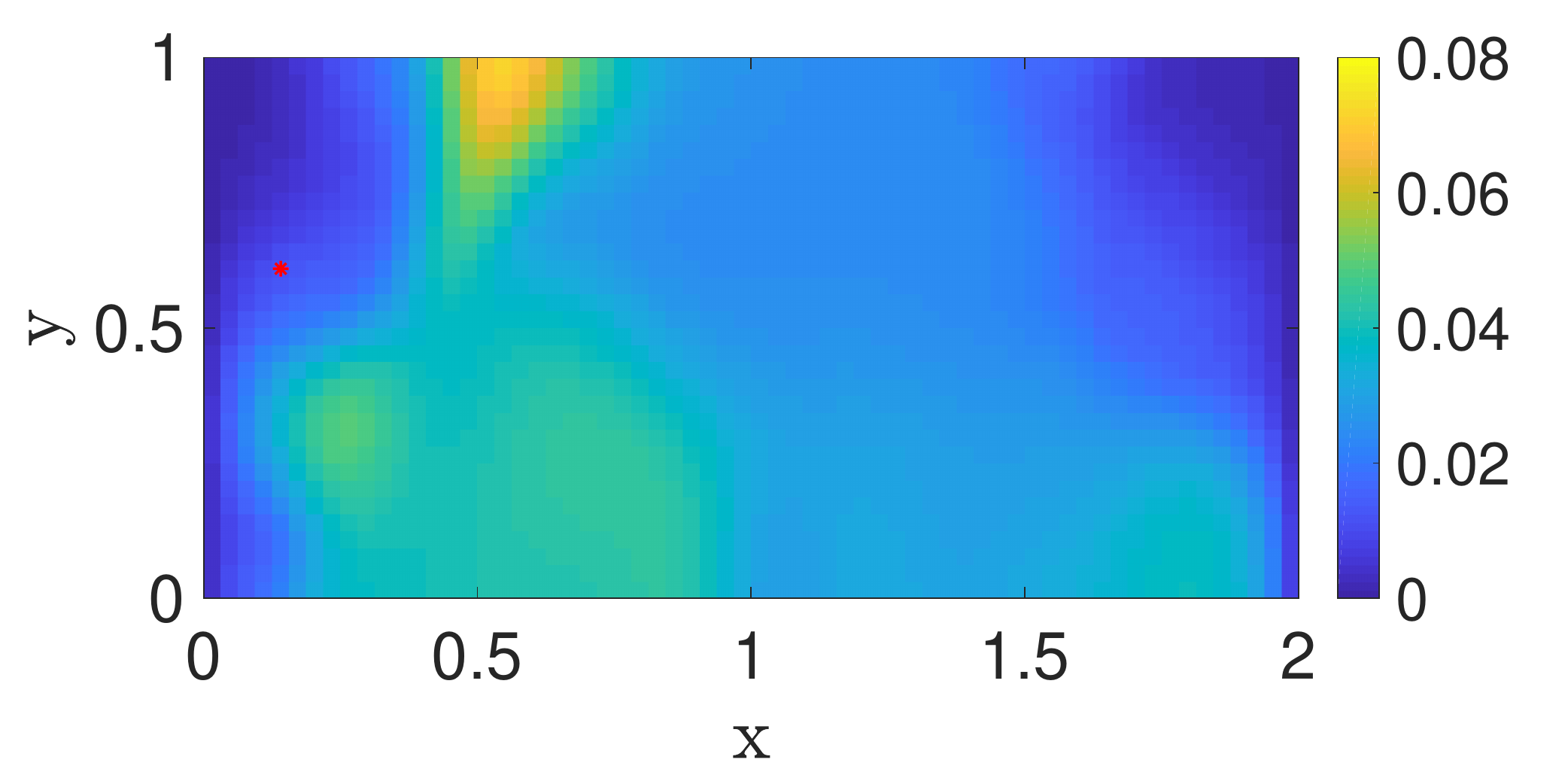}
        \caption{$N_{am}=1$, Method 1} \label{RK:fig:u_sdev_exact_nobs40_mrst_x_star_sdev_1_03_1}
    \end{subfigure}    
     \begin{subfigure}[t]{0.48\textwidth}
        \centering
        \includegraphics[scale=.3]{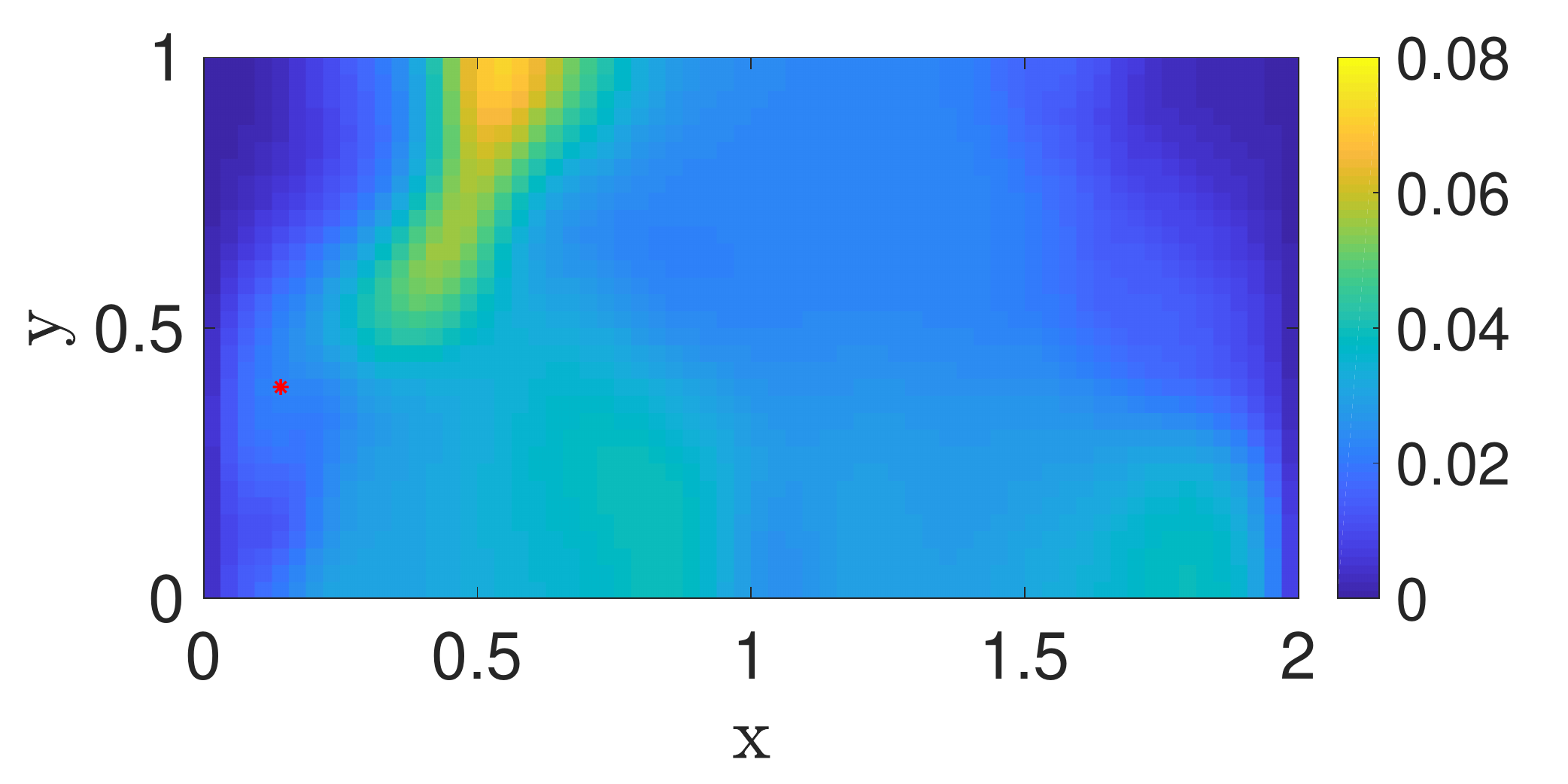}
        \caption{$N_{am}=1$, Method 2} \label{RK:fig:u_sdev_exact_nobs40_mrst_x_star_k_u_sdev_1_03_1}
    \end{subfigure}    
    \begin{subfigure}[t]{0.48\textwidth}
        \centering
        \includegraphics[scale=.3]{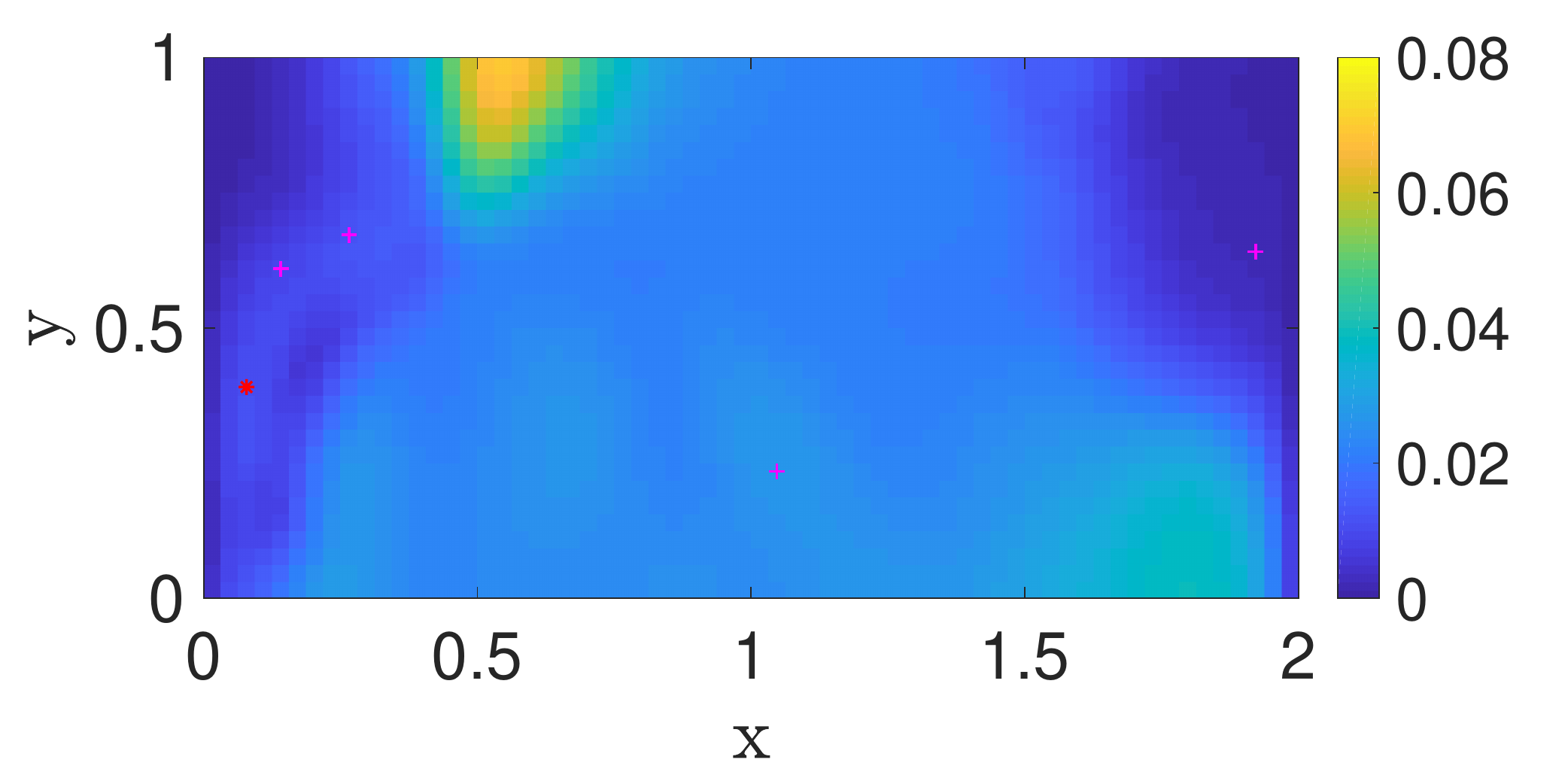}
        \caption{$N_{am}=5$, Method 1} \label{RK:fig:u_sdev_exact_nobs40_mrst_x_star_sdev_1_03_5}
    \end{subfigure}        
    \begin{subfigure}[t]{0.48\textwidth}
        \centering
        \includegraphics[scale=.3]{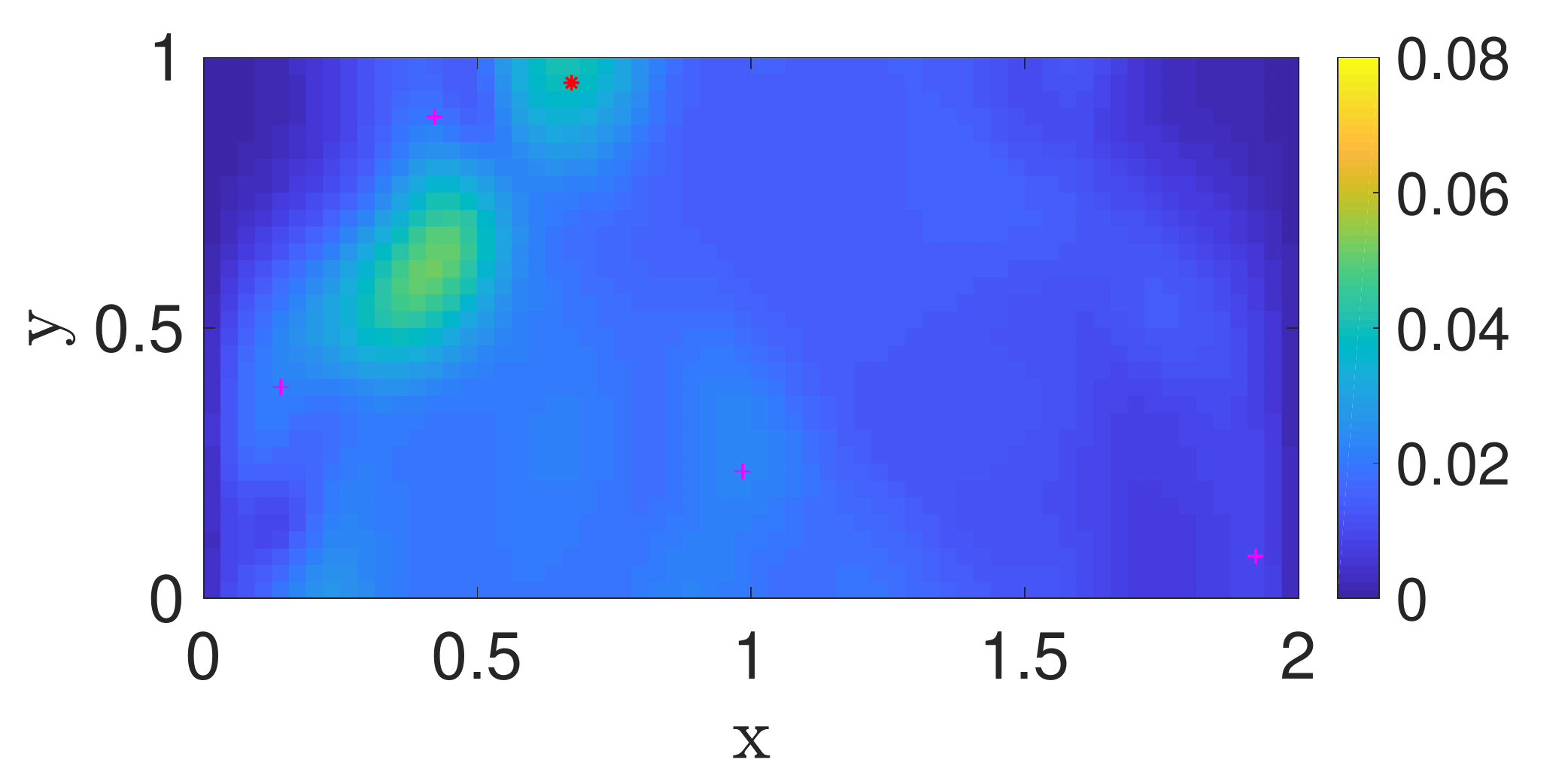}
        \caption{$N_{am}=5$, Method 2} \label{RK:fig:u_sdev_exact_nobs40_mrst_x_star_k_u_sdev_1_03_5}
    \end{subfigure}  
    \begin{subfigure}[t]{0.48\textwidth}
        \centering
        \includegraphics[scale=.3]{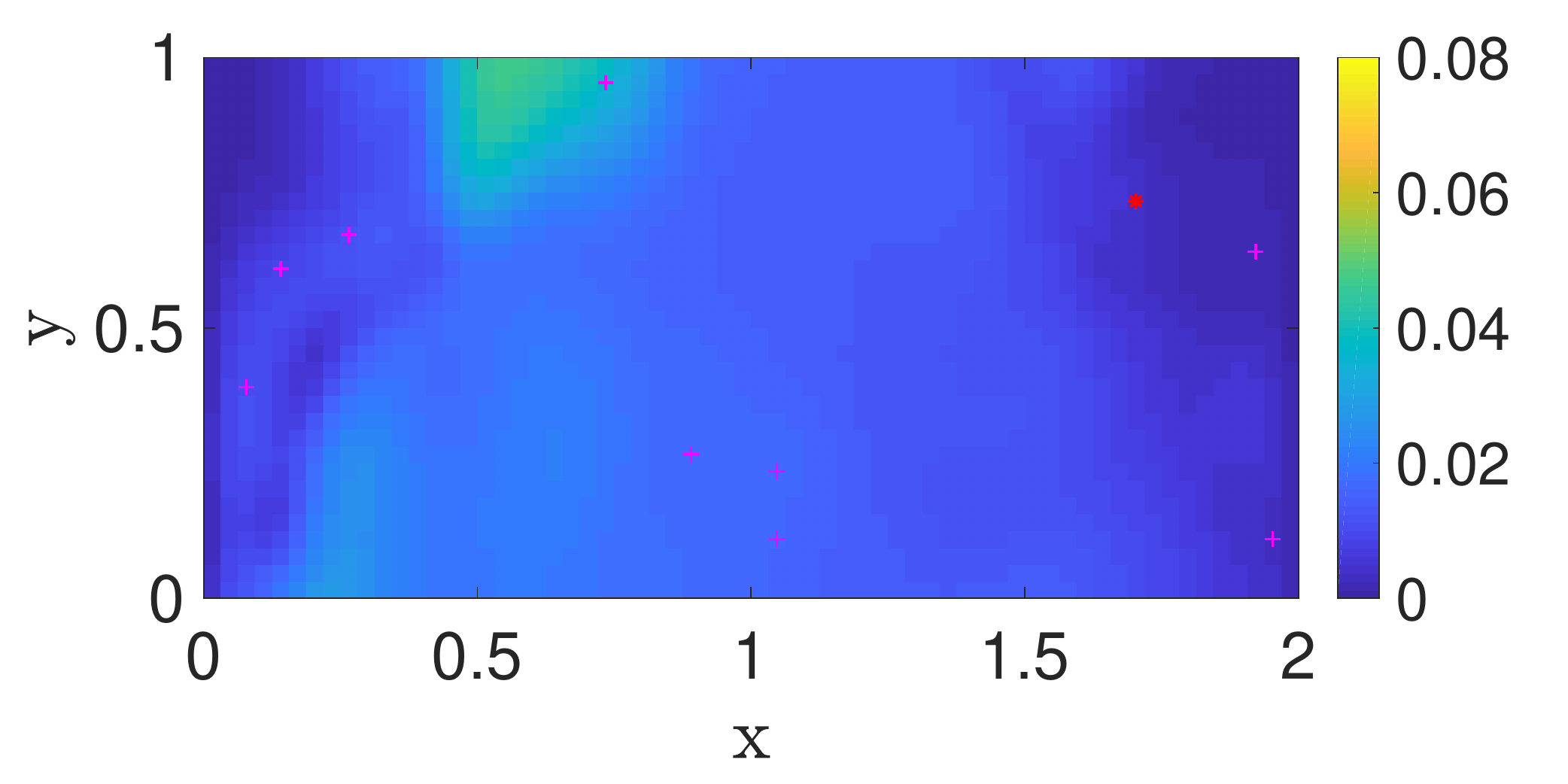}
        \caption{$N_{am}=10$, Method 2} \label{RK:fig:u_sdev_exact_nobs40_mrst_x_star_sdev_1_03_10}
    \end{subfigure}    
    \begin{subfigure}[t]{0.48\textwidth}
        \centering
        \includegraphics[scale=.3]{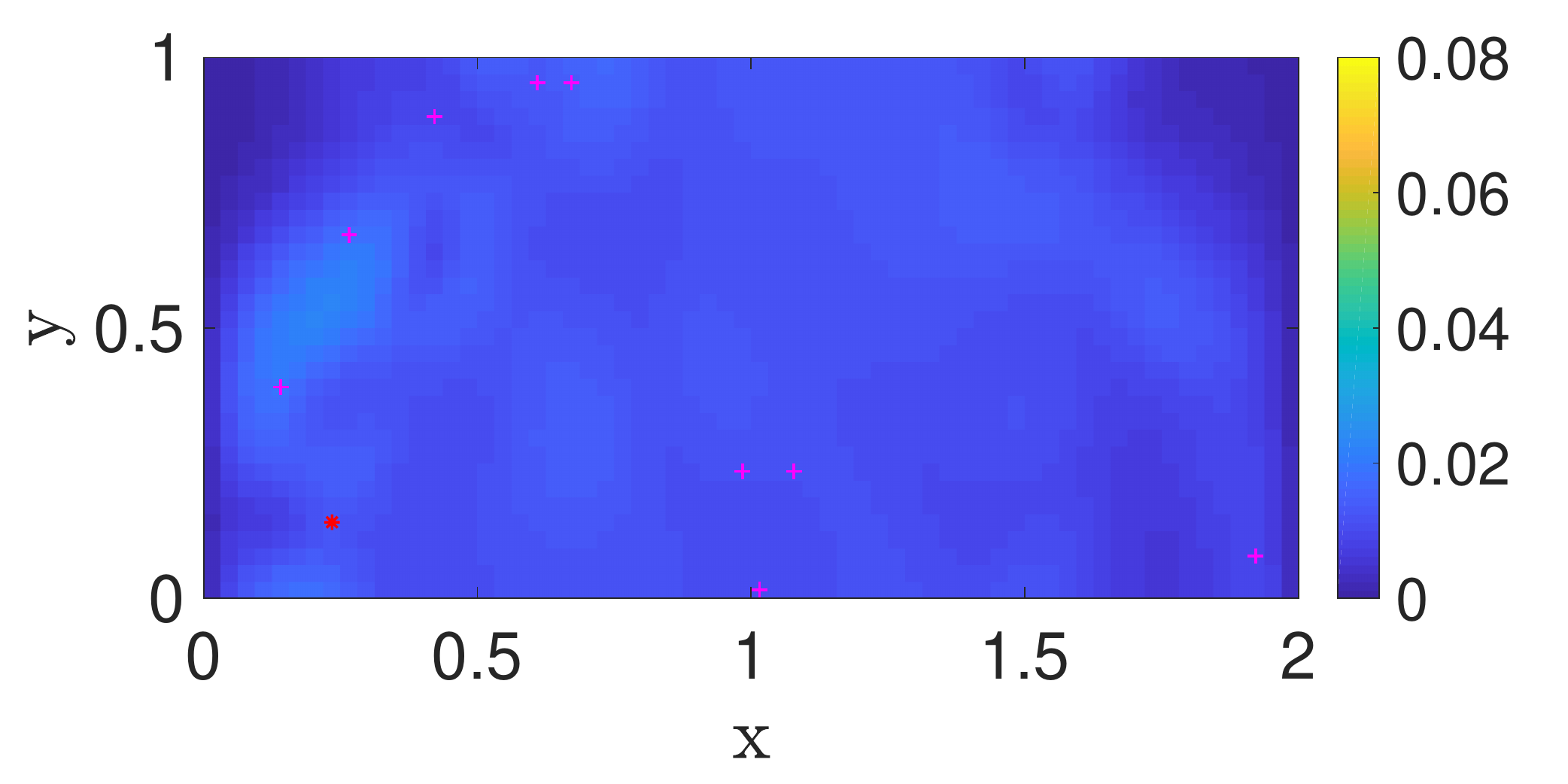}
        \caption{$N_{am}=10$, Method 2} \label{RK:fig:u_sdev_exact_nobs40_mrst_x_star_k_u_sdev_1_03_10}
    \end{subfigure}    
    \caption{(a) Standard deviation of $u^c(\mathbf{x},\omega)$ without additional measurements. Standard deviation of $u^c(\mathbf{x},\omega)$ with additional measurements of $g$: (b) $N_{am}=1$, Method 1; (c). $N_{am}=1$, Method 2; (d) $N_{am}=1$, Method 5; (e) $N_{am}=5$, Method 2; (f) $N_{am}=10$, Method 1; and (g) $N_{am}=10$, Method 2. Additional observation locations of $g$ are shown in magenta and the original measurements are shown in red. Standard deviation of unconditional $g(\mathbf{x},\omega)$ is $\sigma_g = 1.63$.} \label{u_sdev_exact_nobs40_mrst_x_star_k_u_sdev_1_03}
\end{figure}

For the case $\sigma_g=0.65$, Figures~\ref{C_pg_exact_nobs40_mrst_x_star} and \ref{u_sdev_exact_nobs40_mrst_x_star_k_u} show the $\sigma^c_g$ and $\sigma^c_u$ fields for $N_{am} = 1$, $5$, and $10$ obtained using both active learning methods.
Similarly,  $\sigma^c_g$ and $\sigma^c_u$ for $\sigma_g=1.3$ are shown in Figures \ref{C_pg_exact_nobs40_mrst_x_star_sdev_1_03} and \ref{u_sdev_exact_nobs40_mrst_x_star_k_u_sdev_1_03}, respectively. 
As expected, the locations of the additional observations obtained with Methods 1 and 2 are different. Our results show that if the main objective is to predict states, i.e., $u(\mathbf{x})$, rather than coefficients, i.e., $k(\mathbf{x})$, than Method 2 is more efficient than Method 1.  


% plots for sdev 1.3

\section{Conclusions}
\label{sec:conclusions}

We presented two methods for constructing finite-dimensional conditional  Karhunen-Lo\`eve (KL) expansions of partially known  parameters in PDE problems.
We demonstrated that conditioning on data reduces the dimensionality of KL expansions and, most importantly, reduces uncertainty in the solution of PDE problems.
Finally, conditioning on data reduces the computational cost of solving stochastic PDEs. 
We also present a new active learning strategy for acquiring new observations of the data based on minimizing variance of the conditional PDE solution (referred to as Method 2 in the paper) and compared it with the standard active learning method based on minimizing the conditional variance of the partially known parameters (referred to in the paper as Method 1).  

In the first approach for constructing finite-dimensional conditional GPs, presented in \cite{li2014conditional}, the parameter field is conditioned first on data and then discretized by computing its KL expansion.
In the second approach, presented in~\cite{ossiander2014}, the parameter field is discretized first by computing its KL expansion, and then the resulting KL expansion is conditioned on data.
For the second approach, we demonstrated that conditioning leads to dimension reduction of the conditional representation, and we proposed a method for constructing a reduced representation in terms of the effective number of random dimensions.
For a linear diffusion SPDE with uncertain log-normal coefficient, we show that Approach 1 provides a more accurate approximation of the conditional log-normal coefficient and solution of the SPDE than Approach 2 for the same number of random dimensions in a conditional KL expansion. 
Furthermore, Approach 2 provides a good estimate for the number of terms of the truncated KL expansion of the conditional field of Approach 1.

Finally, we demonstrate that the proposed active learning method (Method 2) is more efficient for reducing uncertainty in the solution of the SPDE under consideration (i.e., it leads to a larger reduction of the variance) than the standard active learning method (Method 1).
The difference between two methods increases as the variance of the partially known coefficient increases.

\bibliographystyle{elsarticle-num} 

\bibliography{conditional_jcp}

\begin{thebibliography}{10}
\expandafter\ifx\csname url\endcsname\relax
  \def\url#1{\texttt{#1}}\fi
\expandafter\ifx\csname urlprefix\endcsname\relax\def\urlprefix{URL }\fi
\expandafter\ifx\csname href\endcsname\relax
  \def\href#1#2{#2} \def\path#1{#1}\fi

\bibitem{RK:Ghanem1991}
R.~Ghanem, P.~Spanos, Stochastic Finite Elements: A Spectral Approach,
  Springer-Verlag, 1991.

\bibitem{Babuka2002}
I.~Babu{\v{s}}ka, P.~Chatzipantelidis, On solving elliptic stochastic partial
  differential equations, Computer Methods in Applied Mechanics and Engineering
  191~(37-38) (2002) 4093--4122.

\bibitem{RK:Xiu2002}
D.~Xiu, G.~E. Karniadakis, The wiener--askey polynomial chaos for stochastic
  differential equations, SIAM Journal on Scientific Computing 24 (2002)
  619--644.

\bibitem{Babuka2010}
I.~Babu{\v{s}}ka, F.~Nobile, R.~Tempone, A stochastic collocation method for
  elliptic partial differential equations with random input data, {SIAM} Rev.
  52~(2) (2010) 317--355.

\bibitem{Venturi2013JCP}
D.~Venturi, D.~Tartakovsky, A.~Tartakovsky, G.~Karniadakis, Exact \{PDF\}
  equations and closure approximations for advective-reactive transport,
  Journal of Computational Physics 243~(0) (2013) 323 -- 343.

\bibitem{Lin2010JSC}
G.~Lin, A.~M. Tartakovsky, Numerical studies of three-dimensional stochastic
  darcy's equation and stochastic advection-diffusion-dispersion equation,
  Journal of Scientific Computing 43~(1) (2010) 92--117.

\bibitem{Lin2010JCP}
G.~Lin, A.~Tartakovsky, D.~Tartakovsky, Uncertainty quantification via random
  domain decomposition and probabilistic collocation on sparse grids, Journal
  of Computational Physics 229~(19) (2010) 6995--7012.

\bibitem{Lin2009AWR}
G.~Lin, A.~Tartakovsky, An efficient, high-order probabilistic collocation
  method on sparse grids for three-dimensional flow and solute transport in
  randomly heterogeneous porous media, Advances in Water Resources 32~(5)
  (2009) 712 -- 722.

\bibitem{barajassolano-2016-stochastic}
D.~A. Barajas-Solano, D.~M. Tartakovsky, Stochastic collocation methods for
  nonlinear parabolic equations with random coefficients, SIAM/ASA J. Uncert.
  Quantif. 4 (2016) 475--494.
\newblock \href {http://dx.doi.org/10.1137/130930108}
  {\path{doi:10.1137/130930108}}.

\bibitem{RK:Tipireddy2014}
R.~Tipireddy, R.~Ghanem, Basis adaptation in homogeneous chaos spaces, Journal
  of Computational Physics 259 (2014) 304--317.

\bibitem{TIPIREDDY2017203}
R.~Tipireddy, P.~Stinis, A.~Tartakovsky, Basis adaptation and domain
  decomposition for steady-state partial differential equations with random
  coefficients, Journal of Computational Physics.

\bibitem{tipireddy2017stochastic}
R.~Tipireddy, P.~Stinis, A.~Tartakovsky, Stochastic basis adaptation and
  spatial domain decomposition for pdes with random coefficients, arXiv
  preprint arXiv:1709.02488.

\bibitem{li2016inverse}
W.~Li, G.~Lin, B.~Li, Inverse regression-based uncertainty quantification
  algorithms for high-dimensional models: Theory and practice, Journal of
  Computational Physics 321 (2016) 259--278.

\bibitem{yang2017sliced}
X.~Yang, W.~Li, A.~Tartakovsky, Sliced-inverse-regression-aided rotated
  compressive sensing method for uncertainty quantification, arXiv preprint
  arXiv:1709.07937.

\bibitem{yang2016enhancing}
X.~Yang, H.~Lei, N.~A. Baker, G.~Lin, Enhancing sparsity of hermite polynomial
  expansions by iterative rotations, Journal of Computational Physics 307
  (2016) 94--109.

\bibitem{matheron1963principles}
G.~Matheron, Principles of geostatistics, Economic geology 58~(8) (1963)
  1246--1266.

\bibitem{zhang-2002-stochastic}
D.~Zhang, Stochastic Methods for Flow in Porous Media, Academic Press, 2002.

\bibitem{cressie-2015-geostatistics}
N.~A.~C. Cressie, Geostatistics, John Wiley \& Sons, Inc., 2015, pp. 27--104.
\newblock \href {http://dx.doi.org/10.1002/9781119115151.ch2}
  {\path{doi:10.1002/9781119115151.ch2}}.

\bibitem{neuman1993prediction}
S.~P. Neuman, S.~Orr, Prediction of steady state flow in nonuniform geologic
  media by conditional moments: Exact nonlocal formalism, effective
  conductivities, and weak approximation, Water resources research 29~(2)
  (1993) 341--364.

\bibitem{morales2006non}
E.~Morales-Casique, S.~P. Neuman, A.~Guadagnini, Non-local and localized
  analyses of non-reactive solute transport in bounded randomly heterogeneous
  porous media: Theoretical framework, Advances in water resources 29~(8)
  (2006) 1238--1255.

\bibitem{LU2004859}
Z.~Lu, D.~Zhang,
  \href{http://www.sciencedirect.com/science/article/pii/S0309170804001198}{Conditional
  simulations of flow in randomly heterogeneous porous media using a kl-based
  moment-equation approach}, Advances in Water Resources 27~(9) (2004) 859 --
  874.
\newblock \href
  {http://dx.doi.org/https://doi.org/10.1016/j.advwatres.2004.08.001}
  {\path{doi:https://doi.org/10.1016/j.advwatres.2004.08.001}}.
\newline\urlprefix\url{http://www.sciencedirect.com/science/article/pii/S0309170804001198}

\bibitem{ossiander2014}
M.~E. Ossiander, M.~Peszynska, V.~S. Vasylkivska,
  \href{https://doi.org/10.1155/2014/652594}{Conditional stochastic simulations
  of flow and transport with karhunen-loève expansions, stochastic
  collocation, and sequential gaussian simulation}, J. Appl. Math. 2014 (2014)
  21 pages.
\newblock \href {http://dx.doi.org/10.1155/2014/652594}
  {\path{doi:10.1155/2014/652594}}.
\newline\urlprefix\url{https://doi.org/10.1155/2014/652594}

\bibitem{li2014conditional}
H.~Li, Conditional simulation of flow in heterogeneous porous media with the
  probabilistic collocation method, Communications in Computational Physics
  16~(4) (2014) 1010--1030.

\bibitem{zhao2017active}
L.~Zhao, Z.~Li, B.~Caswell, J.~Ouyang, G.~E. Karniadakis, Active learning of
  constitutive relation from mesoscopic dynamics for macroscopic modeling of
  non-newtonian flows, arXiv preprint arXiv:1709.06228.

\bibitem{raissi2017inferring}
M.~Raissi, P.~Perdikaris, G.~E. Karniadakis, Inferring solutions of
  differential equations using noisy multi-fidelity data, Journal of
  Computational Physics 335 (2017) 736--746.

\bibitem{Ghanem1999}
R.~Ghanem, The nonlinear gaussian spectrum of log-normal stochastic processes
  and variables, Journal of Applied Mechanics 66~(4) (1999) 964.

\bibitem{rasmussen2006gaussian}
C.~E. Rasmussen, C.~K. Williams, Gaussian processes for machine learning,
  Vol.~1, MIT press Cambridge, 2006.

\bibitem{cohn-1996-active}
D.~A. Cohn, Z.~Ghahramani, M.~I. Jordan,
  \href{http://dl.acm.org/citation.cfm?id=1622737.1622744}{Active learning with
  statistical models}, J. Artif. Int. Res. 4~(1) (1996) 129--145.
\newline\urlprefix\url{http://dl.acm.org/citation.cfm?id=1622737.1622744}

\end{thebibliography}

\end{document}